\documentclass[leqno]{amsart}

\newcommand{\klassentyp}{ams} 

\newcommand{\klassenart}{artikel}

\usepackage{ifthen}
\ifthenelse{ 
		 \equal{\klassentyp}{ams}
		}
		{ 
		}
		{ 
		 \usepackage[centertags,sumlimits,intlimits,namelimits]{amsmath} 
		}
\usepackage[fixamsmath,disallowspaces]{mathtools} 
\ifthenelse{
		 \equal{\klassentyp}{ams}
		}
		{ 
		 \usepackage[overload]{empheq} 
		 \usepackage{thmtools}
		}
		{ 
		 \usepackage[overload,ntheorem]{empheq} 
		 \usepackage[amsmath,thmmarks,hyperref]{ntheorem} 
		 \usetagform{default} 
		}
\usepackage{amsxtra}
\usepackage{amssymb} 
\usepackage[mathscr]{eucal} 
\usepackage{upgreek}
\usepackage{upref}
\usepackage{hyperref} 
\usepackage{setspace} 
\ifthenelse{ 
		 \equal{\klassentyp}{ams}
		}
		{ 
		}
		{ 
		 \usepackage{tocloft}
		}
\usepackage[headsepline,footsepline,headtopline,footbotline,ilines,komastyle,nouppercase,autooneside]{scrpage2} 
\makeatletter
\@ifundefined{@chapapp}{\let\@chapapp\chaptername}{}
\makeatother
\usepackage{geometry} 
\ifthenelse{
		 \equal{\klassentyp}{koma}
		}
		{ 
		}
		{
		 \usepackage{typearea}
		 \usepackage{scrextend}
		}

\usepackage{microtype} 

\newcommand{\verfasser}{Clemens Kienzler} 
\newcommand{\titel}{Flat Fronts and Stability for the Porous Medium Equation} 

\ifthenelse{ 
		 \equal{\klassentyp}{ams}
		}
		{ 
		 \title{Flat Fronts and Stability for the Porous Medium Equation} 
		 \author{\verfasser} 
		 \address{Mathematisches Institut, Universit\"at Bonn, Endenicher Allee 60, 53115 Bonn, Germany}
		 \email{kienzler@math.uni-bonn.de} 
		}
		{ 
		 \author{\large \verfasser} 
		 \title{\Large \bfseries \titel} 
		 \date{} 
		 \ifthenelse{ 
				 \equal{\klassentyp}{koma}
				}
				{ 
				 }
				 { 
				 }
		}

\newcommand{\inputencoding}{utf8} 

\newcommand{\sprache}{english} 

\usepackage{mathpazo} 
\ifthenelse{ 
		 \equal{\klassentyp}{koma}
		}
		{ 

		 \KOMAoption{titlepage}{false} 

		 \KOMAoption{abstract}{false} 

		 \KOMAoption{toc}{bib} 

		 \KOMAoption{bibliography}{oldstyle} 

		 \KOMAoption{numbers}{noendperiod} 

		 \KOMAoption{headings}{normal} 
		}
		{ 
		}

\KOMAoption{fontsize}{10pt} 

\KOMAoption{draft}{false} 

\KOMAoption{BCOR}{0mm}

\KOMAoption{paper}{a4} 

\KOMAoption{pagesize}{auto} 

\KOMAoption{twoside}{semi} 

\KOMAoption{headinclude}{false} 
\KOMAoption{headlines}{1.25} 
\KOMAoption{footinclude}{false} 
\KOMAoption{mpinclude}{false} 

\KOMAoption{DIV}{10} 

\newlength{\realtopmargin} 
\hypersetup{colorlinks, 
    citecolor=black,%
    filecolor=black,%
    linkcolor=black,%
    urlcolor=black}

\automark[section]{section} 

\clearscrheadfoot 

\lofoot[]{} 
\cofoot[]{} 
\rofoot[]{} 
\lohead[]{} 
\cohead[]{\titel} 
\rohead[]{\pagemark} 

\lefoot[]{} 
\cefoot[]{} 
\refoot[]{} 
\lehead[]{\pagemark} 
\cehead[]{\verfasser} 
\rehead[]{} 

\setheadwidth[0pt]{text} 
\setfootwidth[0pt]{text} 

\setheadsepline{0.0pt}[] 
\setfootsepline{0.0pt}[] 
\setheadtopline{0.0pt}[] 
\setfootbotline{0.0pt}[] 

\ifthenelse{ 
		 \equal{\klassentyp}{koma}
		}
		{ 
		 \setkomafont{pagehead}{\normalfont\footnotesize\scshape} 
		 \setkomafont{pagefoot}{\normalfont\footnotesize\slshape} 
		 \setkomafont{pagenumber}{\normalfont\footnotesize\slshape} 
		}
		{ 
		}

\ifthenelse{ 
		 \equal{\klassentyp}{koma}
		}
		{ 
		}
		{ 
		}

\vfuzz \hfuzz

\numberwithin{equation}{section}  
\numberwithin{figure}{section}    

\newcommand{\theoremcounterreset}{section} 

\newcounter{Nr}[\theoremcounterreset] 
\renewcommand*{\theNr}{\thesection.\arabic{Nr}} 

\makeatletter
\def\@tocline#1#2#3#4#5#6#7{\relax
  \ifnum #1>\c@tocdepth 
  \else
    \par \addpenalty\@secpenalty\addvspace{#2}%
    \begingroup \hyphenpenalty\@M
    \@ifempty{#4}{%
      \@tempdima\csname r@tocindent\number#1\endcsname\relax
    }{%
      \@tempdima#4\relax
    }%
    \parindent\z@ \leftskip#3\relax \advance\leftskip\@tempdima\relax
    \rightskip\@pnumwidth plus4em \parfillskip-\@pnumwidth
    #5\leavevmode\hskip-\@tempdima
      \ifcase #1
       \or\or \hskip 1em \or \hskip 2em \else \hskip 3em \fi%
      #6\nobreak\relax
    \dotfill\hbox to\@pnumwidth{\@tocpagenum{#7}}\par 
    \nobreak
    \endgroup
  \fi}
\makeatother

\ifthenelse{ 
		 \equal{\klassentyp}{ams}
		}
		{ 
		}
		{ 

		}		

\newcommand{\zitierstil}{alpha} 

\ifthenelse{ 
		 \equal{\zitierstil}{authoryear}
		}
		{ 
		 \newcommand{\authoryearzitierstil}{numbers} 
		 \newcommand{\authoryearbraces}{round} 
		}
		{ 
		}

\newcommand{\bibstil}{my} 

\newcommand{\journalnames}{Journals} 
\newcommand{\authornames}{Authors} 
\newcommand{\database}{Database} 

\DeclareMathOperator{\Id}{id}

\newcommand{\R}{\mathbb{R}}

\newcommand{\N}{\mathbb{N}}

\newcommand{\Rn}{\R ^n}
\newcommand{\Rnm}{\R ^{n-1}}

\newcommand{\Lp}{L ^p}
\newcommand{\Lq}{L ^q}

\newcommand{\Lone}{L ^1}
\newcommand{\Ltwo}{L ^2} 
\newcommand{\Linfty}{L ^\infty} 
\newcommand{\Lploc}{L ^p_{loc}} 
\newcommand{\Loneloc}{L ^1_{loc}} 
\newcommand{\Ltwoloc}{L ^2_{loc}}

\newcommand{\Wonetwoloc}{W ^{1,2}_{loc}}

\newcommand{\stetig}{C} 
\newcommand{\glatt}{\stetig^\infty} 

\newcommand{\distrib}{\mathcal{D}'}

\newcommand{\grad}{\ensuremath{\nabla}}
\newcommand{\gradx}{\grad_{\negthickspace x}}
\newcommand{\grady}{\grad_{\negthickspace y}}

\newcommand{\divx}{\gradx \cdot} 
\newcommand{\divy}{\grady \cdot} 
\newcommand{\lap}{\ensuremath{\Delta}} 
\newcommand{\lapx}{\lap_x} 
\newcommand{\lapy}{\lap_y} 

\newcommand{\Leb}{\mathcal{L}}

\newcommand{\Lebn}{\Leb^n}
\newcommand{\Lebnp}{\Leb^{n+1}}

\newcommand{\menge}[1]
  {\left\{ #1 \right\}}

\newcommand{\norm}[2][] 
  { \ifthenelse{ 
                \equal{#1}{}
               }
               { 
                \ensuremath{\left\lVert#2\right\rVert}
               }
               { 
                \ensuremath{\left\lVert#2\right\rVert_{#1}}
               }
  }

\newcommand{\ip}[3][] 
  { \ifthenelse{ 
                \equal{#1}{}
               }
               { 
                \ensuremath{\boldsymbol{(} \: #2 \: \boldsymbol{|} \: #3 \: \boldsymbol{)}}
               }
               { 
                \ensuremath{\boldsymbol{(} \: #2 \: \boldsymbol{|} \: #3 \: \boldsymbol{)}_{#1}}
               }
  }

\def\Xint#1{\mathchoice
	{\XXint\displaystyle\textstyle{#1}}%
	{\XXint\textstyle\scriptstyle{#1}}%
	{\XXint\scriptstyle\scriptscriptstyle{#1}}%
	{\XXint\scriptscriptstyle\scriptscriptstyle{#1}}%
	\!\int}
\def\XXint#1#2#3{{\setbox0=\hbox{$#1{#2#3}{\int}$}
	\vcenter{\hbox{$#2#3$}}\kern-.5\wd0}}
\def\mint{\Xint-}

\newcommand{\inte}[4][] 
{ 	\ifthenelse{ 
                     \equal{#1}{}
                   }
                   { 
		     \ensuremath{\underset{#2}{\int}#3 \ d#4}
                   }
                   { 
                     \ensuremath{\overset{#1}{\underset{#2}{\int}}#3 \ d#4}
                   }
}	

\newcommand{\minte}[4][] 
{ 	\ifthenelse{ 
                     \equal{#1}{}
                   }
                   { 
		     \ensuremath{\underset{#2}{\mint}#3 \ d#4}
                   }
                   { 
                     \ensuremath{\overset{#1}{\underset{#2}{\mint}}#3 \ d#4}
                   }
}	

\newcommand{\limes}[2] 
  {\ \underset{{\scriptscriptstyle #1 \rightarrow #2}}{\lim \ }}

\newcommand{\limessup}[2] 
  {\underset{{\scriptscriptstyle #1 \rightarrow #2}}{\lim\sup \; }}

\newcommand{\limesinf}[2] 
  {\underset{{\scriptscriptstyle #1 \rightarrow #2}}{\lim\inf \; }}

\newcommand{\goesto}[1][] 
  {  \ifthenelse{ 
                     \equal{#1}{}
                   }
                   { 
		     \longrightarrow 
                   }
                   { 
		     \stackrel{\norm[#1]{\cdot}}{\longrightarrow} 
		   }
   }

\newcommand{\schwachgg}[1][] 
  {  \ifthenelse{ 
                     \equal{#1}{}
                   }
                   { 
		     \rightharpoonup
                   }
                   { 
		     \stackrel{\norm[#1]{\cdot}}{\rightharpoonup} 
		   }
   }

\newcommand{\summe}[2][] 
  {  \ifthenelse{ 
                  \equal{#1}{}
                }
                { 
		  \ensuremath{ \underset{#2}{\sum} }
                } 
                { 
                  \ensuremath{ \overset{#1}{\underset{#2}{\sum}} }
                }
  }

\newcommand{\reihe}[1] 
	{ \summe[\infty]{#1}}

\newcommand{\produkt}[2][] 
  {  \ifthenelse{ 
                  \equal{#1}{}
                }
                { 
		  \ensuremath{ \underset{#2}{\prod} }
                } 
                { 
                  \ensuremath{ \overset{#1}{\underset{#2}{\prod}} }
                }
  }

\newcommand{\vereinig}[2][] 
  {  \ifthenelse{ 
                  \equal{#1}{}
                }
                { 
		  \ensuremath{ \underset{#2}{\bigcup} }
                } 
                { 
                  \ensuremath{ \overset{#1}{\underset{#2}{\bigcup}} }
                }
  }

\newcommand{\schnitt}[2][] 
  {  \ifthenelse{ 
                  \equal{#1}{}
                }
                { 
		  \ensuremath{ \underset{#2}{\bigcap} }
                } 
                { 
                  \ensuremath{ \overset{#1}{\underset{#2}{\bigcap}} }
                }
  }

\newcommand{\ohne}{\smallsetminus}

\newcommand{\supre}[1][] 
  {\underset{{\scriptscriptstyle #1}}{\sup} \;}

\newcommand{\infi}[1][] 
  {\underset{{\scriptscriptstyle #1}}{\inf} \;}

\newcommand{\esssupre}[1][] 
  {\underset{{\scriptscriptstyle #1}}{\textup{ess}\sup} \;}

\newcommand{\maxi}[1][] 
  {\underset{{\scriptscriptstyle #1}}{\max} \!}

\newcommand{\mini}[1][] 
  {\underset{{\scriptscriptstyle #1}}{\min} \;}

\newcommand{\dt}{\ensuremath{\partial_t}}

\newcommand{\ds}{\ensuremath{\partial_s}}
\newcommand{\dshut}{\ensuremath{\partial_{\hat{s}}}}

\newcommand{\dx}{\ensuremath{\partial_x}}

\newcommand{\dy}{\ensuremath{\partial_y}}
\newcommand{\dystrich}{\ensuremath{\partial_{y'}}}
\newcommand{\dyhut}{\ensuremath{\partial_{\hat{y}}}}
\newcommand{\dz}{\ensuremath{\partial_z}}

\newcommand{\dxn}{\ensuremath{\partial_{x_n}}}

\newcommand{\dyj}{\ensuremath{\partial_{y_j}}} 
 
\newcommand{\dyn}{\ensuremath{\partial_{y_n}}}

\newcommand{\ldef}{\vcentcolon=} 
\newcommand{\rdef}{=\vcentcolon} 

\DeclareMathOperator{\supp}{supp} 

\newcommand{\Yon}[1]{\norm[Y^{\small{on}}_\theta(p)]{#1}}
\newcommand{\Yoff}[1]{\norm[Y^{\small{off}}_{\varepsilon_1,\, \varepsilon_2}(p)]{#1}}
\newcommand{\Yoffone}[1]{\norm[Y^{\small{off}}_{\varepsilon_1,\, 1}(p)]{#1}}
\newcommand{\Yoffoneth}[1]{\norm[Y^{\small{off}}_{3 - \theta, \, 1}(p)]{#1}}
\newcommand{\Pos}{\mathcal{P}}

\newcommand{\Interf}{\mathcal{I}}
\newcommand{\Beu}{B^{eu}}

\newcommand{\Iquer}{{\overline{I}}}

\newcommand{\Ilinks}{{\leftoverbracket{I}}} 
\newcommand{\Irechts}{{\rightoverbracket{I}}} 
 
\newcommand{\Hquer}{{\overline{H}}}
\newcommand{\en}{{\vec{e}_n}}
\newcommand{\siga}{\sigma_{|\alpha|}}
\newcommand{\musig}{\mu_\sigma}
\newcommand{\musigp}{\mu_{1 + \sigma}}

\newcommand{\Lpsig}{\Lp_\sigma}
\newcommand{\Lonesig}{\Lone_\sigma}

\newcommand{\Ltwosig}{\Ltwo_\sigma}
\newcommand{\Ltwosigp}{\Ltwo_{1 + \sigma}}

\newcommand{\yzeron}{y_{0,n}}
\newcommand{\zzeron}{z_{0,n}}
\newcommand{\tquer}{\overline{t}}
\newcommand{\squer}{\overline{s}}
\newcommand{\taquer}{\overline{\tau}}
\newcommand{\xquer}{\overline{x}}
\newcommand{\yquer}{\overline{y}}
\newcommand{\zquer}{\overline{z}}
\newcommand{\Lpsiga}{\Lp_{\siga}}
\newcommand{\Wmp}{W^{m,p}}

\newcommand{\Wmpsigvec}{\Wmp_{\vec{\sigma}}}

\newcommand{\Pnoned}{(\: \cdot \:)}
\newcommand{\Pnonedhut}{(\: \hat{\cdot} \:)}
\newcommand{\Pn}{\Pnoned_n}
\newcommand{\Pnhut}{\Pnonedhut_n}

\newcommand{\dist}{d}
\newcommand{\quasi}{\tilde{\dist}}
\newcommand{\Dist}{D}
\newcommand{\homLip}{\dot{C}^{0,1}}
\newcommand{\constd}{c_d}
\newcommand{\constl}{c_L}
\newcommand{\constpsi}{\zeta_\Psi}
\newcommand{\epsi}{\varepsilon_\Psi}

\newcommand{\Pnl}{\Pn^l}
\newcommand{\Pnlhut}{\Pnhut^l}

\newcommand{\dta}{\ensuremath{\partial_\tau}}

\newcommand{\ipr}[3][] 
  { \ifthenelse{ 
                \equal{#1}{}
               }
               { 
                \ensuremath{\left\langle \: #2 \: \left| \: #3 \right. \: \right\rangle}
               }
               { 
                \ensuremath{\left\langle \: #2 \: \left| \: #3 \right. \: \right\rangle_{#1}}
               }
  }

\makeatletter
\def\leftoverbracketfill{$\m@th\makesm@sh{\llap{\vrule\@height.3\p@\@depth1.5\p@\@width.3\p@}}
\leaders\vrule\@height.3\p@\hskip 0.6ex 
\makesm@sh{\rlap{\vrule\@height.0\p@\@depth0\p@\@width.0\p@}}$} 

\def\rightoverbracketfill{$\m@th\makesm@sh{\llap{\vrule\@height.0\p@\@depth0.0\p@\@width.0\p@}}
\hskip 1ex \leaders\vrule\@height.3\p@\hskip 0.6ex
\makesm@sh{\rlap{\vrule\@height.3\p@\@depth1.5\p@\@width.3\p@}}$} 

\def\leftoverbracket#1{\mathop{\vbox{\ialign{##\cr\noalign{\kern1.3\p@}
\leftoverbracketfill\cr\noalign{\kern1.5\p@\nointerlineskip}
$\hfil\displaystyle{#1}\hfil$\cr}}}\limits}

\def\rightoverbracket#1{\mathop{\vbox{\ialign{##\cr
\rightoverbracketfill\cr\noalign{\kern1.5\p@\nointerlineskip}
$\hfil\displaystyle{#1}\hfil$\cr}}}\limits}

\makeatother

\ifthenelse{ 
                \equal{\sprache}{english}
               }
               { 
                \newcommand{\theore}{Theorem}
                \newcommand{\koro}{Corollary}
                \newcommand{\beme}{Remark}
                \newcommand{\vermu}{Conjecture}
                \newcommand{\beisp}{Example}
                \newcommand{\proo}{Proof}
                \newcommand{\lemm}{Lemma}
                \newcommand{\defini}{Definition}
                \newcommand{\propos}{Proposition}
               }
               { 
                \newcommand{\theore}{Satz}
                \newcommand{\koro}{Korollar}
                \newcommand{\beme}{Bemerkung}
                \newcommand{\vermu}{Vermutung}
                \newcommand{\beisp}{Beispiel}
                \newcommand{\proo}{Beweis}
                \newcommand{\lemm}{Lemma}
                \newcommand{\defini}{Definition}
                \newcommand{\propos}{Proposition}
               }

\makeatletter
\ifcsname phantomsection\endcsname
    \newcommand*{\qrr@gobblenexttocentry}[5]{}
\else
    \newcommand*{\qrr@gobblenexttocentry}[4]{}
\fi
\newcommand*{\addsubsection}{%
    \addtocontents{toc}{\protect\qrr@gobblenexttocentry}%
    \subsection}
\makeatother

\ifthenelse{ 
		 \equal{\klassentyp}{ams}
		}
		{ 

			\newenvironment*{theo}[1][] 
			  { 
			   \begin{sloppypar}\bigskip\refstepcounter{Nr}
			   \ifthenelse{ 
			                \equal{#1}{}
			               }
			               { 
			                \textbf{\theNr{}~\theore} 
			               }
			               { 
			                \textbf{\theNr{}~\theore \ (#1)} 
			               }
			   \begin{it} 
			  }
			  { 
			   \end{it}
			   \end{sloppypar}\bigskip
			  }
			
			\newenvironment*{prop}[1][] 
			  { 
			   \begin{sloppypar}\bigskip\refstepcounter{Nr}
			   \ifthenelse{ 
			                \equal{#1}{}
			               }
			               { 
			                \textbf{\theNr{}~\propos} 
			               }
			               { 
			                \textbf{\theNr{}~\propos \ (#1)} 
			               }
			   \begin{it} 
			  }
			  { 
			   \end{it}
			   \end{sloppypar}\bigskip
			  }
			
			\newenvironment*{kor}[1][] 
			  { 
			   \begin{sloppypar}\bigskip\refstepcounter{Nr}
			   \ifthenelse{ 
			                \equal{#1}{}
			               }
			               { 
			                \textbf{\theNr{}~\koro} 
			               }
			               { 
			                \textbf{\theNr{}~\koro \ (#1)} 
			               }
			   \begin{it} 
			  }
			  { 
			   \end{it}
			   \end{sloppypar}\bigskip
			  }
			
			\newenvironment*{lemma}[1][] 
			  { 
			   \begin{sloppypar}\bigskip\refstepcounter{Nr}
			   \ifthenelse{ 
			                \equal{#1}{}
			               }
			               { 
			                \textbf{\theNr{}~\lemm} 
			               }
			               { 
			                \textbf{\theNr{}~\lemm \ (#1)} 
			               }
			   \begin{it} 
			  }
			  { 
			   \end{it}
			   \end{sloppypar}\bigskip
			  }
			
			\newenvironment*{bem}[1][] 
			  { 
			   \begin{sloppypar}\bigskip\refstepcounter{Nr}
			   \ifthenelse{ 
			                \equal{#1}{}
			               }
			               { 
			                \textbf{\theNr{}~\beme} 
			               }
			               { 
			                \textbf{\theNr{}~\beme \ (#1)} 
			               }
			   \begin{it} 
			  }
			  { 
			   \end{it}
			   \end{sloppypar}\bigskip
			  }
			
			\newenvironment*{verm}[1][] 
			  { 
			   \begin{sloppypar}\bigskip\refstepcounter{Nr}
			   \ifthenelse{ 
			                \equal{#1}{}
			               }
			               { 
			                \textbf{\theNr{}~\vermu} 
			               }
			               { 
			                \textbf{\theNr{}~\vermu \ (#1)} 
			               }
			   \begin{it} 
			  }
			  { 
			   \end{it}
			   \end{sloppypar}\bigskip
			  }
			
			\newenvironment*{bsp}[1][] 
			  { 
			   \begin{sloppypar}\bigskip\refstepcounter{Nr}
			   \ifthenelse{ 
			                \equal{#1}{}
			               }
			               { 
			                \textbf{\theNr{}~\beisp} 
			               }
			               { 
			                \textbf{\theNr{}~\beisp \ (#1)} 
			               }
			   \begin{it} 
			  }
			  { 
			   \end{it}
			   \end{sloppypar}\bigskip
			  }
			
			\newenvironment*{defin}[1][] 
			  { 
			   \begin{sloppypar}\bigskip\refstepcounter{Nr}
			   \ifthenelse{ 
			                \equal{#1}{}
			               }
			               { 
			                \textbf{\theNr{}~\defini} 
			               }
			               { 
			                \textbf{\theNr{}~\defini \ (#1)} 
			               }
			   \begin{it} 
			  }
			  { 
			   \end{it}
			   \end{sloppypar}\bigskip
			  }
			
			\newenvironment*{bew}[1][] 
			{ 
			   \begin{sloppypar}
			   \ifthenelse{ 
			                \equal{#1}{}
			               }
			               { 
					\textbf{\proo: }
			               }
			               { 
					\textbf{\proo \text{ of }#1: }
			               }
			  }
			  { 
			  $ \quad~\rule{1ex}{1ex}$ \end{sloppypar}\bigskip
			  } 

		}
		{ 

			\theoremstyle{change} 
			\theorembodyfont{\itshape} 
			\theoremheaderfont{\normalfont\bfseries} 
			\theoremseparator{ } 
			\theoremnumbering{arabic} 
			\theorempreskip{\topsep} 
			\theorempostskip{\topsep} 
			\theoremindent0.0cm 
			\theoremsymbol{} 
			\theoremprework{\smallskip}
			\theorempostwork{\smallskip}
			
			\newtheorem{theo}{\theore}[\theoremcounterreset] 
			
			\theoremprework{\smallskip}
			\theorempostwork{\smallskip}
			\newtheorem{prop}[theo]{\propos} 
			
			\theoremprework{\smallskip}
			\theorempostwork{\smallskip}
			\newtheorem{kor}[theo]{\koro} 
			
			\theoremprework{\smallskip}
			\theorempostwork{\smallskip}
			\ifthenelse{ 
			              	 \equal{\klassentyp}{beamer}
					}
			               	{ 
					 \renewtheorem{lemma}[theo]{\lemm} 
					}
					{ 
					 \newtheorem{lemma}[theo]{\lemm} 
					}
			
			\theoremprework{\smallskip}
			\theorempostwork{\smallskip}
			\newtheorem{bem}[theo]{\beme} 
			
			\theoremprework{\smallskip}
			\theorempostwork{\smallskip}

			\theoremprework{\smallskip}
			\theorempostwork{\smallskip}

			\theoremprework{\smallskip}
			\theorempostwork{\smallskip}
			\newtheorem{defin}[theo]{\defini} 
			
			\theoremstyle{nonumberplain} 
			\theorembodyfont{\normalfont} 
			\theoremheaderfont{\bfseries}
			\theoremseparator{:}
			\theoremsymbol{\rule{1ex}{1ex}} 
			\theorempreskip{0.1 \topsep} 
			\theorempostwork{\smallskip}
			\newtheorem{bew}{\proo} 

			\renewcommand*{\addsubsection}{\subsection} 

		}

\usepackage[T1]{fontenc}
\usepackage[\inputencoding]{inputenc} 
\usepackage[\sprache]{babel} 
\usepackage{csquotes} 
\usepackage{enumitem} 
\usepackage{graphicx} 
\usepackage{eurosym}
\ifthenelse{
		 \equal{\zitierstil}{authoryear}
		}
		{ 
		 \usepackage[\authoryearzitierstil,\authoryearbraces]{natbib}
		}
		{
		 \newcommand{\citet}{\cite}
		 \newcommand{\citep}{(\cite}
		}

\begin{document}


\ifthenelse{ 
		 \equal{\klassenart}{buch}
		}
		{ 
		 \frontmatter 
		}
		{ 
		}


\begin{spacing}{1} 

\ifthenelse{ 
		 \equal{\klassenart}{buch}
		}
		{ 
		 \maketitle
		}
		{ 
		 \ifthenelse{ 
		 		 \equal{\klassentyp}{ams}
				}
				{ 
		 		 \begin{abstract}
				 \noindent This work is concerned with the equation $ \partial_t \rho =  \Delta_x \rho^m  $, $ m > 1 $, known as the porous medium equation. It shows stability of the pressure of solutions close to flat travelling wave fronts in the homogeneous Lipschitz sense that is in a way optimal for the treatment of the equation. This is the first result of this type and implies global regularity estimates for any number of derivatives of the pressure. Consequences include smoothness, analyticity in temporal and tangential directions, and analyticity of the interface between empty and occupied regions. In the course of the argument a Gaussian estimate in an intrinsically arising space of homogeneous type is crucial to obtain linear estimates by means of the non-Euclidean Calder\'{o}n-Zygmund singular integral theory.

		      		 \end{abstract}
		 		 \maketitle 
				}
				{ 
		 		 \maketitle 
		 		 \begin{abstract}
				 { \small 
				 \noindent This work is concerned with the equation $ \partial_t \rho =  \Delta_x \rho^m  $, $ m > 1 $, known as the porous medium equation. It shows stability of the pressure of solutions close to flat travelling wave fronts in the homogeneous Lipschitz sense that is in a way optimal for the treatment of the equation. This is the first result of this type and implies global regularity estimates for any number of derivatives of the pressure. Consequences include smoothness, analyticity in temporal and tangential directions, and analyticity of the interface between empty and occupied regions. In the course of the argument a Gaussian estimate in an intrinsically arising space of homogeneous type is crucial to obtain linear estimates by means of the non-Euclidean Calder\'{o}n-Zygmund singular integral theory.
 }
		 		 \end{abstract}
				}
		}


\tableofcontents
\thispagestyle{empty}

\end{spacing} 


\ifthenelse{ 
		 \equal{\klassenart}{buch}
		}
		{ 
		 \mainmatter 
		}
		{ 
		}

\section{The Porous Medium Equation} \label{section_pme}
\subsection{Background}

For $ m > 1 $ we consider a non-negative distributional solution \[ \rho \in L^m_{loc}((0,T) \times \Rn) \] of the porous medium equation (PME) \[ \dt \rho - \lapx \rho^m = 0 \] on $ (0,T) \times \Rn $, where $ 0 < T \leq \infty $. The naming is explained by the physical origin of the equation as a model for the density of diffusing gas in a porous medium. Thanks to \citet{dibenfried_cnt} and \citet{dahlken_cnt} it is known that any solution is Hölder continuous on $ (0,T) \times \Rn $. Moreover, \citet{aroncaf_initialtrace} showed that $ \rho $ takes on a uniquely determined initial value $ \rho(0,\, \cdot \, ) $ in the sense of distributions that has the property of being a Borel measure with growth restriction 
	\begin{align} \label{eq_initialGrowth} 
		\supre[r > 1] r^{-n -\frac{2}{m-1}} \, \rho(0,B_r(0)) < \infty. \tag{$ \ast $} 
	\end{align}
	The converse problem of finding a solution of the PME to a given initial datum $ \rho_0 $ was solved positively by \citet{bencrandpierre_exsolutions} if $ \rho_0 $ satisfys \eqref{eq_initialGrowth}. Their work provides a maximal existence time $ T $ for which a lower bound can be calculated explicitly. If the initial datum is non-negative and locally integrable with an additional growth bound, one gets $ T = \infty $ and therefore global time solvability. 
That the initial datum determines the solution uniquely was shown by \citet{dahlken_pme} in full generality. \\

Considering the equation in divergence form, namely \[ \dt \rho - (m - 1) \, \divx\left(c_m \, \rho^{m-1} \, \gradx \rho\right) = 0, \] where we abbreviate \[ \frac{m}{m-1} \rdef c_m, \] it becomes evident that the PME is a degenerate parabolic equation: The diffusion coefficient \[ v \ldef c_m \, \rho^{m-1} \] representing the pressure of the gas vanishes as the density $ \rho $ approaches zero. The formal calculation for passing between the two instances of the equation is rigorous in case also \[ \gradx \rho \in L^m_{loc}((0,T) \times \Rn). \] 
As opposed to the situation for uniformely strongly parabolic equations, the degeneracy of the PME implies that any solution whose initial positivity set is compact retains this property for any finite time. 
In physical terms this means that the diffusing gas does not get to every point of space instantaneously, but rather propagates with finite speed. Not only does this paint a more realistic picture of the real world in terms of modelling diffusion processes, but it does also give rise to an interesting mathematical phenomenon: The time-space positivity set $ \Pos(\rho) $ of a solution, open because of the continuity of $ \rho $, has a non-empty boundary that separates it from the time-space region where $ \rho $ vanishes, thus constituting a sharply defined interface \[ \Interf(\rho) \ldef \partial \Pos(\rho) \cap ((0,T) \times \Rn). \] 
The regularity of the interface and the regularity of the solution are closely connected. Note that parabolic regularity theory 
implies the smoothness of solutions on $ \Pos(\rho) $, so regularity is only an issue near $ \Interf(\rho) $. There are, on the other hand, examples of explicit solutions with smooth interface that are not more than time-space Hölder continuous across the interface. \\
In the case of one space dimension, the interface is always Lipschitz regular as was shown in \citet{aronson_lipInterfOneD}. This is optimal according to \citet{aroncafvaz_cornerInterfOneD}, since when starting to move after a waiting time the interface may have a corner. However, once the waiting time has passed, the interface is not only smooth \citep{aronvaz_smoothInterfOneD}), but even real analytic \citep{angenent_analyticInterfOneD}). Furthermore, according to \citet{benilan_temporalLipschitzOneD} and \citet{aronson_regularityOneD}, on top of the Hölder regularity of the density, the pressure is Lipschitz continuous everywhere in time and space. \\
In the general case of arbitrary dimensions, some irregularities can appear. In \citet{cafvazwolan_lipcntn} it was shown that after waiting the time it takes the gas to overflow the ball containing the initial support, the pressure of a solution is a Lipschitz function in time and space and on all of $ \Rn $, that is especially across the interface. It follows that $ \Interf(\rho) $ can be described as a Lipschitz continuous surface for sufficiently large times. As noted above, the temporal constraint solely applies for dimensions $ n > 1 $ and is important in case the positivity set of the initial datum contains one or more holes. Any hole is filled in finite time, but advancing interfaces may hit each other and the velocity of the interface can become unbounded at the focussing time, as described for example in \citet{aronGrav}. \\ 
Under certain conditions on the initial pressure, more regularity can be gained. Imposing \[ \rho_0^{m-1} + |\gradx \rho_0^{m-1}| \geq c > 0 \text{ on } \Pos(\rho_0) \] for a constant $ c > 0 $ ensures that the interface starts to move at all of its points immediately at the initial time. This becomes plausible by the physical interpretation of the equation that suggests to view the spatial derivative of the pressure as the velocity of the extension of gas in space and hence the speed of the interface. If in addition $ \rho_0 $ has a bounded positivity set 
and satisfies \[ \rho_0^{m-1} \in \stetig^1(\overline{\Pos}(\rho_0)) \] as well as \[ \lapx \rho_0^{m-1} \geq - C \, \text{ in } \distrib(\Pos(\rho_0)) \] for a constant $ C > 0 $, we speak of a non-degenerate initial datum. For the density belonging to such a $ \rho_0 $, \citet{cafvazwolan_lipcntn} showed that the spatial gradient of the pressure is bounded away from zero near the interface for sufficiently large times, making it possible for \citet{koch_habil} to establish that for sufficiently large times the interface is smooth and the pressure is smooth up to and including the interface. \\
Short time smoothness of the pressure up to and including the interface as well as smoothness of the interface itself, both before a possible blow-up time, hold according to \citet{daskahamil_regfreebdry} and \citet{koch_habil} for initial data that start to move immediately and have an initial pressure with Hölder continuous derivative. 
To complete the picture, \citet{daskahamilee} found that non-degenerate initial data which in addition possess a weakly concave square root function of the initial pressure generate solutions with convex positivity set for all times and hence smoothness of the pressure up to and including the interface as well as smoothness of the interface follow for all times. 
 
\subsection{Main Results}

The main result of this paper is a global regularity estimate for any number of derivatives of the pressure of a distributional solution, imposing only Lipschitz conditions on the initial pressure. The space of not necessarily bounded Lipschitz functions on an arbitrary set $ \Omega \subset \Rn $ is denoted by $ \homLip(\Omega) $. It is obviously a subset of $ \stetig(\Omega) $ and hence $ \Loneloc(\Omega) $, and the distributional gradient of any Lipschitz function in the interior of $ \Omega $ is itself a regular distribution in $ \Linfty(\Omega) $. 


\begin{theo} \label{theo_main} 
	There exists an $ \varepsilon > 0 $ such that the following holds: \\ 
	If $ \rho_0^{m - 1} \in \homLip(\Rn) $ satisfies \[ \norm[\Linfty(\Pos(\rho_0))]{c_m \, \gradx \rho_0^{m - 1} -  \en} \leq \varepsilon \] and if $ \rho $ is the solution of the PME on $ (0,T) \times \Rn $ with initial value $ \rho_0 $, then \[ \rho^{m - 1} \in \glatt(\overline{\Pos}(\rho)) \] and \[ \supre[(t,x) \in \Pos(\rho)] t^{k+|\alpha|} \left| \dt^k \dx^\alpha \left( c_m \, \gradx \rho^{m - 1}(t,x) - \en \right) \right| \leq c \, \norm[\Linfty(\Pos(\rho_0))]{c_m \, \gradx \rho_0^{m - 1} -  \en} \] for any $ k \in \N_0 $ and $ \alpha \in \N_0^n $ with a constant $ c = c(n,m,k,\alpha) $. Furthermore, all level sets of $ \rho $ are analytic. 
\end{theo}

Note that the theorem contains the analyticity of $ \Interf(\rho) $, which is the level set belonging to zero. We conjecture that also the pressure itself is analytic rather than merely smooth. \\ 

Theorem \ref{theo_main} can be interpreted as a stability result: For the perfectly flat travelling wave front solution \[ \rho_{tw}(t,x) = \left(c_m^{-1} \, (x_n + t)_+\right)^\frac{1}{m-1} \] normed with respect to the $ n $-th coordinate direction, the derivative on its positivity set is given by the $ n $-th unit vector $ \en $. Seen in this light, the result asserts that solutions of the PME whose pressure is initially close to a flat front in the Lipschitz sense stay close to it for all times. \\

In order to prove Theorem \ref{theo_main}, we consider the so-called transformed pressure equation (TPE) \[ \ds w - y_n \, \lapy' w + y_n^{-\sigma} \, \dyn \left( y_n^{1 + \sigma} \, \frac{1 + |\grad_y' w|^2}{\dyn w} \right) = 0  
\text{ on } (0,S) \times \Hquer, \] where \[ 
0 < S \leq \infty, 
\, H\ldef \menge{y \in \Rn \mid y_n > 0 } \text{ and } 
\sigma > - 1. 
\] 
We establish existence of solutions that possess good regularity properties and are Lipschitz stable with respect to the travelling wave solution \[ w_{tw}(s,y) \ldef   y_n - ( 1 + \sigma ) \, s \] with $ \grady w_{tw} = \en $.

\begin{theo} \label{theo_TPE}
	There exists an $ \varepsilon > 0 $ such that for any $ w_0 \in \homLip(\Hquer) $ satisfying \[ \norm[\Linfty(H)]{\grady \left( w_0 - w_{tw}(0,\, \cdot\, ) \right)} \leq \varepsilon \] we can find a solution $ w_* \in \glatt((0,S) \times \Hquer) $ of the TPE on $ (0,S) \times \Hquer $ 
	with initial value $ w_0 $ for which we have \[ \supre[(s,y) \in (0,S) \times H] s^{k+|\alpha|} \left|\ds^k \dy^\alpha \grady \left( w_*(s,y) - w_{tw}(s,y) \right) \right| \leq c \, \norm[\Linfty(H)]{\grady \left( w_0 - w_{tw}(0,\, \cdot\, ) \right)} \] for any $ k \in \N_0 $ and $ \alpha \in \N_0^n $ with a constant $ c = c(n,\sigma,k,\alpha) $. Furthermore, $ w_* $ is analytic in the temporal and tangential directions on $ (0,S) \times \Hquer $ with $ \Lambda > 0 $ and $ C = C(n) > 0 $ such that \[ \supre[(s,y) \in (0,S) \times H] s^{k + |\alpha'|} \, |\ds^k \dystrich^{\alpha'} \grady w_*(s,y)| \leq C \,  \Lambda^{-k-|\alpha'|} \, k! \, \alpha'! \, \norm[\Linfty(H)]{\grady \left(w_0 - w_{tw}(0,\, \cdot\, ) \right)} \] for any $ k \in \N_0 $ and $ \alpha' \in \N_0^{n-1} $ with $ k + |\alpha'| > 0 $.
\end{theo} 

If $ \varepsilon $ is chosen possibly smaller so that for $ k = |\alpha| = 0 $ we have $ c \, \varepsilon < 1 $, then from the global bound on the gradient of $ w_* $ it follows that $ |\dyn w_*| > 1 - c \, \varepsilon $ on $ (0,S) \times \Hquer $. As a consequence, for any $ x_n \in \R $ we can locally everywhere solve the graph equation $ w_*(s,y) = x_n $ for the implicit variable $ y_n $ by means of a unique smooth function $ v_*(s,y',x_n) $ with \[ \grad_{s,y'} v_* = - (\dyn w_*)^{-1} \, \grad_{s,y'} w_* \text{ and } \dxn v_* = (\dyn w_*)^{-1}. \] This suggests a reparametrisation of the graph in terms of a local change of coordinates performed by means of a von-Mises-transformation \[ A: (s,y) \mapsto (s,y',w_*(s,y)) \rdef (t,x) \] interchanging dependent and independent variables. 
For a fixed time $ s \in (0,S) $, Theorem \ref{theo_TPE} also implies \[ \left| (A(s,y) - A(s,\yquer)) - ((s,y) - (s,\yquer)) \right| \leq c \, \varepsilon \, |y - \yquer| \] through the medium of the mean value theorem. Hence both \[ \left| (A(s,y) - A(s,\yquer)) - ((s,y) - (s,\yquer)) \right| \leq \frac{c \varepsilon}{1 - c \varepsilon} \, |y - \yquer| \] and \[ \left| (A(s,y) - A(s,\yquer)) - ((s,y) - (s,\yquer)) \right| \leq \frac{c \varepsilon}{1 - c \varepsilon} \, |A(s,y) - A(s,\yquer)| \] hold, and we can conclude \[ (1 - c \, \varepsilon) \, |y - \yquer| < |A(s,y) - A(s,\yquer)| < (1 - c \, \varepsilon)^{-1} \, |y - \yquer| \] for any $ y, \, \yquer \in \Hquer $. Consequently, $ A $ is an injective quasi-isometry that allows us to reparametrise the graph globally with $ A(s,\Hquer) = \overline{\Pos}(v_*(s)) $. After an additional rescaling in time by $ (1 + \sigma) $ without changing the notation of the time variable $ t $ it becomes evident that $ v_* $ solves the equation \[ \dt v - (m-1) \, v \, \lapx v - | \gradx v |^2 = 0 \] classically on $ \overline{\Pos}(v_*) $ up to and including the boundary, where \[ m = \frac{2 + \sigma}{1 + \sigma}. \] However, this equation is nothing but the porous medium pressure equation (PMPE) with parameter $ m > 1 $, explaining the name of the TPE. Since the level set of $ v_* $ at height $ y_n $ is given by \[ \menge{(s,y',x_n) \mid x_n = w_*(s,y)}, \] the temporal and tangential analyticity of $ w_* $ translates into analyticity of the level sets of $ v_* $. Note that for these transformations to work it is necessary that the distance of the initial data from the initial travelling wave is small in the homogeneous Lipschitz norm, thus in this sense the smallness condition in Theorems \ref{theo_TPE} and \ref{theo_main} is optimal. The repetition of the process just described with $ w_{tw} $ instead of $ w_* $ generates the travelling wave solution \[ v_{tw}(t,x) = x_n + t \] of the PMPE on $ \overline{\Pos}(v_{tw}) $. \\
Setting \[ \rho(t,x) \ldef 
	\begin{cases}
		\left( c_m^{-1} \, v(t,x) \right)^\frac{1}{m-1} &\text{ for any } (t,x) \in \Pos(v), \\
		0  & \text{ else } 
	\end{cases} 
\] for $ v = v_* $ and $ v = v_{tw} $ yields a classical solution $ \rho_* $ of the PME on $ \overline{\Pos}(\rho_*) $ as well as the travelling wave solution $ \rho_{tw} $ of the PME we started out with. But $ \rho_* $ is indeed a distributional solution of the PME on $ (0,T) \times \Rn $ for $ T = (m - 1)^{-1} S $, as can be seen by the simple computation 
\begin{align*}
	\inte{(0,T) \times \Rn}{\rho_* \dt \varphi}{\Lebnp} + \inte{(0,T) \times \Rn}{\rho_*^m \lapx \varphi}{\Lebnp} = - \inte{\Pos(\rho_*)}{\dt \rho_* \varphi}{\Lebnp} + \inte{\Pos(\rho_*)}{\lapx \rho_*^m \varphi}{\Lebnp} = 0 
\end{align*} 
for any $ \varphi \in \glatt_c((0,T) \times \Rn) $, where the spatial boundary terms of the integrations by parts vanish since $ \rho_* $ does so at the boundary of its positivity set and we use the fact \[ \gradx \rho_*^m = c_m^{-c_m} \gradx v_*^{c_m} = \rho_* \gradx v_*. \] In conjunction with the existence and uniqueness results for the PME, this shows Theorem \ref{theo_main}. 
The remainder of this work is concerned with proving Theorem \ref{theo_TPE}.
 
\section{The Perturbation Equation} \label{section_pe}
\subsection{The Perturbational Setting} \label{subsection_perturbSetting}

When considering the perturbed travelling wave $ w_{tw} + u $, formally the perturbation $ u $ is itself a solution of the TPE if it satisfies the perturbation equation (PE) \[ \ds u - L_\sigma u = f[u] \] on $ \omega $ with linear spatial part \[ y_n \, \lap_y u + ( 1 + \sigma ) \, \dyn u \rdef L_\sigma u \] and non-linearity \[ - ( 1 + \sigma ) \, \frac{|\grady u|^2}{\dyn u + 1} - y_n \, \dyn \frac{|\grady u|^2}{\dyn u + 1} \rdef f[u]. \] It is worth pointing out that the PE possesses an invariant scaling: If $ u $ is a solution with respect to $ (s,y) $ with initial value $ u_0 $, then $ \frac{1}{\lambda} \left(u \circ A_\lambda \right) $ is a solution with respect to $ (\hat{s},\hat{y}) $ with initial value $ \hat{u}_0 $, where \[ A_\lambda: (\hat{s},\hat{y}) \mapsto (\lambda \, \hat{s},\lambda \, \hat{y}) \rdef (s,y) \] and $ \hat{u}_0(\hat{y}) \ldef u_0(\lambda \, \hat{y}) $. \\

The main theorem for the PE is the following.

\begin{theo} \label{theo_PE} 
	There exists an $ \varepsilon > 0 $ such that for any $ u_0 \in \homLip(\Hquer) $ with \[ \norm[\Linfty(H)]{\grady u_0} \leq \varepsilon \] we can find a solution $ u_* \in \glatt((0,S) \times \Hquer) $ of the  PE on $ (0,S) \times \Hquer $ with initial value $ u_0 $ satisfying \[ \supre[(s,y) \in (0,S) \times H] s^{k+|\alpha|} \left|\ds^k \dy^\alpha \grady u_*(s,y) \right| \leq c \, \norm[\Linfty(H)]{\grady u_0} \] for any $ k \in \N_0 $ and $ \alpha \in \N_0^n $ with a constant $ c = c(n,\sigma,k,\alpha) $. Furthermore, $ u_* $ is analytic in the temporal and tangential directions on $ (0,S) \times \Hquer $ with a $ \Lambda > 0 $ and a $ C = C(n) > 0 $ such that \[ \supre[(s,y) \in (0,S) \times H] s^{k + |\alpha'|} \, |\ds^k \dystrich^{\alpha'} \grady u_*(s,y)| \leq C \, \Lambda^{-k-|\alpha'|} \, k! \, \alpha'! \, \norm[\Linfty(H)]{\grady u_0} \] for any $ k \in \N_0 $ and $ \alpha' \in \N_0^{n-1} $.
\end{theo}

This generates a solution $ w_* = w_{tw} + u_* $ of the TPE, implying Theorem \ref{theo_TPE} immediately.

\begin{bem} \label{bem_fullAnalytic}
	It is shown in \citet{koch_habil} that the solution $ u_* $ and all its temporal derivatives are also analytic on $ (0,S) \times \Hquer $ in any spatial direction including the vertical one. We conjecture that $ u_* $ is indeed analytic on $ (0,S) \times \Hquer $ in time and space. This would ultimately also show the analyticity of the pressure. 
\end{bem}

The proof of Theorem \ref{theo_PE} is given in Section \ref{section_nonlinear}. It closely follows \citet{kochlamm} and uses a fixed point argument in a special function space. This also requires a thorough treatment of the linearised PE (LPE), where we simply ignore the dependence of $ f[u] $ on $ u $. It has the same invariant scaling as the non-linear PE. Moreover, translations in any temporal and spatial direction save the $ y_n $-direction commute with the differential operator. Starting from a suitable definition of weak solution, we will obtain the linear estimates needed in Section \ref{section_linear}. \\
Note also that we can express both the spatial part of the operator and the non-linearity in divergence form as \[ L_\sigma u =  y_n^{-\sigma} \, \divy \left( y_n^{1 + \sigma} \, \grady u \right) \] and \[ f[u] = - y_n^{-\sigma} \dyn \left( y_n^{1 + \sigma} \frac{|\grady u|^2}{\dyn u + 1} \right). \] Although we will not make use of this fact structurally, it motivates the definition of a weak notion of solution and makes the use of weighted measures natural.


\subsection{Weighted Measures}
Set $ \musig(y) \ldef y_n^\sigma \, d\Lebn(y) $. Then $ \musig $ is a countably finite Radon measure possessing the same nullsets as the Lebesgue measure. We denote $ \musig(\Omega) \rdef |\Omega|_\sigma $. As an abbreviation for the Lebesgue spaces with respect to $ \musig $ on an arbitrary set $ \Omega \subset \Rn $ we set \[ \Lp(\Omega,\musig) \rdef \Lpsig(\Omega) \text{ for } 1 \leq p \leq \infty \] and drop the measure from the notation in case of the Lebesgue measure $ \sigma = 0 $. Since the nullsets coincide, we have \[ \Linfty_\sigma(\Omega) = \Linfty(\Omega). \] We would also like to introduce Sobolev spaces that allow for different weights in every order of derivatives by defining \[ \Wmpsigvec(\Omega) \ldef \menge{ u \in \Lp_{\sigma_0}(\Omega) \mid \dy^\alpha u \in \Lpsiga(\Omega) \text{ for all } |\alpha| \leq m} \] for $ m \in \N $, $ 1 \leq p < \infty $ and an open set $ \Omega \subset H $, where the weight exponents $ \sigma_0, \ldots, \sigma_m > - 1 $ are understood as the vector $ (\sigma_0,\ldots,\sigma_m) = \vec{\sigma} $. Here the derivatives are taken in the distributional sense, which is possible since on an open set $ \Omega $ that is contained in $ H $ we have \[ \Lpsig(\Omega) \subset \Lploc(\overline{\Omega},\musig) \subset \Lploc(\Omega,\musig) \subset \Loneloc(\Omega,\musig) \subset \Loneloc(\Omega) \]  for any $ 1 \leq p \leq \infty $. Both weighted Lebesgue and Sobolev spaces possess all the usual functional theoretical properties. The inner product on $ \Ltwosig(\Omega) $ is denoted as $ \ip[\Ltwosig(\Omega)]{\, \cdot \,}{\, \cdot \,} $.

\subsection{Intrinsic Metric}

A crucial observation is that the linear spatial operator $ L_\sigma $ gives rise to a Carnot-Caratheodory-metric $ \dist $ on $ \Hquer $. This intrinsic metric is described in detail in \citet{daskahamil_regfreebdry} in two dimensions, and in \citet{koch_habil} and \cite{diss} in general. 
As a manifestation of the degeneracy of the equation, it is singular towards $ \partial H $. An equivalent characterisation in terms of an explicit expression is given for all $ y, \, z \in \Hquer $ by \[ \constd^{-1} \, \dist(y,z) \leq \frac{|y-z|}{\sqrt{y_n} + \sqrt{z_n} + \sqrt{|y-z|}} \leq \dist(y,z) \] with $ \constd \ldef 12 $. It is sometimes convenient to consider the quasi-metric \[ \quasi(y,z) \ldef \frac{|y-z|}{(y_n^2 + z_n^2 + |y-z|^2)^\frac{1}{4}} \] instead. Then we have \[ \constd^{-1} \, \dist(y,z) \leq \quasi(y,z) \leq 4 \, \dist(y,z) \] for any $ y, \, z \in \Hquer $. All balls $ B_r(z) $ will be understood with respect to $ \dist $, while $ \Beu_r(z) $ denotes Euclidean balls intersected with the closed upper half plane $ \Hquer $. Their interplay is displayed by the inclusions \[ \Beu_{\constd^{-2}r(r + \sqrt{z_n})}(z) \subset B_r(z) \subset \Beu_{2r(r + 2  \sqrt{z_n})}(z). \] Close to the boundary of the upper half plane, that is for $ r > \sqrt{z_n} $, we observe that $ B_r(z) $ behaves like $ \Beu_{r^2}(z) $, while for $ r < \sqrt{z_n} $ the intrinsic geometry is comparable to the Euclidean one. This reflects the uniform parabolicity of our equation away from the free boundary. For $ y_0 \in \Hquer $ and a parameter $ \delta_1 > 0 $, it is also possible to construct a cutoff function \[ \eta \in \glatt_c(B_{\delta_1 r}(y_0)) \] with \[ \eta = 1 \text{ on } \overline{B}_{\delta_2 r}(y_0) \] and \[ |\dy^\alpha \eta| \lesssim r^{- |\alpha|} \, (r + \sqrt{\yzeron})^{-|\alpha|} \text{ on } B_{\delta_1 r}(y_0) \] if $ \delta_2 $ is small enough compared to $ \delta_1 $. \\ 
In terms of the weighted measure we have \[ |B_r(z)|_\sigma \eqsim_{n,\sigma} r^n (r + \sqrt{z_n})^{n+2\sigma}. \] Note also that the intrinsic metric turns the weighted measure space $ (H,\musig) $ into a space of homogeneous type, satisfying the doubling condition \[ |B_{\kappa r}(z)|_\sigma \leq 
c \, \kappa^n \, (1 + \kappa^{n + 2 \sigma}) |B_r(z)|_\sigma \] for $ \kappa > 1 $ with $ c = c(n, \sigma) > 0 $.


\subsection{Function Spaces}

The fixed point argument in the proof of Theorem \ref{theo_PE} takes place in a special function subspace $ X $ of the space of spatial Lipschitz functions that has to be complemented by an intermediate function space $ Y $ for the linear considerations. In \citet{kochtata}, the choice of these function spaces was motivated by the square function characterisation of BMO. Since we aim at a Lipschitz setting, we first look for a bound of the homogeneous Lipschitz norm $ \norm[\Linfty((0,S) \times H)]{\grady u} $ of a solution $ u $ to the linear equation in terms of $ \norm[\Linfty(H)]{\grady u_0} $ and the $ Y $-norm $ \norm[Y]{f} $ of the inhomogeneity. It turns out that the right space to consider $ f $ in is split into an on-diagonal and an off-diagonal part given by
\begin{align*}
	\Yon{f} &\ldef \supre[0 < r^2 < S \atop z \in H] r^\theta \, |Q_r(z)|^{-\frac{1}{p}} \, \norm[L^{p}(Q_r(z))]{\grady f}, \\
	\Yoff{f} &\ldef \supre[0 < r^2 < S \atop z \in H] r^{2 - \varepsilon_1} \, (r + \sqrt{z_n})^{-\varepsilon_2} \, |Q_r(z)|^{-\frac{1}{p}} \norm[\Lp(Q_r(z))]{f},
\end{align*}
where $ Q_r(z) \ldef (\frac{1}{2} \, r^2,r^2) \times B_r(z) $ are intrinsic parabolic cylinders that stay away from the initial time. We set $ Y(p) \ldef Y^{\small{on}}_2 \cap Y^{\small{off}}_{1,\, 1} $. \\
In order to close the argument, we also need an estimate for the non-linearity. This is only possible when adding additional terms to the homogeneous Lipschitz norm and considering
\begin{align*} 
	\norm[X^{(1)}(p)]{u} &\ldef \norm[\Linfty((0,S) \times H)]{\grady u} \\
	&+ \supre[0 < r^2 < S \atop z \in H] \hspace{-0.1cm} \left(r \, (r + \sqrt{z_n}) \, |Q_r(z)|^{-\frac{1}{p}} \, \norm[\Lp(Q_r(z))]{D^2_y u} + r^2 \, |Q_r(z)|^{-\frac{1}{p}} \, \norm[\Lp(Q_r(z))]{y_n \, D^3_y u} \right).
\end{align*}

Thanks to the special structure of $ f[u] $, we can then establish mapping properties of $ u \mapsto f[u] $ that include the reverse inequality.

\begin{lemma} \label{lemma_YggX}
	If $ 1 \leq p \leq \infty $ and $ 0 < r < 1 $, then the mapping \[ \overline{B}^X_R \ldef \menge{u \in X^{(1)}(p) \mid \norm[X^{(1)}]{u} \leq R} \ni u \mapsto f[u] \in Y(p) \] is analytic and there exists a constant $ c = c(n,p) > 0 $ such that \[ \norm[Y(p)]{f[u]} \leq c \,  \frac{1}{(1-R)^3} \, \norm[X^{(1)}(p)]{u}^2 \text{ for all } u \in \overline{B}^X_R \] as well as \[ \norm[Y(p)]{f[u_1] - f[u_2]} \leq c \, \frac{R}{(1-R)^6} \, \norm[X^{(1)}(p)]{u_1 - u_2} \text{ for all } u_1, \, u_2 \in \overline{B}^X_R. \] 
\end{lemma}
\begin{bew}
	We introduce the symbolic notation $ D^{m_1}_y u \star D^{m_2}_y u $ to denote any arbitrary linear combination of products of derivatives of orders $ m_1 $ and $ m_2 $. 
	The iterated application of $ \star $ onto the same order of derivatives, as in $ D^m_y u \star \ldots \star D^m_y u $ ($ j $ times), is abbreviated by $ (D^m_y u)^{j \star} $ with the usual conventions $ (D^m_y u)^{1 \star} = 1 \star D^m_y u $ and $ (D^m_y u)^{0 \star} = 1 $. \\
	A direct calculation shows that \[ f[u] = f_1[u] \star \grady u \star \grady u + f_2[u] \star \grady u \star \Pn \, D^2_y u \] and, for $ j = 1, \ldots, n $, \[ \dyj f[u] = f_2[u] \star \grady u \star D^2_y u + f_2[u] \star \grady u \star \Pn \, D^3_y u + f_3[u] \star D^2_y u \star \Pn D^2_y u \] with 
	\begin{align*}
		f_i[u] \ldef \summe[i]{k=1} (\dyn u +1)^{-k} \star (\grady u)^{(k-1) \star}, \ i = 1, 2, 3.
	\end{align*}
	We now consider $ \norm[Y(p)]{f[u]} $ term by term, here concentrating on three examplatory cases. More details can be found in \citet{diss}. \\
	On the one hand, for arbitrary $ r > 0 $ and $ z \in \Hquer $ we get
	\begin{align*}
		r \, (r + \sqrt{z_n})^{-1} \, |Q_r(z)|^{-\frac{1}{p}} \, &\norm[L^{p}(Q_r(z))]{f_1[u] \star \grady u \star \grady u} \\
		&\leq \norm[\Linfty((0,S) \times H)]{f_1[u]} \, \norm[\Linfty((0,S) \times H)]{\grady u}^2,
	\end{align*} 
	since $ \norm[L^{p}(Q_r(z))]{\grady u} \leq \norm[\Linfty((0,S) \times H)]{\gradx u} \, |Q_r(z)|^\frac{1}{p} $ and $ (r + \sqrt{z_n})^{-1} \leq r $, and
	\begin{align*}
		r \, (r + \sqrt{z_n})^{-1} \, &|Q_r(z)|^{-\frac{1}{p}} \, \norm[L^{p}(Q_r(z))]{f_2[u] \star \grady u \star \Pn \, D^2_y u} \\
		&\lesssim \norm[\Linfty((0,S) \times H)]{f_2[u]} \, \norm[\Linfty((0,S) \times H)]{\grady u} \, r \, (r + \sqrt{z_n}) \, |Q_r(z)|^{-\frac{1}{p}} \, \norm[L^{p}(Q_r(z))]{D^2_y u},
	\end{align*} 
	since $ \sqrt{y_n} \lesssim r + \sqrt{z_n} $ for any $ y \in B_r(z) $. \\
	On the other hand, the identity \[ \norm[\Lp(Q_r(z))]{\Pn \, |D^2_y u|^2} = \norm[L^{2p}(Q_r(z))]{\Pn^\frac{1}{2} \, D^2_y u}^2 \] shows that
	\begin{align*}
		r^2 \, |Q_r(z)|^{-\frac{1}{p}} \, &\norm[\Lp(Q_r(z))]{f_3[u] \star D^2_y u \star \Pn \, D^2_y u } \\
		&\leq \norm[\Linfty((0,S) \times H)]{f_3[u]} \, r^2 \, |Q_r(z)|^{-\frac{1}{p}} \, \norm[L^{2p}(Q_r(z))]{\Pn^\frac{1}{2} \, D^2_y u}^2.
	\end{align*}
	Applying the weighted Gagliardo-Nirenberg interpolation inequality onto $ \eta \, \grady u $, where $ \eta $ is a cutoff function adapted to the intrinsic metric, gives the upper bound 
	\begin{multline*}
		\norm[\Linfty((0,S) \times H)]{f_3[u]} \, \Bigl(\norm[\Linfty((0,S) \times H)]{\grady u} + r^2 \, |Q_r(z)|^{-\frac{1}{p}} \, \norm[\Lp(Q_r(z))]{\Pn \, D^3_y u} \\
		+ r \, (r + \sqrt{z_n}) \, |Q_r(z)|^{-\frac{1}{p}} \, \norm[\Lp(Q_r(z))]{D^2_y u} \Bigr)
	\end{multline*}
	with a constant depending on $ n $ and $ p $ only. \\
	We proceed similarily with the other terms and see that
	\[ \norm[Y(p)]{f[u]} \lesssim \left( \norm[\Linfty((0,S) \times H)]{f_1[u]} + \norm[\Linfty((0,S) \times H)]{f_2[u]} + \norm[\Linfty((0,S) \times H)]{f_3[u]} \right) \norm[X^{(1)}(p)]{u}^2, \]
	But by assumption we have \[ \norm[\Linfty((0,S) \times H)]{\grady u} \leq R \text{ and } \norm[\Linfty((0,S) \times H)]{\dyn u + 1}^{-k} \leq (1 - R)^{-k} \] and hence \[ \norm[\Linfty((0,S) \times H)]{f_i[u]} \leq \summe[i]{k=1} \frac{R^{k-1}}{(1-R)^k}. \] 
	The estimate for the difference follows by akin considerations. 
\end{bew}

Considering the linear equation again, we finally see that a solution $ u $ to an inhomogeneity $ f \in Y(p) $ is not only contained in $ X^{(1)} $, but in an even smaller space $ X(p) \ldef X^{(1)}(p) \cap X^{(2)}(p) $ with \[ X^{(2)}(p) \ldef \norm[\Linfty((0,S) \times H)]{\sqrt{s} \, \sqrt{y_n} \, D_y^2 u} + \supre[0 < r^2 < S \atop z \in H] r^2 \, |Q_r(z)|^{-\frac{1}{p}} \, \norm[\Lp(Q_r(z))]{\grady \ds u}. \] Obviously, $ \norm[X(p)]{\, \cdot \,} $ is a norm modulo constants. As an intersection of complete spaces, $ X(p) $ and $ Y(p) $ are also complete. Note that the notation does not reflect the dependence of the spaces on the interval. 
The term \enquote{solution} in the following theorem is to be understood in a sense that will be made precise later.

\begin{theo} \label{theo_linear}
	If $ \max\menge{2 \, (n + 1),(1+\sigma)^{-1}} < p < \infty $, $ f \in Y(p) $, $ u_0 \in \homLip(\Hquer) $ and $ u $ is a solution of the LPE to inhomogeneity $ f $ on $ (0,S) \times \Hquer $ with initial value $ u_0 $, then we have \[ \norm[X(p)]{u} \leq c \left(\norm[\Linfty(H)]{\grady u_0} + \norm[Y(p)]{f} \right) \] for a constant $ c = c(n,\sigma,p) > 0 $.
\end{theo}

While the estimate of the local $ \Lp $-norm in $ X^{(2)} $ here is a byproduct of the theory that will not be used any more in the course of the argument, the bound of the pointwise norm in $ X^{(2)} $ becomes important in the proof of Theorem \ref{theo_PE}. The proof of Theorem \ref{theo_linear} is given in Proposition \ref{prop_roughInitial}, Remark \ref{bem_defRoughInitial} and Proposition \ref{prop_XggY}.
 
\section{Linear Estimates} \label{section_linear}
\subsection{Definition of Solutions}
We first define a suitable weak notion of solution of the LPE on arbitrary time-space cylinders possibly touching the spatial boundary. Our particular interest in the behaviour of solutions towards $ \menge{y_n = 0} $ is the motivation to carry our considerations to the boundary of the upper half plane by applying a wider class of test functions that can attain non-zero values at $ \partial H $. For an arbitrary -- and not necessarily open -- set $ \Omega \subset \Rn $ we thus define
\begin{align*}
	\glatt_c(\Omega) \ldef \menge{ \varphi: \Omega \to \R \mid \supp \varphi \text{ is a compact subset of } \Omega, \, \varphi \text{ has an extension } \breve{\varphi} \in \glatt_c(\Rn)},
\end{align*}
and similarly $ \glatt_c(\omega) $ for an arbitrary $ \omega \subset \R \times \Rn $.

\begin{bem} \label{bem_density}
	It is clear that $ \glatt_c(\Omega) $ is a dense subset of $ \Lpsig(\Omega) $ for any $ \Omega \subset \Hquer $ and $ 1 \leq p < \infty $. By a modification of the usual arguments we can show that for $ \Omega = \Hquer $ we have that $ \glatt_c(\Hquer) $ is dense in $ \Wmpsigvec(H) $, given that $ 1 \leq p < \infty $ and $ \sigma_{i+1} - \sigma_i \leq 1 $ for $ i = 0, \ldots, m-1 $ as well as $ \sigma_m \geq 0 $. See also \citet{adams} and \citet{kufner_weightedSobolevSpaces}.
\end{bem}

In the following, $ \Omega $ will always be a relatively open subset of $ \Hquer $, while $ I \ldef (s_1,s_2) \subset \R $ denotes an arbitrary open interval with $ -\infty \leq s_1 < s_2 \leq \infty $. By the notation $ \Ilinks $ and $ \Irechts $ we mean the closure of $ I $ only at its left or right end point, respectively, where $ \pm \infty $ are always excluded from being an element of the interval. Evaluations at non-finite points are to be understood as limits. \\
For a reasonable theory we equip solutions with some additional regularity properties that allow for energy techniques. In this context, however, the conditions we set are the weakest possible.

\begin{defin} \label{defin_energySolution}
	\vspace{-0.2cm}
	\begin{itemize}
		\item
		Considering $ f \in \Loneloc\bigl(I; \, \Ltwosig(\Omega)\bigr) $, we say that $ u $ is a $ \sigma $-solution to $ f $ on $ I \times \Omega $, if and only if $ u \in \Ltwoloc\bigl(I; \, \Ltwosig(\Omega)\bigr) $, $ \grady u \in \Ltwo\bigl(I; \, \Ltwosigp(\Omega)\bigr) $ and \[ - \inte{I}{\ip[\Ltwosig(\Omega)]{u}{\ds \varphi}}{\Leb} + \inte{I}{\ip[\Ltwosigp(\Omega)]{\grady u}{\grady \varphi}}{\Leb} = \inte{I}{\ip[\Ltwosig(\Omega)]{f}{\varphi}}{\Leb} \] for all $ \varphi \in \glatt_c(I \times \Omega) $. 
		\item
		Considering $ f \in \Loneloc\bigl(\, \Ilinks; \, \Ltwosig(\Omega)\bigr) $, we say that $ u $ is a $ \sigma $-solution to $ f $ on $ I \times \Omega $ with initial value $ u_0 \in \Ltwosig(\Omega) $, if and only if $ u \in \Ltwoloc\bigl(\, \Ilinks ; \, \Ltwosig(\Omega)\bigr) $, $ \grady u \in \Ltwo\bigl(I; \, \Ltwosigp(\Omega)\bigr) $ and 
		\begin{align*}
			- \inte{I}{\ip[\Ltwosig(\Omega)]{u}{\ds \varphi}}{\Leb} + \inte{I}{\ip[\Ltwosigp(\Omega)]{\grady u}{\grady \varphi}}{\Leb} = \inte{I}{&\ip[\Ltwosig(\Omega)]{f}{\varphi}}{\Leb} \\
			&+ \ip[\Ltwosig(\Omega)]{u_0}{\varphi(s_1)} 
		\end{align*}
		for all $ \varphi \in \glatt_c(\, \Ilinks  \times \Omega) $. 
	\end{itemize}
\end{defin}

Out of consistency considerations, in case of $ s_1 = - \infty $ we only allow zero initial values and implicitely set $ u_0 = 0 $. 

\subsection{Energy Theory} \label{section_energyTheory}

For any relatively open $ \Omega \subset \Hquer $, existence of a $ \sigma $-solution with initial value $ u_0 \in \Ltwosig(\Omega) $ for the linear PE to $ f \in \Lone\bigl(I; \, \Ltwosig(\Omega)\bigr) $ can be shown by a Galerkin approximation. Uniqueness, however, is only guaranteed on cylinders without a spatial boundary.  

\begin{prop} \label{prop_eId}
	If $ f \in \Lone\bigl(I; \, \Ltwosig(H)\bigr) $ and $ u_0 \in \Ltwosig(H) $, then there exists a unique $ \sigma $-solution to $ f $ on $ I \times \Hquer $ with initial value $ u_0 $, for which $ u \in C_b\bigl(\Iquer; \, \Ltwosig(H)\bigr) $ with $ u(s_1) = u_0 $ holds as well as the energy identity \[ \frac{1}{2} \, \norm[\Ltwosig(H)]{u(s_2)}^2 + \inte{I}{\norm[\Ltwosigp(H)]{\grady u}^2}{\Leb} = \frac{1}{2} \, \norm[\Ltwosig(H)]{u(s_1)}^2 + \inte{I}{\ip[\Ltwosig(H)]{f}{u}}{\Leb}. \]
\end{prop} 
\begin{bew}
	First let $ s_1 > - \infty $. Since $ \Omega = \Hquer $, we can use Remark \ref{bem_density} to see that any \[ \varphi \in \Ltwo\bigl(I; \, \Ltwosig(H)\bigr) \text{ with } \grady \varphi \in \Ltwo\bigl(I; \, \Ltwosigp(H)\bigr) \text{ and } \ds \varphi \in \Ltwo\bigl(I; \, \Ltwosig(H)\bigr) \] serves as an admissible test function for the initial value problem, if in addition it has compact support in time contained in $ \Ilinks $. This last requirement, together with the temporal square-integrability of both $ \varphi $ and $ \ds \varphi $, implies $ \varphi \in \Linfty\bigl(I; \, \Ltwosig(H)\bigr) $, and the defining equation remains reasonable with such test functions for the class of inhomogeneities considered in Definition \ref{defin_energySolution}. \\
	Fixing $ \widetilde{s}_1 < \widetilde{s}_2 \in I $, so that $ (\widetilde{s}_1,\widetilde{s}_2) \rdef \widetilde{I} \subsetneq I $ has finite end points, a formal calculation with $ \chi_{\widetilde{I}} \, u $ as a test function implies the energy identity on $ \widetilde{I} $ immediatley. However, both $ \chi_{\widetilde{I}} $ and $ u $ do not possess the $ \Ltwo $-regularity of the time derivative needed for a justification. We regularise in time for $ h \geq 0 $ by \[ u^{{\scriptsize \pm} h}(s) \ldef 
	\begin{cases}
		{\scriptsize \pm} h^{-1} \inte[{\scriptsize \pm} h]{0}{u(s+\tau)}{\Leb(\tau)} &\text{ for all } s \in I^{{\scriptsize \pm} h} \\
		0 &\text{ for all } s \in I \ohne I^{{\scriptsize \pm }h},
	\end{cases} \]
	where $ I^h \ldef (s_1,s_2 - h) $ and $ I^{-h} \ldef (s_1 + h, s_2) $. By convention we set $ I^{-h} = I $ if $ s_1 = - \infty $ and $ I^h = I $ if $ s_2 = \infty $. For arbitrary $ h \in \R $ it is obvious that $ I^h + h = I^{-h} $. Choosing $ h > 0 $ so small that $ \widetilde{I} \subset I^h \cap I^{-h} $, the regularisation $ (\eta \, u^h)^{-h} $ is an admissible test function in the equation for $ u $ on $ I \times \Hquer $ for any $ \eta \in \Wonetwoloc(I) $ with $ \supp \eta \subset I^h \cap I^{-h} $. \\
	Now specify such a cutoff function $ \eta $ by defining $ \eta_\varepsilon \ldef \eta_{\varepsilon_1} \, \eta_{\varepsilon_2} $ for suitably small $ \varepsilon_1, \ \varepsilon_2 > 0 $, with \[ \eta_{\varepsilon_1}(s) \ldef 
	\begin{cases} 
		0 &\text{ for all } s \in (s_1,\widetilde{s}_1) \\
		\frac{s - \widetilde{s}_1}{\varepsilon_1} &\text{ for all } s \in (\widetilde{s}_1,\widetilde{s}_1 + \varepsilon_1) \\
		1 &\text{ for all } s \in (\widetilde{s}_1 + \varepsilon_1,s_2)
	\end{cases} 
	\] and \[ \eta_{\varepsilon_2}(s) \ldef 
	\begin{cases}
		1 &\text{ for all } s \in (s_1, \widetilde{s}_2 - \varepsilon_2) \\
		- \frac{s - \widetilde{s}_2}{\varepsilon_2} &\text{ for all } s \in (\widetilde{s}_2 - \varepsilon_2,\widetilde{s}_2) \\
		0 &\text{ for all } s \in (\widetilde{s}_2,s_2). 
	\end{cases} 
	\]
	A regularised version of the energy identity on $ \widetilde{I} $ in terms of $ h $, $ \varepsilon_1 $ and $ \varepsilon_2 $ follows, and the limit $ h \to 0 $ poses no difficulties. Moreover, one shows that \[ \widetilde{s}_2 \mapsto (\norm[\Ltwosig(H)]{u(\widetilde{s}_2)}^2)^{-\varepsilon_2} \] constitutes a Cauchy sequence in the function space $ C((\widetilde{s}_1,s_2)) $ and thus has a continuous limit \[ \widetilde{s}_2 \mapsto \norm[\Ltwosig(H)]{u(\widetilde{s}_2)}^2. \] 
	After the same considerations for $ \widetilde{s}_1 $, the norm continuity on $ [\widetilde{s}_1,\widetilde{s}_2] $ can be deduced as well as the energy identity on $ \widetilde{I} $. But then \[ \frac{1}{4} \, \supre[\widetilde{s}_2 \in (\widetilde{s}_1,s_2)] \norm[\Ltwosig(H)]{u(\widetilde{s}_2)}^2 + \inte{(\widetilde{s}_1,s_2)}{\norm[\Ltwosigp(H)]{\grady u}^2}{\Leb} \leq \frac{1}{2} \, \norm[\Ltwosig(H)]{u(\widetilde{s}_1)}^2 + \left( \inte{(\widetilde{s}_1,s_2)}{\norm[\Ltwosig(H)]{f}}{\Leb} \right)^2. \] Boundedness on $ (\widetilde{s}_1,s_2) $ is an immediate consequence, and hence both norm continuity and energy identity continue to hold for $ \widetilde{s}_2 \to s_2 $, reaching the end point if it is finite. \\ 

	If we repeat the whole process with the test function $ \eta_\varepsilon \, \widetilde{\varphi} $ for $ \widetilde{\varphi} \in \glatt_c(\Hquer) $, it becomes clear that \[ s \mapsto \ip[\Ltwosig(H)]{u(s)}{\widetilde{\varphi}} \] is continuous on $ \overline{[\widetilde{s}_1,s_2)} $. By the density from Remark \ref{bem_density}, this is nothing but weak continuity of \[ s \mapsto u(s) \in \Ltwosig(H) \text{ on } \overline{[\widetilde{s}_1,s_2)}, \] amounting to continuity on $ \Irechts $. \\ 
 
	For the inclusion of the initial time $ s_1 > - \infty $,  we now take into account that $ u $ possesses an initial value, fix an $ s_0 \in I $ and define \[ \eta_\varepsilon \ldef 
	\begin{cases}
		1 &\text{ for all } s \in (s_1, s_0) \\
		1 - \frac{s - s_0}{\varepsilon} &\text{ for all } s \in (s_0, s_0 + \varepsilon) \\ 
		0 &\text{ for all } s \in (s_0 + \varepsilon,s_2) 
	\end{cases}
	\]
	for small $ \varepsilon > 0 $. Then $ \widetilde{\varphi} \, \eta_\varepsilon $ is an admissible test function, and similar calculations as above lead to \[ \limes{s_0}{s_1} \ip[\Ltwosig(H)]{u(s_0)}{\widetilde{\varphi}} = \ip[\Ltwosig(H)]{u_0}{\widetilde{\varphi}}. \] This proves that weak continuity can be extended into $ s_1 $ with $ u(s_0) \rightharpoonup u_0 $ for $ s_0 \to s_1 $. But we know by the first part that the norm is uniformely bounded. Full continuity in $ \Ltwosig(H) $ down to $ s_1 $ with $ u(s_1) = u_0 $ is therefore proven, and we can let $ \widetilde{s}_1 \to s_1 $ in the energy identity. The uniqueness of $ \sigma $-solutions to the initial value problem follows directly. \\
	
	Finally, let $ s_1 = - \infty $. If we consider $ s^{(j)} \in I $ with $ s^{(j)} \goesto -\infty $ for $ j \goesto \infty $, the results obtained so far ensure the unique existence of $ \sigma $-solutions $ u^{(j)} $ on $ (s^{(j)},s_2) \times \Hquer $ with zero initial value. Extending them onto all of $ I $ by zero, we can apply the energy identity to see that those $ u^{(j)} $ form a Cauchy sequence in $ C(I; \Ltwosig(H)) $ whose limit will be denoted by $ u $. Then $ u $ is a $ \sigma $-solution on $ (-\infty,s_2) \times \Hquer $ with zero initial value. Continuity, energy identity, and hence uniqueness are implied immediately.
\end{bew}
%
\vspace{-0.2cm}
\begin{bem} \label{bem_decrease}
	Since the first part of the proof does not see the initial time, we have also shown that $ \sigma $-solutions without initial value are in $ C(\Irechts;\Ltwosig(H)) $, and both boundedness and the energy identity hold on $ (\widetilde{s}_1,s_2) $ for any $ \widetilde{s}_1 \in I $. Consequences of the energy identity, such as the temporal decrease of the $ \Ltwosig $-norm of $ \sigma $-solutions of the homogeneous linear PE on the whole space, are thus true on $ \Irechts $ in any case, and can be extended to $ \Iquer $ for $ \sigma $-solutions with initial value.
\end{bem}
\vspace{-0.2cm}
\begin{bem} \label{bem_duality}
	We can use the same method of proof to get the duality equality \[ \ip[\Ltwosig(H)]{u^{(1)}(s_1)}{u^{(2)}(s_2)} = \ip[\Ltwosig(H)]{u^{(1)}(s_2)}{u^{(2)}(s_1)} \text{ for all } s_1, \, s_2 \in \Iquer, \] where $ u^{(1)} $ and $ u^{(2)} $ are $ \sigma $-solutions of the homogeneous equation on $ I \times \Hquer $ that possess initial values. 
\end{bem}


If $ f \in \Ltwo \bigl(I; \, \Ltwosig(H)\bigr) $, additional regularity properties of $ \sigma $-solutions are revealed by the following weighted energy estimates. 
In order to prove them,  we consider only solutions with zero initial values $ u_0 = 0 $. This could be generalised under some regularity conditions on $ u_0 $. 
The statement alongside with a formal proof is already contained in \citet{koch_habil}. 

\begin{prop} \label{prop_energyEst} 
	If $ f \in \Ltwo \bigl(I; \, \Ltwosig(H)\bigr) $ and $ u $ is a $ \sigma $-solution to $ f $ on $ I \times \Hquer $ with zero initial value, then \[ s \mapsto \norm[\Ltwosigp(H)]{\grady u(s)} \in C_b(\Iquer) \text{ with } \grady u(s_1) = 0 \] and there exists a constant $ c = c(n,\sigma) > 0 $ such that \[ \inte{I}{\norm[\Ltwosig(H)]{\ds u}^2}{\Leb} + \inte{I}{\norm[\Ltwosig(H)]{\grady u}^2}{\Leb} + \inte{I}{\norm[\Ltwo_{2 + \sigma}(H)]{D^2_y u}^2}{\Leb} \leq c \, \inte{I}{\norm[\Ltwosig(H)]{f}^2}{\Leb}. \] 
\end{prop} 
\begin{bew}
	Fix $ \widetilde{s}_1 < \widetilde{s}_2 \in I $ and write $ \widetilde{I} \ldef (\widetilde{s}_1,\widetilde{s}_2) $. 
	In the formal proof for the temporal derivative we test the equation with $ \chi_{\widetilde{I}} \, \ds u $. This can be made rigorous similar as above, using on $ u^h $ that any $ \sigma $-solution $ \widetilde{u} $ with the additional property $ \ds \widetilde{u} \in  \Ltwoloc \bigl(I; \, \Ltwosig(H)\bigr) $ equivalently satisfies \[ \inte{I}{\ip[\Ltwosig(H)]{\ds \widetilde{u}}{\varphi}}{\Leb} + \inte{I}{\ip[\Ltwosigp(H)]{\grady \widetilde{u}}{\grady \varphi}}{\Leb} = \inte{I}{\ip[\Ltwosig(H)]{f}{\varphi}}{\Leb} \] for all $ \varphi \in \Ltwo \bigl(I; \, \Ltwosig(H)\bigr) $ with $ \grady \varphi \in \Ltwo \bigl(I; \, \Ltwosigp(H)\bigr) $ and compact temporal support. From this we get the continuity of \[ s \mapsto \norm[\Ltwosigp(H)]{\grady u(s)} \text{ on } \Irechts \] as well as boundedness of this function and the inequality for the temporal derivative on $ (\widetilde{s}_1,s_2) $ with an additional summand \[ \norm[\Ltwosigp(H)]{\grady u(\widetilde{s}_1)}^2 \] on the right hand side. The zero initial value now makes it possible to extend the $ \sigma $-solutions by zero onto a bigger interval and thus let $ \widetilde{s}_1 \to s_1 $ in all considerations. \\
	This weak regularity gain in the temporal derivative now enables us to restrict ourselves to elliptic equations: If $ u $ is a $ \sigma $-solution to $ f $ on $ I \times \Hquer $, then for almost all $ s \in I $ we have that $ u(s) $ satisfies \[ - L_\sigma u(s) = f(s) - \ds u(s) \rdef \widetilde{f}(s) \text{ on } \Hquer. \] 
	In an abuse of notation we will thus supress the time dependence and consider $ u \in \Ltwosig(H) $ with $ \grady u \in \Ltwosigp(H) $ and  \[ \inte{H}{\grady u \cdot \grady \varphi}{\musigp} = \inte{H}{\widetilde{f} \, \varphi}{\musig}, \, \varphi \in \glatt_c(\Hquer), \] 
	 as the energy formulation of the elliptic equation. As before, the density statement in Remark \ref{bem_density} implies that any function in $ \Ltwosig(H) $ with derivative in $ \Ltwosigp(H) $ is an admissible test function. \\
	Formally, the estimate for both the tangential and vertical spatial derivatives follows by testing the elliptic equation with $ \dyn u $, while for the second order derivative we consider $ y_n \lapy u $ as a test function. To make this rigorous we 
	perform a Fourier transformation in the tangential directions with Fourier variables $ \xi' \in \Rnm $ without renaming the functions $ u $ and $ \widetilde{f} $, coupled with a linear transformation $ z \ldef |\xi'| \, y_n $, to get 
	\begin{align} \tag{$\ast$} \label{ode}
		z \, \dz^2 u + (1 + \sigma) \, \dz u - z \, u = -\widetilde{f} \text{ on } (0,\infty),
	\end{align}  
	where we considered $ \xi' $ as parameters and set $ u(z) \ldef |\xi'| \, u(\xi',z) $.  A solution $ u $ of the homogeneous version of this equation defines a solution $ v = z^\frac{\sigma}{2} \, u $ of the modified Bessel equation with parameter $ \frac{\sigma}{2} $, that is \[ z^2 \, \dz^2 v + z \, \dz v - z^2 \, v - \frac{\sigma^2}{4} \, v = 0, \] and vice versa. The modified Bessel functions $ I_\frac{\sigma}{2} $ and $ K_\frac{\sigma}{2} $, described in detail in \citet{BesselFunctions}, form a fundamental system of this ordinary differential equation, hence a fundamental system for the homogeneous equation (\ref{ode}) is given by \[ \Psi_1(z) \ldef z^{-\frac{\sigma}{2}} \, I_\frac{\sigma}{2} \] and \[ \Psi_2(z) \ldef z^{-\frac{\sigma}{2}} \, K_\frac{\sigma}{2}. \] The asymptotics of the modified Bessel functions are known, and up to constants depending on $ \sigma $, for $ \sigma > 0 $ we get
	\begin{align*}
		\Psi_1(z) \sim 1 \quad (z \goesto 0), \qquad & \Psi_1(z) \sim z^{-\frac{1+\sigma}{2}} e^z \quad (z \goesto \infty) \\
		\Psi_2(z) \sim z^{-\sigma} \quad (z \goesto 0), \qquad & \Psi_2(z) \sim z^{-\frac{1+\sigma}{2}} e^{-z} \quad (z \goesto \infty). 
	\end{align*} 
	For $ \sigma < 0 $ we have $ \Psi_2(z) \sim 1 \; (z \goesto 0) $, and  $ \Psi_2(z) \sim \ln z \; (z \goesto 0) $ for $ \sigma = 0 $, while the other three relations remain as above. The Wronskian of $ \Psi_1 $ and $ \Psi_2 $ can be computed to be $ z^{-1-\sigma} $. \\
	All this leads to the fundamental solution \[ k(z,x) \ldef
	\begin{cases}
		x^\sigma \, \Psi_1(z) \, \Psi_2(x), \quad &z < x \\
		x^\sigma \, \Psi_1(x) \, \Psi_2(z), \quad &z > x, \\
	\end{cases} \]
	with first order derivative having a jump discontinuity of the type $ x^{-1} $ at $ z = x $. Therefore, solutions $ u $ to (\ref{ode}) are characterised by the representation \[ z^l \, u(z) = - \inte{(0,\infty)}{z^l \, k(z,x) \, \widetilde{f}(x)}{\Leb(x)} \] for any $ l \in \R $. We rewrite this to get the operator \[ \widetilde{f} \mapsto z^l \, u = - \inte{(0,\infty)}{z^l \, k(z,x) \, x^{-\sigma} \, \widetilde{f}(x)}{\musig(x)}. \] The definition of the fundamental solution and the asymptotic expansions of $ \Psi_1 $ and $ \Psi_2 $ ensure that \[ \supre [z \in (0,\infty)] \inte{(0,\infty)}{z^l \, |k(z,x)| \, x^{-\sigma}}{\musig(x)} < \infty \] and \[ \supre [x \in (0,\infty)] \inte{(0,\infty)}{z^l \, |k(z,x)| \, x^{-\sigma}}{\musig(z)} < \infty \] in case of $ l \in \menge{0, 1} $, verifying the conditions for the application of Schur's Lemma \citep{folland_RealAnalysis}). For solutions $ u $ of (\ref{ode}) it follows that \[ \norm[\Ltwosig((0,\infty))]{u} + \norm[\Ltwosig((0,\infty))]{z \,  u} \lesssim\norm[\Ltwosig((0,\infty))]{\widetilde{f}} \] with a constant depending on $ \sigma $ only. \\
	Now it makes sense to incorporate $ z \, u $ into the right hand side of (\ref{ode}). This results in a first order ordinary differential equation for $ \dz u \rdef w $, namely
	\begin{align} \tag{$ \ast \, \ast $} \label{ode2}
		z \, \dz w + (1 + \sigma) \, w  = -\widetilde{f} + z \, u \rdef \overline{f}.
	\end{align}
	A solution to the homogeneous equation is clearly given by $ z \mapsto z^{-1-\sigma} $, and so for this equation we get the fundamental solution \[ \overline{k}(z,x) \ldef 
	\begin{cases}
		x^\sigma \, z^{-1-\sigma}, \quad &z > x \\
		0, \quad &z < x. 
	\end{cases} \] 
	This time we consider the operator \[ z^\delta \, \overline{f}(z) \mapsto z^\delta \, w(z) = \inte{(0,\infty)}{\left(\frac{z}{x}\right)^\delta \, \overline{k}(z,x) \, x^\delta \, \overline{f}(x)}{\Leb(x)}. \] For $ \delta = \delta_1 $ we have 
	\begin{align*}
		\supre[z \in (0,\infty)] \inte{(0,\infty)}{\left(\frac{z}{x}\right)^{\delta_1} \, |\overline{k}(z,x)|}{\Leb(x)} &= \supre[z \in (0,\infty)] z^{\delta_1 - 1 - \sigma} \inte{(0,z)}{x^{\sigma -\delta_1}}{\Leb(x)} = \frac{1}{1 + \sigma - \delta_1}  
	\end{align*}
	if $ \delta_1 < 1 + \sigma $, showing that  \[ \norm[\Linfty((0,\infty))]{z^{\delta_1} \, w} \lesssim \norm[\Linfty((0,\infty))]{z^{\delta_1} \, \overline{f}} \] with a constant depending only on $ \sigma $ and $ \delta_1 $. \\
	Similarly, for $ \delta = \delta_2 $ we see 
	\begin{align*}
		\supre[x \in (0,\infty)] \inte{(0,\infty)}{\left(\frac{z}{x}\right)^{\delta_2} \, |\overline{k}(z,x)|}{\Leb(z)} &= \supre[x \in (0,\infty)] x^{\sigma + \delta_2} \inte{(x,\infty)}{z^{\delta_2 - 1 - \sigma}}{\Leb(z)} = - \frac{1}{\delta_2 - \sigma } 
	\end{align*}
	if $ \delta_2 < \sigma $, thus generating \[ \norm[\Lone((0,\infty))]{z^{\delta_2} \, w} \lesssim \norm[\Lone((0,\infty))]{z^{\delta_2} \, \overline{f}}, \] this time with a constant depending on $ \sigma $ and $ \delta_2 $. \\
	An interpolation yields \[ \norm[\Ltwo((0,\infty))]{z^{\frac{\delta_1 + \delta_2}{2}} \, w} \lesssim \norm[\Ltwo((0,\infty))]{z^{\frac{\delta_1 + \delta_2}{2}} \, \overline{f}} \] for $ \delta_1 < 1 + \sigma $ and $ \delta_2 < \sigma $. This condition allows us to choose $ \delta_1 $ and $ \delta_2 $ such that $ \delta_1 + \delta_2 = \sigma  $ if only $ \sigma < 1 + 2 \, \sigma $ or equivalently $ \sigma > - 1 $. Hence we get
	\begin{align*}
		\norm[\Ltwosig((0,\infty))]{\dz u} = \norm[\Ltwosig((0,\infty))]{w} \lesssim\norm[\Ltwosig((0,\infty))]{\overline{f}} \lesssim\norm[\Ltwosig((0,\infty))]{\widetilde{f}}, 
	\end{align*}
	where the constants here and in the following only depend on $ \sigma $. \\
	An immediate consequence of (\ref{ode2}) is then 
	\begin{align*}
		\norm[\Ltwo_{2 + \sigma}((0,\infty))]{\dz^2  u}	= \norm[\Ltwosig((0,\infty))]{z \, \dz^2  u} \eqsim \norm[\Ltwosig((0,\infty))]{\overline{f} + \dz u} \lesssim \norm[\Ltwosig((0,\infty))]{\widetilde{f}}.
	\end{align*}	
	Summing up, after reverting the notation back to the starting point,
	the retransformation from $ z $ to $ y_n $ and an integration in the $ \xi' $ direction combined with Plancherel's theorem in the reverse Fourier transformation reveals that 
	\[ \norm[\Ltwosig(H)]{\grady'  u} + \norm[\Ltwo_{2 + \sigma}(H)]{\lapy' u} + \norm[\Ltwosig(H)]{\dyn u} + \norm[\Ltwo_{2 + \sigma}(H)]{\dyn^2 u} \lesssim \norm[\Ltwosig(H)]{\widetilde{f}} \] with an additional dependence of the constant on $ n $.
	Finally, the mixed second order derivatives can be gained thanks to the formula \[ \norm[\Ltwo_{2 + \sigma}(H)]{D^2_y u}^2 \eqsim\norm[\Ltwo_{2 + \sigma}(H)]{\lapy u}^2 + \norm[\Ltwosig(H)]{\grady' u}^2 \] that holds up to a constant once again depending only on $ \sigma $, by means of integration by parts and the density statement from Remark \ref{bem_density}. An integration in time 
	finishes the proof.
\end{bew}

\begin{bem} \label{bem_LtwoCZO}
	It follows from Proposition \ref{prop_energyEst} that for $ l \geq 0 $, $ k \in \N_0 $ and $ \alpha \in \N_0^n $ with \[ (l,k,|\alpha|) \in \menge{(0,1,0), \, (0,0,1), \, (1,0,2)}, \] the mappings \[ f \mapsto y_n \, \ds^k \dy^\alpha u \] 
	that send the inhomogeneity to certain weighted derivatives of a $ \sigma $-solution on $ I \times \Hquer $ with zero initial value are bounded operators from $ \Ltwo(I; \, \Ltwosig(H)) $ to $ \Ltwo(I; \, \Ltwosig(H)) $. \\ 
	Note that the above set of exponents exactly contains those $ l \geq 0 $, $ k \in \N_0 $ and $ \alpha \in \N_0^n $ that satisfy the conditions $ l - k - |\alpha| = -1 $ and $ 2 \, l - |\alpha| \leq 0 $.
\end{bem}

\begin{bem} \label{bem_globalAux}
	A straight forward regularisation in the fashion of the proof of Proposition \ref{prop_eId} shows that for $ \sigma $-solutions with zero initial values there exists a constant $ c = c(n,\sigma) $ such that 
	\begin{align*}
		\inte{I}{\norm[\Ltwosigp(H)]{\grady \grady' u}^2}{\Leb} \leq c \, \left( \inte{I}{\norm[\Ltwosig(H)]{\grady' f}^2}{\Leb} + \inte{I}{\norm[\Ltwosig(H)]{f}^2}{\Leb} \right)
	\end{align*}
	if in addition to $ f \in \Ltwo \bigl(I; \, \Ltwosig(H)\bigr)  $ we also demand $ \grady' f \in \Ltwo \bigl(I; \, \Ltwosig(H)\bigr) $. 
\end{bem}

\begin{bem} \label{bem_iteratedSigmaSol}
	We can now iterate the notion of $ \sigma $-solution in terms of their derivatives rigorously: 
	Let $ k \in \N_0 $ and $ \alpha \in \N_0^n $. If for any $ 0 \leq j \leq k $, $ 0 \leq \beta' \leq \alpha' $, and $ 0 \leq |\gamma| \leq \alpha_n $ we have that \[ \ds^j \partial_{y'}^{\beta' + \gamma'} \dyn^{\gamma_n} f \in \Ltwo(I; \, \Ltwo_{\gamma_n + \sigma}(H)), \] and if $ u $ is a $ \sigma $-solution to $ f $ on $ I \times \Hquer $, then \[ \ds^k \dy^\alpha u \text{ is an } (\alpha_n + \sigma) \text{-solution to } \ds^k \dy^\alpha f + \alpha_n \, \lapy' \ds^k \partial_{y'}^{\alpha'} \dyn^{\alpha_n - 1} u \text{ on } (\widetilde{s}_1,s_2) \times \Hquer \] with natural initial value for any $ \widetilde{s}_1 \in I $. 
\end{bem}

\subsection{Local Estimates} \label{section_localEstimates}

For a localisation of these considerations, 
we shift our paradigm concerning the time interval and characterise it by its length $ r > 0 $ and its initial point $ s_1 \in \R $. We use the abbreviation \[ I_{r,\, \varepsilon}(s_1) \ldef (s_1 + \varepsilon \, r^2,s_1 + r^2) \] for $ \varepsilon \in [0,1) $, thus considering finite intervals that are possibly bounded away from their initial time, and drop the second index if $ \varepsilon = 0 $. \\
Note also that the behaviour of an intrinsic ball depends on the relative position of its centre point $ y_0 \in \Hquer $ with respect to $ \partial H $. For any point $ y \in B_r(y_0) $ we have that $ y_n \lesssim (r + \sqrt{\yzeron})^2 $, where the constant is absolute. Thus, locally, lowering weight exponents is possible on the expense of an extra factor. If $ r \leq \sqrt{\yzeron} $, we can also increase weight exponents by $ y_n \gtrsim (r + \sqrt{\yzeron})^2 $ for any $ y \in B_r(y_0) $, again without dependence on any parameter. \\

The fact that there exists a small $ \delta_0 \in (0,\delta_1) $ such that for any $ \delta_2 \in (0,\delta_0) $ a certain statement is true will be expressed by saying that it holds \enquote{for any $ \delta_2 \in (0,\delta_1) $ small enough}. 

\begin{prop} \label{prop_localEnergyEst}
	Let $ \varepsilon_1 \in [0,1) $ and $ \delta_1 \in (0,1] $. If $ f \in \Ltwo\bigl(I_{r,\varepsilon_1}(s_1); \, \Ltwosig(B_{\delta_1 r}(y_0))\bigr) $ and $ u $ is a $ \sigma $-solution to $ f $ on $ I_{r,\, \varepsilon_1}(s_1) \times B_{\delta_1 r}(y_0) $, then for any $ \varepsilon_2 \in (\varepsilon_1,1) $ as well as any $ \delta_2 \in (0,\delta_1) $ small enough there exists a constant $ c = c(n, \sigma, \varepsilon_1, \varepsilon_2, \delta_1, \delta_2) $ such that 
	\begin{align*}
		\inte{I_{r,\, \varepsilon_2}(s_1)}{\!\!\!\!\!\! \norm[\Ltwosig(B_{\delta_2 r}(y_0))]{\ds u}^2}{\Leb} &+ r^{-2} \, (r + \sqrt{\yzeron})^2 \!\!\!\!\!\! \inte{I_{r,\, \varepsilon_2}(s_1)}{\!\!\!\!\!\! \norm[\Ltwosig(B_{\delta_2 r}(y_0))]{\grady u}^2}{\Leb} + \!\!\!\!\!\! \inte{I_{r,\, \varepsilon_2}(s_1)}{\!\!\!\!\!\! \norm[\Ltwo_{2 + \sigma}(B_{\delta_2 r}(y_0))]{D^2_y u}^2}{\Leb} \\
		&\leq c \, \left( r^{-4} \!\!\!\!\!\! \inte{I_{r,\, \varepsilon_1}(s_1)}{\!\!\!\!\!\! \norm[\Ltwosig(B_{\delta_1 r}(y_0))]{u}^2}{\Leb} + \!\!\!\!\!\! \inte{I_{r,\, \varepsilon_1}(s_1)}{\!\!\!\!\!\! \norm[\Ltwosig(B_{\delta_1 r}(y_0))]{f}^2}{\Leb} \right).
	\end{align*}
\end{prop}
\begin{bew}
	In the whole proof, constants possibly depend on $ n, \sigma $ and any other parameters involved in the calculations, but not on $ r $ and $ \yzeron $. Fix $ \varepsilon_2 \in (\varepsilon_1,1) $ and $ \delta_2 \in (0,\delta_1) $. We can extend the concept of spatial cutoff function adapted to the intrinsic metric and include time to get \[ \eta \in \glatt_c(\Irechts_{r,\, \varepsilon_1}(s_1) \times B_{\delta_1 r}(y_0)) \] with \[ \eta = 1 \text{ on } \Iquer_{r,\, \varepsilon_2}(s_1) \times \overline{B}_{\delta_2 r}(y_0) \] and \[ |\ds^k \dy^\alpha \eta| \lesssim r^{-2k - |\alpha|} \, (r + \sqrt{\yzeron})^{-|\alpha|} \text{ on } \Irechts_{r,\, \varepsilon_1}(s_1) \times B_{\delta_1 r}(y_0) \] if $ \delta_2 $ is small enough. 
	It follows that $ \eta \, u $ is a $ \sigma $-solution on $ I_{\varepsilon_1,\, r}(s_1) \times \Hquer $ with zero initial value to the inhomogeneity $ F[u] + \ds \eta \, u $, where \[ F[u] \ldef \eta \, f - L_\sigma \eta \, u - 2 \, y_n \, \grady \eta \cdot \grady u. \] The calculation involves a spatial integration by parts that also works if the ball touches the boundary and thus both the test functions and $ \eta $ can have non-zero values at $ y_n = 0 $, since then $ y_n^{1 + \sigma} = 0 $ holds there. By Hölder's and Young's inequalities, we get \[ \ip[\Ltwosig(H)]{F[u](s)}{\eta(s) \, u(s)} \lesssim r^2 \, \norm[\Ltwosig(H)]{\eta(s) \, f(s)}^2 + r^{-2} \, \norm[\Ltwosig(H)]{\eta(s) \, u(s)}^2 + \norm[\Ltwosigp(H)]{\grady \eta(s) \, u(s)}^2 \] on $ \overline{I}_{r,\varepsilon_1}(s_1) $ with a constant that does not depend on any parameter. The properties of $ \eta $, $ f $ and $ u $ then show that \[ \inte{I_{r,\, \varepsilon_1}(s_1)}{\ip[\Ltwosig(H)]{F[u] + \ds \eta \, u}{\eta \, u}}{\Leb} \lesssim r^{-2} \inte{I_{r,\, \varepsilon_1}(s_1)}{\norm[\Ltwosig(B_{\delta_1 r}(y_0))]{u}^2}{\Leb} + r^2 \inte{I_{r,\, \varepsilon_1}(s_1)}{\norm[\Ltwosig(B_{\delta_1 r}(y_0))]{f}^2}{\Leb}, \] where we also lowered the weights appropriately. \\
	Along the same lines we see that 
	\begin{align*}
		&\inte{I_{r,\, \varepsilon_1}(s_1)}{\norm[\Ltwosig(H)]{F[u] + \ds \eta \, u}^2}{\Leb} \\
		&\lesssim r^{-4} \inte{I_{r,\, \varepsilon_1}(s_1)}{\norm[\Ltwosig(B_{\delta_1 r}(y_0))]{u}^2}{\Leb} + r^{-2} \inte{I_{r,\, \varepsilon_1}(s_1)}{\norm[\Ltwosigp(B_{\delta_1 r}(y_0))]{\grady u}^2}{\Leb} + \inte{I_{r,\, \varepsilon_1}(s_1)}{\norm[\Ltwosig(B_{\delta_1 r}(y_0))]{f}^2}{\Leb}. 
	\end{align*}
	It thus becomes clear that \[ F[u] + \ds \eta \, u \in \Ltwo\bigl(I_{r,\, \varepsilon_1}(s_1); \, \Ltwosig(B_{\delta_1 r}(y_0))\bigr) \subset \Lone\bigl(I_{r,\, \varepsilon_1}(s_1);\, \Ltwosig(B_{\delta_1 r}(y_0))\bigr). \] We can hence apply the global energy identity from Proposition \ref{prop_eId} on $ \eta \, u $ with zero initial value to obtain
	\begin{align*}
		\inte{I_{r,\, \varepsilon_2}(s_1)}{\norm[\Ltwosigp(B_{\delta_2 r}(y_0))]{\grady u}^2}{\Leb} &= \inte{I_{r,\, \varepsilon_2}(s_1)}{\norm[\Ltwosigp(B_{\delta_2 r}(y_0))]{\grady (\eta \, u)}^2}{\Leb} \\ 
		&\leq \inte{I_{r,\, \varepsilon_1}(\tau)}{\ip[\Ltwosig(H)]{F[u] + \ds \eta \, u}{\eta \, u}}{\Leb} \\
		&\lesssim r^{-2} \inte{I_{r,\, \varepsilon_1}(\tau)}{\norm[\Ltwosig(B_{\delta_1 r}(y_0))]{u}^2}{\Leb} + r^2 \inte{I_{r,\, \varepsilon_1}(\tau)}{\norm[\Ltwosig(B_{\delta_1 r}(y_0))]{f}^2}{\Leb}.
	\end{align*}
	But $ u $ is also a $ \sigma $-solution on the smaller set $ I_{r,\, \varepsilon_2}(s_1) \times B_{\delta_2 r}(y_0) $. For $ \widetilde{\varepsilon}_2 > \varepsilon_2 $ and $ \widetilde{\delta}_2 < \delta_2 $ small enough, we can choose a cutoff function $ \widetilde{\eta} $ and reiterate all the considerations from above to get a spatially global solution $ \widetilde{\eta} \, u $ with zero initial values and the inhomogeneity bound 
	\begin{align*}
		&\inte{I_{r,\, \varepsilon_2}(s_1)}{\norm[\Ltwosig(H)]{F[u] + \ds \eta \, u}^2}{\Leb} \\
		&\lesssim r^{-4} \inte{I_{r,\, \varepsilon_2}(s_1)}{\norm[\Ltwosig(B_{\delta_2 r}(y_0))]{u}^2}{\Leb} + r^{-2} \inte{I_{r,\, \varepsilon_2}(s_1)}{\norm[\Ltwosigp(B_{\delta_2 r}(y_0))]{\grady u}^2}{\Leb} + \inte{I_{r,\, \varepsilon_2}(s_1)}{\norm[\Ltwosig(B_{\delta_2 r}(y_0))]{f}^2}{\Leb} \\
		&\lesssim r^{-4} \inte{I_{r,\, \varepsilon_1}(s_1)}{\norm[\Ltwosig(B_{\delta_1 r}(y_0))]{u}^2}{\Leb} + \inte{I_{r,\, \varepsilon_1}(s_1)}{\norm[\Ltwosig(B_{\delta_1 r}(y_0))]{f}^2}{\Leb}
	\end{align*}
	thanks to the bound for the gradient proven before. By an application of the global energy estimates from Proposition \ref{prop_energyEst} on $ \widetilde{\eta} \, u $ we find
	\begin{align*}
		&\inte{I_{r,\, \widetilde{\varepsilon}_2}(s_1)}{\norm[\Ltwosig(B_{\widetilde{\delta}_2 r}(y_0))]{\ds u}^2}{\Leb} + \inte{I_{r, \, \widetilde{\varepsilon}_2}(s_1)}{\norm[\Ltwosig(B_{\widetilde{\delta}_2 r}(y_0))]{\grady u}^2}{\Leb} + \inte{I_{r,\, \widetilde{\varepsilon_2}}(s_1)}{\norm[\Ltwo_{2 + \sigma}(B_{\widetilde{\delta}_2 r}(y_0))]{D^2_y u}^2}{\Leb} \\
		&\quad \lesssim r^{-4} \inte{I_{r,\, \varepsilon_1}(s_1)}{\norm[\Ltwosig(B_{\delta_1 r}(y_0))]{u}^2}{\Leb} + \inte{I_{r,\, \varepsilon_1}(s_1)}{\norm[\Ltwosig(B_{\delta_1 r}(y_0))]{f}^2}{\Leb},
	\end{align*}
	while the energy identity yields \[ r^{-2} \inte{I_{r,\, \widetilde{\varepsilon}_2}(s_1)}{\norm[\Ltwosigp(B_{\widetilde{\delta}_2 r}(y_0))]{\grady u}^2}{\Leb} \lesssim r^{-4} \inte{I_{r,\, \varepsilon_1}(s_1)}{\norm[\Ltwosig(B_{\delta_1 r}(y_0))]{u}^2}{\Leb} + \inte{I_{r,\, \varepsilon_1}(s_1)}{\norm[\Ltwosig(B_{\delta_1 r}(y_0))]{f}^2}{\Leb}. \]
	Now remember that for $ r \leq \sqrt{\yzeron} $ we can increase weight exponents and hence get \[ r^{-2} \inte{I_{r,\, \widetilde{\varepsilon}_2}(s_1)}{\norm[\Ltwosigp(B_{\widetilde{\delta}_2 r}(y_0))]{\grady u}^2}{\Leb} \gtrsim r^{-2} (r + \sqrt{\yzeron})^2 \inte{I_{r,\, \widetilde{\varepsilon}_2}(s_1)}{\norm[\Ltwosig(B_{\widetilde{\delta}_2 r}(y_0))]{\grady u}^2}{\Leb}. \] In the converse case, however, we see that \[ \inte{I_{r,\, \widetilde{\varepsilon}_2}(s_1)}{\norm[\Ltwosig(B_{\widetilde{\delta}_2 r}(y_0))]{\grady u}^2}{\Leb} \gtrsim r^{-2} (r + \sqrt{\yzeron})^2 \inte{I_{r,\, \widetilde{\varepsilon}_2}(s_1)}{\norm[\Ltwosig(B_{\widetilde{\delta}_2 r}(y_0))]{\grady u}^2}{\Leb}, \] since $ r \geq \sqrt{\yzeron} $ implies $ 1 \gtrsim r^{-2} \, (r + \sqrt{\yzeron})^2 $ with a constant that once more does not depend on any parameter.
\end{bew}

\begin{bem} \label{bem_localAuxEst}
	The same methods of proof combined with an application of the local energy estimate from Proposition \ref{prop_localEnergyEst} also imply a local version of the global auxiliary estimate from Remark \ref{bem_globalAux}, namely
	\begin{align*}
		&\inte{I_{r,\, \varepsilon_2}(s_1)}{\norm[\Ltwosigp(B_{\delta_2 r}(y_0))]{\grady \grady' u}^2}{\Leb} \\
		&\lesssim (1 + r^{-2} \, (r + \sqrt{\yzeron})^{-2}) \, r^{-4} \inte{I_{r,\, \varepsilon_1}(s_1)}{\norm[\Ltwosig(B_{\delta_1 r}(y_0))]{u}^2}{\Leb} \\
		&\quad + (1 + r^{-2} \, (r + \sqrt{\yzeron})^{-2}) \inte{I_{r,\, \varepsilon_1}(s_1)}{\norm[\Ltwosig(B_{\delta_1 r}(y_0))]{f}^2}{\Leb} + \inte{I_{r, \, \varepsilon_1}(s_1)}{\norm[\Ltwosig(B_{\delta_1 r}(y_0))]{\grady' f}^2}{\Leb},
	\end{align*} 
	if only in addition we have $ \grady' f \in \Ltwo\bigl(I_{r,\, \varepsilon_1}(s_1); \, \Ltwosig(B_{\delta_1 r}(y_0))\bigr) $.
\end{bem}


We now iterate the local energy estimates, adjust the weights with a Hardy inequality and then adapt a Morrey-type inequality in time-space to our metric and measure to generate a pointwise estimate for arbitrary derivatives. This shows that local $ \sigma $-solutions are indeed smooth away from their initial time and spatial boundary. In favour of a clearer presentation we set $ f = 0 $. \\ 

\begin{prop} \label{prop_pointwiseEstimate}
	Let $ k \in \N_0 $ and $ \alpha \in \N_0^n $. If $ u $ is a $ \sigma $-solution to $ f = 0 $ on $ I_r(s_1) \times B_r(y_0) $, then for any $ \varepsilon \in (0,1) $ 
	and any $ \delta \in (0,1) $ small enough there exists a constant $ c = c(n,\sigma,k,\alpha,\varepsilon,\delta) > 0 $ such that \[ |\ds^k \dy^\alpha u(s,y)|^2 \leq \,  c \, r^{- 4 k -  2 |\alpha| -2} \, (r + \sqrt{\yzeron})^{- 2 |\alpha|} \, |B_r(y_0)|_\sigma^{-1} \inte{I_r(s_1)}{\norm[\Ltwosig(B_r(y_0))]{u}^2}{\Leb} \] for any $ (s,y) \in \Iquer_{r,\varepsilon}(s_1) \times B_{\delta r}(y_0) $. 
\end{prop}
\begin{bew}
	As before, constants in the inequalities in this proof depend on the parameters involved in the calculations, but not on $ r $ and $ \yzeron $. To iterate the local energy estimate, we use that $ u $ is also a $ \sigma $-solution on any smaller set and start with $ \delta_1 \in (0,1] $ and $ \varepsilon_1 \in [0,1) $ on the set $ I_{r,\varepsilon_1}(s_1) \times B_{\delta_1 r}(y_0) $ for arbitrary $ r > 0 $, $ s_1 \in \R $ and $ y_0 \in \Hquer $. We make the sets on the left hand side of the inequalities smaller by choosing $ \varepsilon_2 $ closer to $ 1 $ and $ \delta_2 $ smaller, and small enough, whenever this is necessary, but will not distinguish these sets in the notation and merely write $ I $ and $ B $ on both sides. \\ 

	We start with the tangential directions. For first order derivatives it is immediately clear by the local energy estimate from Proposition \ref{prop_localEnergyEst} and the auxiliary estimate from Remark \ref{bem_localAuxEst} that the right regularity for being a $ \sigma $-solution is given on the smaller sets for which those estimates hold. Therefore, $ \grady' u $ is a $ \sigma $-solution of the homogeneous equation to which both the local energy estimate and the local auxiliary estimate can be applied. Inductively, for any $ \alpha' \in \N_0^{n-1} $ we get that $ \dystrich^{\alpha'} u $ is a $ \sigma $-solution to $ 0 $ on $ I_{r,\varepsilon_2}(s_1) \times B_{\delta_2 r}(y_0) $ for $ \varepsilon_2 \in (\varepsilon_1,1) $ and $ \delta_2 \in (0,\delta_1) $ small enough, with the tangential energy estimate
	\begin{align*}
		\inte{I}{\norm[\Ltwosig(B)]{\ds \dystrich^{\alpha'} u}^2}{\Leb} &+ r^{-2} \, (r + \sqrt{\yzeron})^2 \inte{I}{\norm[\Ltwosig(B)]{\grady \dystrich^{\alpha'} u}^2}{\Leb} + \inte{I}{\norm[\Ltwo_{2 + \sigma}(B)]{D_y^2 \dystrich^{\alpha'} u}^2}{\Leb} \\
		&\lesssim r^{-4 - 2 |\alpha'|} \, (r + \sqrt{\yzeron})^{- 2 |\alpha'|} \inte{I}{\norm[\Ltwosig(B)]{u}^2}{\Leb}. 
	\end{align*}
	The iterated tangential auxiliary estimate carries an extra factor $ (1 + r^{-2} \, (r + \sqrt{\yzeron})^{-2}) $ on the right hand side. \\

	In the vertical direction, for a given $ \alpha_n \in \N $ we assume as an induction hypothesis that $ \dyn^{\alpha_n} u $ is an $ (\alpha_n + \sigma) $-solution to $ \dyn^{\alpha_n -1} \lapy' u $ on a smaller set, and that we have both the estimates 
	\begin{equation} \label{inequality_energy}
		\begin{split}
			\inte{I}{\norm[\Ltwo_{\alpha_n + \sigma}(B)]{\ds \dyn^{\alpha_n} u}^2}{\Leb} &+ r^{-2} \inte{I}{\norm[\Ltwo_{\alpha_n + \sigma}(B)]{\grady \dyn^{\alpha_n} u}^2}{\Leb} + \inte{I}{\norm[\Ltwo_{2 + \alpha_n + \sigma}(B)]{D_y^2 \dyn^{\alpha_n} u}^2}{\Leb} \\
			&\lesssim r^{-4} \, \left(1 + r^{-2} + r^{- 2} \, (r + \sqrt{\yzeron})^{-2} \right)^{\alpha_n} \inte{I}{\norm[\Ltwosig(B)]{u}^2}{\Leb}
		\end{split} \tag{$ E_{\alpha_n} $}
	\end{equation}
	and 
	\begin{equation} \label{inequality_aux}
		\begin{split}
			&\inte{I}{\norm[\Ltwo_{1 + \alpha_n + \sigma}(B)]{\grady \grady' \dyn^{\alpha_n} u}^2}{\Leb} \\
			&\quad \lesssim (1 + r^{-2} \, (r + \sqrt{\yzeron})^{-2}) \, r^{-4} \, \left(1 + r^{-2} + r^{- 2} \, (r + \sqrt{\yzeron})^{-2} \right)^{\alpha_n} \inte{I}{\norm[\Ltwosig(B)]{u}^2}{\Leb}.
		\end{split} \tag{$ A_{\alpha_n} $}
	\end{equation}
	The local energy estimate from Proposition \ref{prop_localEnergyEst} will be referred to as $ (E_0) $, and the local auxiliary estimate from Remark \ref{bem_localAuxEst} as $ (A_0) $.  In case of the base clause $ \alpha_n = 1 $, by a decrease of weight exponents and $ (E_0) $ it is clear that $ \dyn u $ has the regularity required for a $ (1 + \sigma) $-solution. $ (A_0) $ ensures that also the inhomogeneity $ \lapy' u $ is in the right space. 
	Therefore, in view of Remark \ref{bem_iteratedSigmaSol},  $ \dyn u $ is a $ (1 + \sigma) $-solution to $ \lapy' u $. We can thus apply $ (E_0) $ onto this solution with this inhomogeneity and this weight, and then use $ (E_0) $ onto $ u $ after decreasing the weight exponents once, as well as $ (A_0) $ onto $ u $, to get \eqref{inequality_energy} for $ \alpha_n = 1 $. Similarly, for \ref{inequality_aux} and $ \alpha_n = 1 $ the use of $ (A_0) $ on the newly found solution results in summands that can be controlled by $ (E_0) $ for $ u $ after lowering a weight exponent, as well as by $ (A_0) $ once for $ u $ and once for $ \grady' u $, which are both $ \sigma $-solutions to the homogeneous problem. \\ 
	Likewise, in the inductive step $ \alpha_n + 1 $ the right regularity follows from \eqref{inequality_energy} after lowering weights, and from \eqref{inequality_aux}. 
	Then $ \dyn^{\alpha_n +1} u $ is an $ (\alpha_n + 1 + \sigma) $-solution to $ \dyn^{\alpha_n} \lapy' u $. We apply $ (E_0) $ on this solution 
	and decrease weight exponents before using \eqref{inequality_energy} and \eqref{inequality_aux} for $ u $. For the iterated auxiliary estimate we first apply $ (A_0) $ onto $ \dyn^{\alpha_n +1} u $, then lower weight exponents and apply \eqref{inequality_energy} after another reduction of the weight exponent, and \eqref{inequality_aux} once for $ u $ itself and once for $ \grady' u $. This finishes the induction for the vertical direction. \\
	Restricting ourselves to $ r \leq 1 $ and $ \yzeron \leq 1 $, we gain \[ \left(1 + r^{-2} + r^{- 2} \, (r + \sqrt{\yzeron})^{-2} \right)^{\alpha_n} \lesssim r^{- 2 \alpha_n} \, (r + \sqrt{\yzeron})^{-2 \alpha_n} \] and thus recover the same shape of the estimate as before for the tangential case. \\

	Finally, we deal with the time direction now. Given both the estimates for the first order tangential and vertical derivatives, a simple induction as in the tangential case shows that $ \ds^k u $ is a $ \sigma $-solution to $ 0 $ and we have 
	\begin{align*}
		\inte{I}{\norm[\Ltwosig(B)]{\ds \ds^k u}^2}{\Leb} &+ r^{-2} \, (r + \sqrt{\yzeron})^2 \inte{I}{\norm[\Ltwosig(B)]{\grady \ds^k u}^2}{\Leb} + \inte{I}{\norm[\Ltwo_{2 + \sigma}(B)]{D_y^2 \ds^k u}^2}{\Leb} \\
		&\lesssim r^{-4 - 4 k} \inte{I}{\norm[\Ltwosig(B)]{u}^2}{\Leb}. 
	\end{align*}

	We can then view the case of general derivatives as $ \ds^k \dy^\alpha u = \dyn^{\alpha_n}(\dystrich^{\alpha'} (\ds^k u)) $ and apply first the $ y_n $-directional result, and then the other two estimates subsequently to obtain a combined version of the iterated energy estimates and thus
	\[ \inte{I}{\norm[\Ltwo_{\alpha_n + \sigma}(B)]{\ds^k \dy^\alpha u}^2}{\Leb} \lesssim r^{- 4 k - 2 |\alpha|} \, (r + \sqrt{\yzeron})^{-2 |\alpha|} \inte{I}{\norm[\Ltwosig(B)]{u}^2}{\Leb} \] for $ r $, $ \yzeron \leq 1 $. \\
	
	In the next step the weights on the left hand side are treated. We use Hardy's inequality $ \alpha_n $ times onto $ \eta \, \ds^k \dy^\alpha u $, where $ \eta $ is a purely spatial intrinsic cutoff function. After shrinking the spatial domain of integration on the left hand side once more, this leads to 
	\begin{align*}
		\inte{I}{\norm[\Ltwosig(B)]{\ds^k \dy^\alpha u}^2}{\Leb} &\lesssim \inte{I}{\norm[\Ltwo_{2 \alpha_n + \sigma}(B)]{\dyn^{\alpha_n}(\eta \, \ds^k \dy^\alpha  u)}^2}{\Leb} \\
		& \lesssim \summe{\gamma_n \leq \alpha_n} \inte{I}{\norm[\Ltwo_{2 \alpha_n + \sigma}(B)]{\dyn^{\alpha_n - \gamma_n} \eta \, \ds^k \dystrich^{\alpha'} \dyn^{\alpha_n + \gamma_n} u}^2}{\Leb}.
	\end{align*}
	The derivatives of $ \eta $ are bounded by $ r^{-(\alpha_n - \gamma_n)} \, (r + \sqrt{\yzeron})^{-(\alpha_n - \gamma_n)} $, entering squared here, with a constant depending on $ \alpha $ and the shrinking parameters only. But we also have $ \alpha_n = \alpha_n - \gamma_n + \gamma_n $ and can consequently lower the weight exponent $ y_n^{2 \alpha_n + \sigma} $ to $ y_n^{\alpha_n + \gamma_n + \sigma} $ on expense of a constant $ (r + \sqrt{\yzeron})^{2(\alpha_n - \gamma_n)} $. An application of the iterated energy estimate with additional weights for $ r $, $ \yzeron \leq 1 $ thus yields \[ \inte{I}{\norm[\Ltwosig(B)]{\ds^k \dy^\alpha u}^2}{\Leb} \lesssim r^{- 4k - 2 |\alpha| - 2 \alpha_n} \, (r + \sqrt{\yzeron})^{-2 |\alpha| - 2 \alpha_n} \inte{I}{\norm[\Ltwosig(B)]{u}^2}{\Leb}. \] 
	If in addition we assumed $ r \leq \sqrt{\yzeron} $, weights can also be increased directly, and we get the improved estimate \[ \inte{I}{\norm[\Ltwosig(B)]{\ds^k \dy^\alpha u}^2}{\Leb} \lesssim r^{- 4k - 2 |\alpha|} \, (r + \sqrt{\yzeron})^{-2 |\alpha| - 2 \alpha_n} \inte{I}{\norm[\Ltwosig(B)]{u}^2}{\Leb}. \]

	To get rid of the weights on the left hand side completely, we can use the same argument again, this time applying Hardy's inequality $ m_1 $ times, where $ m_1 \in \N_0 $ be the smallest integer with $ 2 \, m_1 \geq \sigma $. Then we have for any $ r $, $ \yzeron \leq 1 $ that \[ \inte{I}{\norm[\Ltwo(B)]{\ds^k \dy^\alpha u}^2}{\Leb} \lesssim r^{- 4k - 2 |\alpha| - 2 \alpha_n - 4 m_1} \, (r + \sqrt{\yzeron})^{-2 |\alpha| - 2 \alpha_n - 2 \sigma} \inte{I}{\norm[\Ltwosig(B)]{u}^2}{\Leb}, \] while for the region where in addition $ r \leq \sqrt{\yzeron} $ we get \[ \inte{I}{\norm[\Ltwo(B)]{\ds^k \dy^\alpha u}^2}{\Leb} \lesssim r^{- 4k - 2 |\alpha|} \, (r + \sqrt{\yzeron})^{- 2 |\alpha| - 2 \alpha_n - 2 \sigma} \inte{I}{\norm[\Ltwosig(B)]{u}^2}{\Leb}. \]
	We now apply a Morrey-type inequality adapted to the intrinsic metric. For $ m_2 \in \N $ with $ m_2 > \frac{n+1}{2} $ and any $ (s,y) $ in a once again smaller set, 
	this reads
	\begin{align*}
		|\ds^k \dy^\alpha u(s,y)|^2 &\lesssim \summe{j + |\gamma|\leq m_2} r^{4 j + 2 |\gamma| - n - 2} \, (r + \sqrt{\yzeron})^{2 |\gamma| - n} \, \inte{I}{\norm[\Ltwo(B)]{\ds^{k+j} \dy^{\alpha + \gamma} u}^2}{\Leb} \\
		&\lesssim r^{- 4 k - 2 |\alpha| - 2 \alpha_n - 4 m_1 - 2 m_2 - n - 2} \, (r + \sqrt{\yzeron})^{- 2 |\alpha| - 2 \alpha_n - 2 \sigma - 2 m_2 - n} \inte{I}{\norm[\Ltwosig(B)]{u}^2}{\Leb},
	\end{align*}
	where we applied the iterated energy estimates with unweighted left hand side for $ r $, $ \yzeron \leq 1 $. In the improved version we gain \[ |\ds^k \dy^\alpha u(s,y)|^2 \lesssim r^{- 4k - 2 |\alpha| - n - 2} \, (r + \sqrt{\yzeron})^{- 2 |\alpha| - 2 \alpha_n - 2 \sigma - 2 m_2 - n} \inte{I}{\norm[\Ltwosig(B)]{u}^2}{\Leb}, \] if in addition $ r \leq \sqrt{\yzeron} $. \\
	Note that if either $ r = 1 $ or $ \yzeron = 1 $, any term containing $ (r + \sqrt{\yzeron}) $ is bounded by an absolute constant regardless of the exponent it carries, and thus drops out. Therefore, both for $ r \leq 1 = \sqrt{\yzeron} $ and $ \sqrt{\yzeron} \leq 1 = r $ we finally get \[ |\ds^k \dy^\alpha u(s,y)|^2 \lesssim r^{-4 k - 2 |\alpha| - n - 2} \inte{I_{r,\varepsilon_1}(s_1)}{\norm[\Ltwosig(B_{\delta_1 r }(y_0))]{u}^2}{\Leb} \] for any $ (s,y) \in I_{r,\varepsilon_2}(s_1) \times B_{\delta_2 r}(y_0) $, with $ \varepsilon_2 \in (\varepsilon_1,1) $ and $ \delta_2 \in (0,\delta_1) $ small enough. \\

	We now use that the equation is invariant under the scaling $ A_\lambda $ and adjust the sets to the intrinsic setting: If $ u $ is a $ \sigma $-solution on $ I_r(s_1) \times B_r(y_0) $, then $ u \circ A_\lambda $ is a $ \sigma $-solution on $ I_{\hat{r},\varepsilon_1}(\hat{s}_0) \times B_{\delta_1 \hat{r}}(\hat{y}_0) $ with $ \hat{r} = \frac{r}{\sqrt{\lambda}} $, $ \hat{s}_0 = \frac{s_1}{\lambda} $ and $ \hat{y}_0 = \frac{1}{\lambda} y_0 $ for a $ \delta_1 < 1 $ small enough and $ \varepsilon_1 = 0 $. An application of the above result onto $ u \circ A_\lambda $ in conjunction with the integral transformation formula then reveals that for any $ (s,y) \in I_{r,\varepsilon_2}(s_1) \times B_{\delta_2 r}(y_0) $ with $ \varepsilon_2 \in (0,1) $ and $ \delta_2 < \delta_1 $ small enough we get \[ |\ds^k \dy^\alpha u(s,y)|^2 \lesssim \lambda^{- 2 k - 2 |\alpha| - \sigma - n - 1} \, \hat{r}^{- 4 k - 2 |\alpha| - n - 2} \inte{I_r(s_1)}{\norm[\Ltwosig(B_r(y_0))]{u}^2}{\Leb} \] if only $ \hat{r} $, $ \hat{y}_{0,n} \leq 1 $ and either one of them equals $ 1 $. However, the last condition is satisfied for any $ r $ and $ y_0 $ if we choose $ \lambda = r^2 $ in case $ \sqrt{\yzeron} \leq r $, and $ \lambda = \yzeron $ in case $ r \leq \sqrt{\yzeron} $. In both cases, by means of the formula for the $ \sigma $-measure of an intrinsic ball the resulting factor can be estimated to the one stated.
\end{bew}

\subsection{Estimates Against Initial Values} \label{section_estInitial}

We use the local regularity to generate pointwise exponential estimates for global $ \sigma $-solutions on $ I \times \Hquer $. 
To this end, consider the function \[ \Psi(y; \constpsi,\epsi,z_\Psi) \ldef \constpsi \frac{\quasi(y,z_\Psi)^2}{\sqrt{\epsi + \quasi(y,z_\Psi)^2}} \text{ for any } y \in \Hquer, \] where $ \constpsi \in \R $, $ \epsi > 0 $ are constants and $ z_\Psi \in \Hquer $ is arbitrary, but fixed. 

\begin{bem} \label{bem_psi}
	A straightforward calculation shows that $ \Psi(\, \cdot \, ;\constpsi,\epsi,z_\Psi) \in \stetig^1(H) $ and \[ \sqrt{y_n} \, |\grady \Psi(y; \constpsi,\epsi,z_\Psi)| \leq \constl \, |\constpsi| \text{ for any } y \in \Hquer \] with $ \constl \ldef 2^6 $. This implies that $ \Psi $ is a Lipschitz function in terms of the intrinsic metric $ \dist $ with the estimate \[ |\Psi(y;\constpsi,\epsi,z_\Psi) - \Psi(x;\constpsi,\epsi,z_\Psi)| \leq \constl \, |\constpsi| \, \dist(y,x). \] 
\end{bem}

For the proof of the following exponential $ \Linfty $-$ \Ltwo $-estimate, we first extend the norm decrease of global $ \sigma $-solutions of the homogeneous equation mentioned in Remark \ref{bem_decrease} by an exponential factor involving $ \Psi $. 

\begin{prop} \label{prop_LtwoLinftyExpo}
	If $ u $ is a $ \sigma $-solution to $ f = 0 $ on $ I \times \Hquer $ with initial value $ u_0 \in \Ltwosig(H) $, then there exists a constant $ c = c(n,\sigma,k,\alpha) > 0 $ such that 
	\begin{align*}
		&|\ds^k \dy^\alpha u(s,y)| \\
		&\quad\leq c \, r(s)^{-2 k - |\alpha|} \, (r(s) + \sqrt{y_n})^{-|\alpha|} \, |B_{r(s)}(y)|_\sigma^{-\frac{1}{2}} \, e^{2 \constl^2 \constpsi^2 r(s)^2 - \Psi(y;\constpsi,\epsi,z_\Psi)} \, \norm[\Ltwosig(H)]{e^{\Psi(\, \cdot \,;\constpsi,\epsi,z_\Psi)} \, u_0} 
	\end{align*}
	 for any $ y \in \Hquer $ and any $ s \in I $, where $ r(s) \ldef \sqrt{s - s_1} $.
\end{prop}
\begin{bew} 
	We abbreviate $ \Psi(\, \cdot \,;\constpsi,\epsi,z_\Psi) \rdef \Psi $. For $ s \in I $ define \[ F(s) \ldef e^{-2 \constl^2 \constpsi^2 (s-s_1)} \norm[\Ltwosig(H)]{e^{\Psi} \, u(s)}^2 + 2 \inte{(s_1,s)}{e^{- 2 \constl^2 \constpsi^2 (\varsigma-s_1)} \norm[\Ltwosigp(H)]{\grady(e^{\Psi} \, u(\varsigma))}^2}{\Leb(\varsigma)}. \] Then obviously we have 
	\begin{align*}
		\ds F(s) = e^{-2 \constl^2 \constpsi^2 (s-s_1)} \, \ds\left(\norm[\Ltwosig(H)]{e^{\Psi} \, u(s)}^2\right) &- 2 \, \constl^2 \, \constpsi^2 \, e^{-2 \constl^2 \constpsi^2 (s-s_1)} \, \norm[\Ltwosig(H)]{e^{\Psi} \, u(s)}^2 \\
		&+ 2 \, e^{-2 \constl^2 \constpsi^2 (s-s_1)} \, \norm[\Ltwosigp(H)]{\grady(e^{\Psi} \, u(s))}^2.
	\end{align*}	
	Note that \[ \ds\left(\norm[\Ltwosig(H)]{e^{\Psi} \, u(s)}^2\right) = 2 \inte{H}{e^{2 \Psi} \, u(s) \, \ds u(s)}{\musig}, \] where the differentiation under the integral is justified by the product rule for bilinear forms. Now, for a fixed time $ s $ we use $ e^{2 \Psi} \, u(s) $ formally as a test function for the $ \sigma $-solution $ u(s) $ in the equivalent formulation for $ \sigma $-solutions with regular temporal derivative described in the proof of Proposition \ref{prop_energyEst}. The formal consideration can be justified rigorously by a bounded approximation of $ \Psi $. 
	A repeated application of the product rule for derivatives results in 
	\begin{align*} 
		\inte{H}{e^{2 \Psi} \, u(s) \, \ds u(s)}{\musig} &= - \norm[\Ltwosigp(H)]{\grady(e^{\Psi} \, u(s))}^2 + \norm[\Ltwosigp(H)]{e^{\Psi} \, u(s) \, \grady \Psi}^2 \\
		&\leq - \norm[\Ltwosigp(H)]{\grady(e^{\Psi} \, u(s))}^2 + \constl^2 \, \constpsi^2 \, \norm[\Ltwosig(H)]{e^{\Psi} \, u(s)}^2,
	\end{align*}
	where the use of the weighted bound for $\grady \Psi $ from Remark \ref{bem_psi} is crucial to reduce the weight exponent in the last step. \\
	This shows that $ \ds F(s) \leq 0 $ on $ I $. Thus $ F $ is monotonically decreasing and we have $ F(s) \leq F(s_1) $, proving that \[ \norm[\Ltwosig(H)]{e^{\Psi} \, u(s)} \leq e^{\constl^2 \constpsi^2 (s - s_1)} \, \norm[\Ltwosig(H)]{e^{\Psi} \, u_0} \] for all $ s \in \Iquer $. \\

	Now our global $ \sigma $-solution $ u $ is also a $ \sigma $-solution to the initial value problem of the homogeneous equation on $ (s_1,s_0) \times B_r(y_0) $ for any  time $ s_0 \in I $ and any point $ y_0 \in \Hquer $ combined with any radius $ r > 0 $. For $ r \ldef r(s_0) \ldef \sqrt{s_0 - s_1} $ we have $ I_r(s_1) = (s_1,s_0) $ and can use the pointwise estimate from Proposition \ref{prop_pointwiseEstimate} in the temporal end point $ s_0 $ and the spatial centre point $ y_0 $ to obtain
	\begin{align*}
		|\ds^k \dy^\alpha u(s_0,y_0)| \lesssim r^{-2 k - |\alpha| -1} \, (r + \sqrt{\yzeron})^{-|\alpha|} \, |B_r(y_0)|_\sigma^{-\frac{1}{2}} \supre[y \in B_r(y_0)] e^{-\Psi(y)} \, r \, \supre[s \in (s_1,s_0)] \norm[\Ltwosig(B_r(y_0))]{e^{\Psi} \, u(s)}. 
	\end{align*}
	The monotone decrease of the exponential $ \Ltwosig $-norm just shown then implies 
	\begin{align*}
		e^{\Psi(y_0)} \, |\ds^k \dy^\alpha u(s_0,y_0)| \lesssim r^{-2 k - |\alpha|} \, (r + \sqrt{\yzeron})^{-|\alpha|} \, |B_r(y_0)|_\sigma^{-\frac{1}{2}} \, \supre[y \in B_r(y_0)] e^{\Psi(y_0) - \Psi(y)} e^{\constl^2 \constpsi^2 r^2} \, \norm[\Ltwosig(H)]{e^{\Psi} \, u_0}
	\end{align*}
	and we can use the Lipschitz property of $ \Psi $ from Remark \ref{bem_psi} and a short computation on the exponent of the radius factor to finish the proof. 
\end{bew}

The exponential factor allows us to gain an estimate in terms of the space $ X(p) $ by rather rough norms of the initial value, proving the first part of Theorem \ref{theo_linear}.

\begin{prop} \label{prop_roughInitial} 
	If $ u $ is a $ \sigma $-solution to $ f = 0 $ on $ (0,S) \times \Hquer $ with initial value $ u_0 \in \Ltwosig(H) $, then there exists a constant $ c = c(n,\sigma) > 0 $ such that \[ \norm[X(p)]{u} \leq c \, \norm[\Linfty(H)]{\grady u_0} \] for any $ 1 \leq p < \infty $.
\end{prop}
\begin{bew}
	Fix a time $ s_0 \in (0,S) $ as well as a point $ y_0 \in \Hquer $. For any constant $ C > 0 $, with $ u $ also $ u - C $ is a $ \sigma $-solution to $ f = 0 $ on $ (0,S) \times \Hquer $, but with initial value $ u_0 - C $, and an application of Proposition \ref{prop_LtwoLinftyExpo} with $ C \ldef u_0(y_0) $ yields 
	\begin{align*}
		\Bigg| \ds^k \dy^\alpha &\Bigl(u(s,y) - u_0(y_0)\Bigr)\Big|_{(s,y)=(s_0,y_0)} \Bigg| \\
		&\lesssim \frac{e^{2 \constl^2 \constpsi^2 s_0 - \Psi(y_0; \constpsi,\epsi,z_\Psi)}}{\sqrt{s_0}^{2 k + |\alpha|} \, (\sqrt{s_0} + \sqrt{\yzeron})^{|\alpha|} \, |B_{\sqrt{s_0}}(y_0)|_\sigma^\frac{1}{2}} \, \norm[\Ltwosig(H)]{e^{\Psi(\, \cdot \,; \constpsi,\epsi,z_\Psi)} \; |\cdot \, - y_0| \, } \, \norm[\Linfty(H)]{\grady u_0},
	\end{align*}
	since \[ |u_0(y) - u_0(y_0)| \leq |y - y_0| \, \norm[\Linfty(H)]{\grady u_0} \] by the fundamental theorem of calculus. \\ 
	To treat the integral we fix a radius $ r > 0 $ and view the upper half plane $ H $ as a disjoint union of annular rings $ B_{jr}(y_0) \ohne \overline{B}_{(j-1)r}(y_0) \rdef R_{r_j}(y_0) $ for $ j \in \N $. We now specify $ \Psi(\, \cdot \,;\constpsi,\epsi,z_\Psi) $ by setting $ \constpsi \ldef - \frac{1}{r} $, choosing $ \epsi $ such that $ \epsi < r^2 $ and letting $ y_0 $ play the role of the parameter point $ z_\Psi $. For $ y \in R_{r_j}(y_0) $ we then have \[ \Psi(y;\constpsi,\epsi,z_\Psi) \leq \frac{-j + 2}{8 \constd^2} \] thanks to an application of the equivalence estimate for the quasi-metric $ \quasi $. This amounts to
	\begin{align*}
		\inte{R_{r_j}(y_0)}{e^{2 \Psi(y;\constpsi,\epsi,y_0)} \, |y - y_0|^2}{\musig(y)} \lesssim j^{4 + n} \, (1 + j^{n + 2 \sigma}) \, e^{- \frac{1}{4 \constd^2} j} \, r^2 \, (r + \sqrt{\yzeron})^2 \, |B_r(y_0)|_\sigma 
	\end{align*}
	by the doubling property of the metric and the relation between intrinsic and Euclidean balls. \\ 
	Summing this over $ j $ and setting $ r\ldef \sqrt{s_0} $, the special choice of $ \Psi $ gives the bound 
	\begin{align*}
		\Big|\ds^k \dy^\alpha u(s_0,y_0) \Big| \lesssim \sqrt{s_0}^{-2 k - |\alpha| + 1} \, (\sqrt{s_0} + \sqrt{\yzeron})^{-|\alpha| + 1} \, \norm[\Linfty(H)]{\grady u_0}
	\end{align*}
 	for $ k + |\alpha| > 0 $ 
	and any $ s_0 \in (0,S) $ as well as $ y_0 \in \Hquer $, 
	where the convergent series over $ j $ 
	is subsumed into the constant. For $ k = 0 $ and $ |\alpha| = 1 $ this estimates the pointwise component of the $ X^{(1)}(p) $- norm of $ u $, while for $ k = 0 $ and $ |\alpha| = 2 $ with an additional multiplication by $ \sqrt{y_n} \leq (\sqrt{y_n} + \sqrt{s}) $ we bound the pointwise component of the $ X^{(2)} $-norm, both by the homogeneous Lipschitz norm of $ u_0 $ on $ \Hquer $. \\
	 
	Now fix $ p \geq 1 $ and $ z \in \Hquer $ as well as $ 0 < r^2 \leq S $. Applying the pointwise bound just obtained for $ (s,y) \in Q_r(z) $ after multiplying both sides by $ y_n^l $ 
	yields 
	\begin{align*}
		 \norm[\Lp(Q_r(z))]{\Pnl \, \ds^k \dy^\alpha u} &\lesssim |Q_r(z)|^\frac{1}{p} \, \supre[(s,y) \in Q_r(z)] \left( \sqrt{s \,}^{-2 k - |\alpha| + 1} \, (\sqrt{s} + \sqrt{y_n})^{-|\alpha| + 1 + 2l}
		\right) \, \norm[\Linfty(H)]{\grady u_0} \\
		&\lesssim  |Q_r(z)|^\frac{1}{p} \, r^{-2k -|\alpha| + 1} \, (r + \sqrt{z_n})^{- |\alpha| + 1 + 2l} \, \norm[\Linfty(H)]{\grady u_0},
	\end{align*}
	since the cylinders are bounded away from the initial time. 
	It is now possible to choose the right orders of derivatives and exponents to fit the definition of $ X(p) $ and thus finish the proof.	
\end{bew}
 
\subsection{Gaussian Estimate} \label{section_gaussianEstimate}

On any arbitrary time interval $ I = (s_1,s_2) $, the results shown so far, namely Propositions \ref{prop_eId} and \ref{prop_LtwoLinftyExpo}, provide us with the existence of a Green function $ G_\sigma $ that characterises any $ \sigma $-solution $ u $ to $ f = 0 $ on $ I \times \Hquer $ with initial value $ u_0 \in \Ltwosig(H) $ by the representation \[ \ds^k \dy^\alpha u(s,y) = \inte{H}{\ds^k \dy^\alpha G_\sigma(s,y,s_1,z) \, u_0(z)}{\Lebn(z)}. \] 
More than that, we will even show in the following that $ G_\sigma $ decays exponentially with a bound that takes on the shape of a Gaussian function with respect to the weighted measure $ \musig $ and the intrinsic metric $ \dist $. Such an estimate is called Gaussian estimate or Aronson-type estimate after one of the first authors exploring this type of inequalities \citep{aronson_gaussian}). For general uniformely strongly parabolic equations their proof was originally given by means of the Harnack inequality contained in \citet{moser} and \citet{mosererr}. This order was reversed by \citet{fabesstroock} and Gaussian estimates were shown directly. The idea was extended by \citet{koch_habil} to cover the degenerate parabolic case with measurable coefficients. Our proof simplifies this approach in a special case of constant coefficients and at the same time adds control over the derivatives of the Green function. 

\begin{prop} \label{prop_gaussian}
	Let $ k, \, j \in \N_0 $ and $ \alpha, \,  \beta \in \N_0^n $. Then there exist constants $ c = c(n,\sigma,k,j,\alpha,\beta) > 0 $ and $ C > 0 $ such that 
	\begin{align*}
		\Bigg| \dta^j \dz^\beta \Bigl( z_n^{- \sigma} \, &\ds^k \dy^\alpha G_\sigma(s,y,\tau,z) \Bigr) \Bigg| \\
		&\leq c \, \sqrt{s - \tau}^{-2 k - 2 j - |\alpha| - |\beta|} \, (\sqrt{s - \tau} + \sqrt{y_n})^{- |\alpha| - |\beta|} \, |B_{\sqrt{s - \tau}}(z)|_\sigma^{-1} \, e^{-\frac{\dist(y,z)^2}{C (s-\tau)}}
	\end{align*}
	for all $ y \neq z \in \Hquer $ and all $ \tau < s \in \Iquer $, and with any possible combination of the points $ y $ and $ z $ in the factors $ (\sqrt{s - \tau} + \sqrt{y_n})^{- |\alpha| - |\beta|} \, |B_{\sqrt{s - \tau}}(z)|_\sigma^{-1} $.
\end{prop}
\begin{bew}
	Fix $ \tau_0 < s_0 \in \Iquer $ as well as $ y_0 \in \Hquer $, and denote $ r \ldef \sqrt{s_0-\tau_0} $ and $ I_0 \ldef (\tau_0,s_0) \subset I $. 
	Define the linear subspaces \[ U_\Psi \ldef \menge{ h \in \Ltwosig(H) \mid \exists \, g \in \Ltwosig(H) \text{ with } h = e^{\Psi(\, \cdot \,;\constpsi,\epsi,z_\Psi)} \, g} \subset \Ltwosig(H) \] and \[ V_\Psi \ldef \menge{ \xi \in \Lonesig(H) \mid e^{\Psi(\, \cdot \,;\constpsi,\epsi,z_\Psi)} \, \xi \in \Ltwosig(H)} \subset \Lonesig(H). \] 
	We define a modified solution operator \[ A: U_\Psi \ni e^{\Psi(\, \cdot \,;\constpsi,\epsi,z_\Psi)} \, g \mapsto e^{\Psi(\, \cdot \,;\constpsi,\epsi,z_\Psi)} \, u\left(\frac{s_0+\tau_0}{2}\right), \] where $ u $ is the $ \sigma $-solution on $ I_0 \times \Hquer $ with initial value $ g $. The multiplication operator on $ \Ltwosig(H) $ with respect to the function $ \left( y \mapsto |B_r(y)|_\sigma^\frac{1}{2} \right) \in \Linfty(H) $ is denoted by $ M $. Then Proposition \ref{prop_LtwoLinftyExpo} at time $ \frac{s_0 + \tau_0}{2} \in \overline{I}_0 $ can be rephrased to \[ \norm[\Linfty(H)]{(M \circ A) (e^{\Psi(\, \cdot \,;\constpsi,\epsi,z_\Psi)} g)} \lesssim e^{\constl^2 \constpsi^2 r^2} \, \norm[\Ltwosig(H)]{e^{\Psi(\, \cdot \,;\constpsi,\epsi,z_\Psi)} g}. \] Thus $ M \circ A $ maps into $ \Linfty(H) $ and has operator norm bounded by $ e^{\constl^2 \constpsi^2 r^2} $ times a constant depending on $ n $ and $ \sigma $ only. Consequently, the dual operator maps $ \Lonesig(H) $ into $ \Ltwosig(H) $ with the same operator norm, so for any $ \xi \in V_\Psi \subset \Lonesig(H) $ we have \[ \norm[\Ltwosig(H)]{\left(M \circ A \right)^* \xi} \lesssim e^{\constl^2 \constpsi^2 r^2} \, \norm[\Lonesig(H)]{\xi}. \] The duality equation from Remark \ref{bem_duality} on the interval with starting point $ \frac{s_0 + t_0}{2} $ implies that \[ A^* \circ M^*: \xi \mapsto \, e^{-{\Psi(\, \cdot \,;\constpsi,\epsi,z_\Psi)}} \, v\left(\frac{s_0 + \tau_0}{2}\right), \] where $ v $ is the $ \sigma $-solution on $ I_0 \times \Hquer $ with initial value $ |B_r( \, \cdot \,)|_\sigma^\frac{1}{2} \,  e^{\Psi(\, \cdot \,;\constpsi,\epsi,z_\Psi)} \, \xi $. \\
	Applying the pointwise exponential estimate from Proposition \ref{prop_LtwoLinftyExpo} onto $ v $ for the time $ s_0 \in \overline{I}_0 $ and with the Lipschitz function \[ \widetilde{\Psi}(\, \cdot \,; \zeta_{\widetilde{\Psi}},\varepsilon_{\widetilde{\Psi}},z_{\widetilde{\Psi}}) \ldef \Psi(\, \cdot \,; - \constpsi,\epsi,z_\Psi) = - \Psi(\, \cdot \,; \constpsi,\epsi,z_\Psi) \] gives 
	\begin{align*}
		&|\ds^k \dy^\alpha v(s_0,y_0)| \\
		&\quad \lesssim r^{- 2 k - |\alpha|} \, (r + \sqrt{\yzeron})^{-|\alpha|} \, |B_r(y_0)|_\sigma^{-\frac{1}{2}} \, e^{\constl^2 \constpsi^2 r^2 + \Psi(y_0;\constpsi,\epsi,z_\Psi)} \, \norm[\Ltwosig(H)]{e^{-\Psi(\, \cdot \,;\constpsi,\epsi,z_\Psi)} \, v\left(\frac{s_0 + \tau_0}{2}\right)}.
	\end{align*}
	The combination with the duality inequality from above then yields \[ |\ds^k \dy^\alpha v(s_0,y_0)| \lesssim r^{- 2 k - |\alpha|} \, (r + \sqrt{\yzeron})^{-|\alpha|} \, |B_r(y_0)|_\sigma^{-\frac{1}{2}} \, e^{2 \constl^2 \constpsi^2 r^2 + \Psi(y_0;\constpsi,\epsi,z_\Psi)} \norm[\Lonesig(H)]{\xi}. \] 
	Since $ \Linfty(H) $ is the dual space to $ \Lonesig(H) $, for almost any $ z \in \Hquer $ we see
	\begin{align*}
		& z_n^{-\sigma} \, e^{\Psi(z;\constpsi,\epsi,z_\Psi)} \, |B_r(z)|_\sigma^\frac{1}{2} \, |\ds^k \dy^\alpha G_\sigma(s_0,y_0,\tau_0,z)| \\
		&\qquad \leq \supre[{\norm[\Lonesig(H)]{\xi} \leq 1}] \left| \inte{H}{e^{\Psi(\, \cdot \,;\constpsi,\epsi,z_\Psi)} \, |B_r(\, \cdot)|_\sigma^\frac{1}{2} \, \ds^k \dy^\alpha G_\sigma(s_0,y_0,\tau_0, \, \cdot \,) \, \xi}{\musig} \right| \\
		&\qquad = \supre |\ds^k \dy^\alpha v(s_0,y_0)|,
	\end{align*}
	where the last supremum is taken over any $ \sigma $-solution $ v $ on $ I_0 \times \Hquer $ with initial value $ |B_r(\, \cdot \,)|_\sigma^\frac{1}{2} \, e^{\Psi(\, \cdot \,;\constpsi,\epsi,z_\Psi)} \, \xi $ and $ \norm[\Lonesig(H)]{\xi} \leq 1 $. We have thus shown 
	\begin{align*}
		&z_n^{-\sigma} \, |\ds^k \dy^\alpha G_\sigma(s_0,y_0,\tau_0,z)| \\
		&\qquad \lesssim r^{-2 k - |\alpha |} \, (r + \sqrt{\yzeron})^{-|\alpha|} \, |B_r(y_0)|_\sigma^{-\frac{1}{2}} \, |B_r(z)|_\sigma^{-\frac{1}{2}} \, e^{\Psi(y_0;\constpsi,\epsi,z_\Psi)-\Psi(z;\constpsi,\epsi,z_\Psi) + 2 \constl^2 \constpsi^2 r^2 }  
	\end{align*}
	for almost all $ z \in \Hquer $. \\
	Now fix $ z = z_0 \neq y_0 $, and specify $ \Psi(\, \cdot \,;\constpsi,\epsi,z_\Psi) $ by setting $ \epsi\ldef 3 \, \quasi(y_0,z_0) $, $ z_\Psi \ldef z_0 $ and $ \constpsi \ldef - c $ for an arbitrary constant $ c > 0  $. 
	In conjunction with the equivalence characterisation of the intrinsic metric we get \[| \zzeron^{-\sigma} \, \ds^k \dy^\alpha G_\sigma(s_0,y_0,\tau_0,z_0)| \lesssim r^{- 2 k - |\alpha|} \, (r + \sqrt{\yzeron})^{-|\alpha|} \, |B_r(y_0)|_\sigma^{-\frac{1}{2}} \, |B_r(z_0)|_\sigma^{-\frac{1}{2}} \, e^{- \frac{c}{2 \constd} \dist(y_0,z_0) + 2 \constl^2 c^2 r^2}. \] Finally, we optimise over the constant $ c $ to see that the exponent takes on a minimum for \[ c_* = \frac{\dist(y_0,z_0)}{8 \constl^2 \constd r^2}. \] Inserting this value into the inequality then gives the Gaussian estimate \[ |\ds^k \dy^\alpha G_\sigma(s_0,y_0,\tau_0,z_0)| \lesssim r^{- 2 k - |\alpha|} \, (r + \sqrt{\yzeron})^{-|\alpha|} \, |B_r(y_0)|_\sigma^{-\frac{1}{2}} \, |B_r(z_0)|_\sigma^{-\frac{1}{2}} \, \zzeron^\sigma \, e^{-\frac{\dist(y_0,z_0)^2}{C r^2}} \] with $ C = 32 \, \constd^2 \, \constl^2 $. \\

	We can replace the measure of the ball centred at $ z_0 $ or the one centred at $ y_0 $ by the mutually other one without loosing more than a factor $ (1 + \frac{\dist(y_0,z_0)}{r})^{n + \sigma} $. A similar remark applies for $ (r + \sqrt{\yzeron})^{-|\alpha|} $. The exponential decay makes the loss in this exchange controllable at the expense of a portion of the decay expressed by the constant in the exponent. Hence, from now on we do not specify this constant any more and merely write $ C $ even if it changes in the course of the argument. \\

	Note now that there exists a constant $ \delta > 0 $ with $ \delta^2 < \frac{\tau_0-s_1}{r^2} $. We set $ \widetilde{r}\ldef \delta \, r $ and $ \widetilde{\tau}_0 \ldef (1 + \delta ^2) \, \tau_0 - \delta^2 \, s_0 > s_1 $ and obtain $ I_{\widetilde{r}}(\widetilde{\tau}_0) = (\widetilde{\tau}_0, \tau_0) \subset I $. 
	The duality identity from Remark \ref{bem_duality} implies the symmetry of the Green function in the form \[ G_\sigma(s,y,\tau,z) = \left(\frac{z_n}{y_n}\right)^\sigma G_\sigma(\tau,z,s,y). \] Thus $ z_n^{-\sigma} \, \ds^k \dy^\alpha G_\sigma(s_0,y_0,\tau,z) $ is a $ \sigma $-solution on $ I_{\widetilde{r}}(\widetilde{\tau}_0) \times B_{\widetilde{r}}(z_0) $ with respect to $ (\tau,z) $ and we can apply the local pointwise estimate \ref{prop_pointwiseEstimate} in the temporal end point $ \tau_0 $ and the spatial centre point $ z_0 $ to get 
	\begin{align*}
		&\left| \dta^j \dz^\beta \left( \zzeron^{- \sigma} \, \ds^k \dy^\alpha G_\sigma(s_0,y_0,\tau_0,z_0) \right) \right| \\
	&\qquad \lesssim r^{- 2 j - |\beta| } \, (r + \sqrt{\zzeron})^{-|\beta|} \supre[\tau \in I_{\widetilde{r}}(\widetilde{\tau}_0)] \supre[z \in B_{\widetilde{r}}(z_0)] \left( z_n^{-\sigma} \, \ds^k \dy^\alpha G_\sigma(s_0,y_0,\tau,z) \right)
	\end{align*} 
	The Gaussian estimate proven above then shows the statement.
\end{bew}

\begin{bem} \label{bem_defRoughInitial}
	The Gaussian estimate makes it possible to give sense to more general initial values and solve the initial value problem for data not contained in $ \Ltwosig(H) $: A $ \sigma $-solution with initial value $ u_0 $ is uniquely given through the representation by the Green function whenever this exists. The exponential decay, however, ensures the convergence of the integral in many cases, as for example for \[ u_0 \in \Lonesig(H) \cup \Linfty(H) \cup \homLip(\Hquer). \] 
	All statements up to now, and especially \ref{prop_roughInitial}, then continue to hold.
\end{bem}

We now state a series of direct consequences of the Gaussian estimate. To incorporate the temporal dimension into our intrinsic metric we proceed as in \citet{fefferman_fundamentalSolutionsSubelliptic} and define \[ \Dist((s,y),(\tau,z)) \ldef \sqrt{|s - \tau|  + \dist(y,z)^2} \text{ on } I \times \Hquer. \] Balls with respect to $ \Dist $ will be denoted $ B^D $, and we set \[ |\, \cdot \, |_{0 \times \sigma} \ldef (\Leb \times \musig)(\, \cdot \, ). \] Furthermore, abbreviate \[ V  \ldef |B^D_{\Dist((s,y),(\tau,z))}(s,y)|_{0 \times \sigma} + |B^D_{\Dist((s,y),(\tau,z))}(\tau,z)|_{0 \times \sigma} \]
	and \[ W \ldef \frac{\Dist((s,y),(\squer,\yquer)) + \Dist((\tau,z),(\taquer,\zquer))}{\Dist((s,y),(\tau,z)) + \Dist((\squer,\yquer),(\taquer,\zquer))}. \] Detailed calculations of the following proof can be found in \citet{diss}.

\begin{kor} \label{kor_CZOkernels}
	If $ l \geq 0 $, $ k \in \N_0 $ and $ \alpha \in \N_0^n $ satisfy the conditions $ l - k - |\alpha| = -1 $ and $ 2 \, l - |\alpha| \leq 0 $, then there exists a constant $ c = c(n,\sigma,l,k,\alpha) > 0 $ such that \[ z_n^{-\sigma} \, y_n^l \, |\ds^k \dy^\alpha G_\sigma(s,y,\tau,z)| \leq c \, V^{-1}  \text{ for any } \tau < s \in I \text{ and } y \neq z \in \Hquer,  \] and \[ \big| z_n^{-\sigma} \, y_n^l \, \ds^k \dy^\alpha G_\sigma(s,y,\tau,z) -  \zquer_n^{-\sigma} \, \yquer_n^l \, \ds^k \dy^\alpha G_\sigma(\squer,\yquer,\taquer,\zquer) \big| \leq c \, V^{-1} \, W \] for any $ s < \tau \in I, \, y \neq z \in \Hquer , \, \squer < \taquer \in I, \, \yquer \neq \zquer \in \Hquer $ with  $ W  \leq \frac{1}{6} $. 
\end{kor}
\begin{bew}
	Abbreviate $ \sqrt{s - \tau} \rdef r $, $ \Dist((s,y),(\tau,z)) \rdef R $ and \[ K(s,y,\tau,z) \ldef  z_n^{-\sigma} \, y_n^l \, \ds^k \dy^\alpha G_\sigma(s,y,\tau,z), \] suppressing the derivatives of $ G_\sigma $ in the notation of $ K $. Note that \[ R^{-2} \left(|B_R(y)|_\sigma + |B_R(z)|_\sigma \right)^{-1} \lesssim V^{-1}. \] However, a calculation based on the Gaussian estimate and the doubling condition of the metric $ \dist $ reveals that \[ |K(s,y,\tau,z)| \lesssim y_n^l \, r^{-2 k - |\alpha|} \, (r + \sqrt{y_n})^{-|\alpha|} \, \left (|B_R(y)|_\sigma + |B_R(z)|_\sigma \right)^{-1} \, e^{-\frac{\dist(y,z)^2}{Cr^2}} \] for any $ \tau < s \in I $ and almost any $ y \neq z \in \Hquer $.  
	On the other hand, we have that \[ y_n^l \, r^{- 2 k - |\alpha|} \, (r + \sqrt{y_n})^{-|\alpha|} \lesssim r^{2 l - 2 k - 2 |\alpha|} \] if $ 2 \, l - |\alpha| \leq 0 $. Thanks to the exponential decay we can modify this to obtain \[ r^{2 l - 2 k - 2 |\alpha|} \, e^{-\frac{\dist(y,z)^2}{C r^2}} \lesssim \sqrt{r^2 + \dist(y,z)^2}^{2(l - k - |\alpha|)} = R^{2(l - k - |\alpha|)}. \] For $ l - k - |\alpha| = - 1 $, these estimates exactly imply the desired bound. \\

	For the second part consider 
	\begin{align*}
		|K(t,x,s,y) - &K(\tquer,\xquer,\squer,\yquer)| \\
		&\leq |K(t,x,s,y) - K(\tquer,x,s,y)| +  |K(\tquer,x,s,y) - K(\tquer,\xquer,s,y)| \\
		&+  |K(\tquer,\xquer,s,y) - K(\tquer,\xquer,\squer,y)| +  |K(\tquer,\xquer,\squer,y) - K(\tquer,\xquer,\squer,\yquer)| 
	\end{align*}
	Exemplarily, we focus on the second term. If $ \Gamma $ is the arc-length parametrised geodesic that connects $ y $ and $ \yquer $ with respect to the Riemannian metric generating $ \dist $, 
	then with the chain rule and the fundamental theorem of calculus we see that
	\begin{align*}
		| K(\squer,y,\tau,z) -  K(\squer,\yquer,\tau,z)| &= \bigl| \inte{(a,b)}{\grad_{\Gamma(\varsigma)} K(\tquer,\Gamma(\varsigma),s,y) \cdot \partial_\varsigma \Gamma(\varsigma)}{\Leb(\varsigma)} \bigr| \\
		&\leq \dist(y,\yquer) \supre[\zeta \in B_{\dist(y,\yquer)}(y)] \zeta_n^\frac{1}{2} \, |\grad_\zeta K(\squer,\zeta,\tau,z)| 
	\end{align*} 
	with $ \zeta = \Gamma(\varsigma) $. The Gaussian estimate and some calculations based on the condition $ W \leq \frac{1}{6} $ show
	\begin{align*}
		| K(\squer,y,\tau,z) -  K(\squer,\yquer,\tau,z)| 
		\lesssim \frac{\dist(y,\yquer)}{R + \Dist((\squer,\yquer),(\taquer,\zquer))} \, R^{2 l - 2 k - 2 |\alpha| + 2} \, V^{-1}.
	\end{align*} 
	We proceed similarly for the other terms, where we also need the generalised version of the Gaussian estimate from Proposition \ref{prop_gaussian}.
\end{bew} 

We now prove an integral bound for the Green function where the conditions on the exponents differ from the ones we considered before. To simplify the appearance of the estimates we specialise to unit time intervals. 

\begin{kor} \label{kor_GreenLq}
	Let $ I \ldef (0,1) $ as well as $ l \geq 0 $ and $ \alpha \in \N_0^n $ such that $ 2 \, l -|\alpha| \leq 0 $ and $ l - |\alpha| > - 1 $.
	\begin{enumerate}
		\item
		If $ 1 \leq q < \frac{n+1}{n + |\alpha| - l} $ and $ \sigma \, q > -1 $, then there exists a constant $ c = c(n,\sigma,k,\alpha,q) $ such that \[ \norm[\Lq((0,s) \times H)]{y_n^l \, \dy^\alpha G_\sigma(s,y,\, \cdot\, ,\, \cdot\, )} \leq c \, (1 + \sqrt{y_n})^{2l - |\alpha|} \, |B_1(y)|^{\frac{1}{q}-1}  \] for all $ s \in \Iquer $ and almost all $ y \in \Hquer $.
		\item
		If $ 1 \leq q < \frac{n+1}{n + |\alpha| - l} $, then there exists a constant $ c = c(n,\sigma,l,k,\alpha,q) $ such that  \[ \norm[\Lq((\tau,1) \times H)]{\Pnl \, \dy^\alpha G_\sigma(\, \cdot\, ,\, \cdot\, ,\tau,z)} \leq c \, (1 + \sqrt{z_n})^{2l - |\alpha|} \, |B_1(z)|^{\frac{1}{q}-1} \, \left( z_n^\sigma \right)^{\chi_{\menge{\sigma < 0}} \chi_{\menge{z_n < 1}}}  \] for all $ \tau \in \Iquer $ and almost all $ z \in \Hquer $.
	\end{enumerate} 
\end{kor}
\begin{bew}
	First consider the spatial integrals only. Fix $ (s,y) \in \Irechts \times \Hquer $ as well as $ 1 \leq q  < \infty $. For $ \tau \in (0,s) $ define $ r \ldef \sqrt{s - \tau} $. The Gaussian estimate implies \[ \inte{H}{|\ds^k \dy^\alpha G_\sigma(s,y,\tau,z)|^q}{\Lebn(z)} \lesssim r^{- 2 k q - |\alpha| q} \, (r + \sqrt{y_n})^{- |\alpha| q} \, |B_r(y)|_\sigma^{-q} \, \summe{j \in \N} e^{-q \frac{(j-1)^2 r^2}{C r^2}} \inte{B_{j r}(y)}{z_n^{\sigma q}}{\Leb(z)}, \]
	where we covered $ H $ with annular rings $ B_{j r}(y) \ohne B_{(j-1) r}(y) $ for $ j \in \N $. At this point we have to assume $ \sigma q > -1 $ to be sure that the integral exists. By the doubling condition we get a convergent series with value $ |B_r(y)|_{\sigma q} $ times a constant depending on $ n, \, \sigma $ and $ q $, and the formula for the measure of the balls shows that \[ |B_r(y)|_\sigma^{-q} |B_r(y)|_{\sigma q} \eqsim |B_r(y)|^{1 - q}. \] This results in 
	\begin{align*}
		\norm[\Lq(H)]{y_n^l \, \ds^k \dy^\alpha G_\sigma(s,y,\tau,\, \cdot\, )} \lesssim r^{- 2 k - |\alpha|} \, (r + \sqrt{y_n})^{2 l - |\alpha|} \, |B_r(y)|^{\frac{1}{q} - 1} 
	\end{align*}
	for all $ \tau < s \in \Iquer $ and almost all $ y \in \Hquer $, and for any $ q $ with $ \sigma \, q > -1 $. \\
	Similarly we get  
		\begin{align*}
			\norm[\Lq(H)]{\Pnl \, \ds^k \dy^\alpha G_\sigma(s,\, \cdot\, ,\tau,z)} \lesssim r^{- 2 k - |\alpha|} \, (r + \sqrt{z_n})^{2 l - |\alpha|} \, |B_r(z)|^{\frac{1}{q} - 1} \left( \frac{\sqrt{z_n}}{r + \sqrt{z_n}} \right)^{2 \sigma} 
		\end{align*}
		 for all $ \tau < s \in \Iquer $ and almost all $ z \in \Hquer $, without any restriction on $ q \in [1,\infty)$, but with an additional dependency of the constant on $ l $. \\

	Integrating the first formula in time with respect to $ \tau $ over $ (0,s) $ then implies
	\begin{align*}
		\norm[\Lq((0,s) \times H)]{y_n^l \, \ds^k \dy^\alpha G_\sigma(s,y,\, \cdot\, ,\, \cdot\, )}^q \lesssim (1 + \sqrt{y_n})^{(2l - |\alpha|)q} \, |B_1(y)|^{1-q} \inte{(0,s)}{\tau^{(l- k - |\alpha| -n)q + n}}{\Leb(\tau)}
	\end{align*}
	for $ 2 \, l - |\alpha| \leq 0 $. Here we used the unit interval in the computation. For any $ q < \frac{1+n}{-l + k + |\alpha|} $, this integral converges to \[ c(n, \alpha, k, l, q) \,  s^{(l - k - |\alpha| - n)q + n + 1} \] which in turns is bounded for $ s \in (0,1] $. It is obvious that the convergence condition can only be satisfied if $ -l +  k + |\alpha| < 1 $. In conjunction with the condition $ 2 \, l - |\alpha| \leq 0 $ that arouse from the calculations above it becomes clear that only $ k = 0  $ is admissible and we need to have $ - l + |\alpha| < 1 $. \\
	For the second part of the statement, the additional factor 
	is bounded by an absolute constant for $ \sigma \geq 0 $, as well as for $ -1 < \sigma < 0 $ if in addition $ \sqrt{z_n} \geq \sqrt{|I|} $. It thus only shows up for negative $ \sigma $ and for $ \sqrt{z_n} < \sqrt{|I|} $, where in the present case of the unit interval it can be bounded by a multiple of $ z_n^\sigma $.
\end{bew} 

\begin{bem} \label{bem_GreenGradientNot}
	It is worth pointing out that no temporal derivative can be treated in the fashion of Corollary \ref{kor_GreenLq}. We cannot even achieve a bound for the gradient of $ G_\sigma $ in $ \Lq $ without extra weights to compensate for the spatial derivative either. In fact, the condition on the exponents $ l $ and $ \alpha $ are satisfied if and only if either $ l = 0 $ and $ \alpha = 0 $ or $ l \in (0,\frac{1}{2}] $ and $ |\alpha| = 1 $. 
\end{bem}

The final corollary is a pointwise estimate of Green's function on $ (0,1) \times \Hquer $ on a rather complicated range of values that follows immediately by the Gaussian estimate on the unit interval.

\begin{kor} \label{kor_GreenOffDiag} 
	Let $ I \ldef (0,1) $ and $ 0 < \delta \leq \frac{1}{2} $. 
	If $ l \geq 0 $, $ k \in \N_0 $ and $ \alpha \in \N_0^n $ with $ 2 \, l -|\alpha| \leq 0 $, then there exists a constant $ c = c(n,\sigma,l,k,\alpha,\delta) > 0 $ such that \[ y_n^l \, |\ds^k \dy^\alpha G_\sigma(s,y,\tau,z)| \leq c \, (1 + \sqrt{y_n})^{2l - |\alpha|} \, |B_1(z)|^{-1} \, (z_n^\sigma)^{\chi_{\menge{\sigma < 0}} \chi_{\menge{z_n < 1}}} \, e^{- \frac{\dist(y,z)}{C} } \] for all $ s \in [2 \, \delta,1] $, almost all $ y \in \Hquer $, and almost all $ (\tau,z) \in \left( [0,s) \times \Hquer \right) \ohne \left( (\delta,s) \times B_1(y) \right) $, and the same statement holds with any possible combination of the points $ y $ and $ z $ in the factors $ (1 + \sqrt{y_n})^{2l - |\alpha|} \, |B_1(z)|^{-1} $. 
\end{kor}
 
\subsection{Estimates Against the Inhomogeneity}

We now abandon the situation where $ f = 0 $ and gain estimates against the inhomogeneity. 

\begin{bem} \label{bem_Duhamel}
	By Duhamel's principle, $ \sigma $-solutions $ u $ to $ f $ on $ I \times \Hquer $ with zero initial value can be expressed in terms of $ G_\sigma $ by \[ u(s,y) = \inte{(s_1,s) \times H}{G_\sigma(s,y,\tau,z) f(\tau,z)}{\Lebnp(\tau,z)} \] whenever the integral makes sense. We can thus consistently extend Definition \ref{defin_energySolution} to more general inhomogeneities 
	and will henceforth always assume the inhomogeneity to be given such that the Duhamel integral converges without specifying the conditions on $ f $.
\end{bem} 

Duhamel's principle also enables us to view the operator that maps inhomogeneities to solutions as integral kernel operators. 
Since $ (I \times \Hquer, \Leb \times \musig,\Dist) $ is a space of homogenous type, the theory of Calder\'{o}n-Zygmund can be applied in this non-Euclidean setting to gain $ \Lp $-estimates. The basic reference for this material is \citet{harmonic_analysis}. It was noted by \citet{coifweiss} that the arguments do not depend on the Euclidean structure of the underlying space. The definitions and statements can be found in \citet{koch_quaderni} and \citet{koch_CZO} as well as \citet{diss}. 

\begin{prop} \label{prop_CZO}
	If $ \max\menge{1,(1+\sigma)^{-1}} < p < \infty $, if $ l \geq 0 $, $ k \in \N_0 $ and $ \alpha \in \N_0^n $ satisfy the conditions $ l - k - |\alpha| = - 2 $ and $ 2 \, l - |\alpha| \leq -1 $, and if $ u $ is a $ \sigma $-solution to $ f $ on $ I \times \Hquer $ with zero initial value, then there exists a constant $ c = c(n,\sigma,l,k,\alpha,p) $ such that \[ \norm[\Lp(I \times H)]{\Pnl \, \ds^k \dy^\alpha u} \leq c \, \norm[\Lp(I \times H)]{\grady f}. \]
\end{prop}
\begin{bew}
	For fixed $ l \geq 0 $, $ k \in \N_0 $ and $ \alpha \in \N_0^n $ consider the operator \[ \Ltwo(I \times H, \Leb \times \musig) \ni f \mapsto \Pnl \, \ds^k \dy^\alpha u \in \Ltwo(I \times H, \Leb \times \musig) \] with kernels \[ K(s,y,\tau,z) = z_n^{-\sigma} \, y_n^l \, \ds^k \dy^\alpha G_\sigma(s,y,\tau,z) \] that are continuous where $ (s,y) \neq (\tau,z) $. By Corollary \ref{kor_CZOkernels}, these kernels satisfy the Calder\'{o}n-Zygmund cancellation conditions on $ ( I \times \Hquer, \Dist, \Leb \times \musig) $ if $ l - k - |\alpha| = - 1 $ and $ 2 \, l - |\alpha| \leq 0 $. As noted in Remark \ref{bem_LtwoCZO}, these are exactly the exponents for which the operators are bounded on $ \Ltwo(I \times H, \Leb \times \musig) $. It follows that those operators are of Calder\'{o}n-Zygmund type and hence extend to bounded operators on $ \Lp(I \times H, \Leb \times \musig) $ for any $ 1 < p < \infty $ with the estimate \[ \inte{I}{\norm[\Lpsig(H)]{\Pnl \, \ds^k \dy^\alpha u}^p}{\Leb} \lesssim \inte{I}{\norm[\Lpsig(H)]{f}^p}{\Leb}. \] 
	
	The theory of Muckenhoupt weights allows us to formulate this result in unweighted spaces. One checks that $ y_n^{\vartheta - \sigma} $ is a $ p $-Muckenhoupt weight with respect to $ \musig $ if and only if \[ -1 < \vartheta < p \, (\sigma + 1) -1. \] Thus any Calder\'{o}n-Zygmund operator with respect to $ \musig $ also extends to all spaces $ \Lp(x_n^{\vartheta-\sigma}d\musig) = \Lp(\mu_\vartheta) $, and the same is true for the corresponding time-space setting. Note that $ \vartheta = 0 $ is always admissible if $ \sigma \geq 0 $, whereas for negative $ \sigma $ we need the extra condition $ p > (1+\sigma)^{-1} $. On the other hand, this does not constitute a restriction for $ p $ in case $ \sigma \geq 0 $. For this range of $ p $ it follows that for $ l - k - |\alpha| = - 1 $ and $ 2 \, l - |\alpha| \leq 0 $ we have \[ \norm[\Lp(I \times H)]{\Pnl \, \ds^k \dy^\alpha u} \lesssim \norm[\Lp(I \times H)]{f}. \] 
	Now that we have gotten rid of the weights, we can also take extra derivatives into consideration. Given that $ \gradx' u $ is a $ \sigma $-solution to $ \gradx' f $ in view of Remark \ref{bem_iteratedSigmaSol}, we have also shown that \[ \norm[\Lp(I \times H)]{\Pnl \, \ds^k \dy^\alpha \grady' u} \lesssim \norm[\Lp(I \times H)]{\grady' f} \] for any exponents $ l $, $ k $ and $ \alpha $ as above, and especially also for $ l = 0 $, $ k = 0 $ and $ |\alpha| = 1 $. \\
	For the vertical direction, however, we can use that $ \dyn u $ is a $ (1+\sigma) $-solution to \[ \dyn f + (1+\sigma) \, \lapy' u \] by Remark \ref{bem_iteratedSigmaSol}. Thus \[ \norm[\Lp(I \times H)]{\Pnl \, \ds^k \dy^\alpha \dyn u} \lesssim \norm[\Lp(I \times H)]{\dyn f} + \norm[\Lp(I \times H)]{\lapy' u} \] which has the required bound because of the previous considerations. 
\end{bew} 

We normalise our interval once more and consider $ I = (0,1) $. For inhomogeneities $ f $ that are supported on the diagonal, more precisely on increased cylinders \[ \widetilde{Q}_r(y_0) \ldef (\frac{1}{4} r^2,r^2) \times B_{2r}(y_0), \] we get a local version of the $ \Lp $-estimates in terms of the space $ Y^{on}_\theta $ introduced in Section \ref{subsection_perturbSetting}. 

\begin{prop} \label{prop_LpOnDiag}
	Let $ I \ldef (0,1) $ and $ y_0 \in \Hquer $ as well as $ \supp f \subset \overline{\widetilde{Q}_1(y_0)} $ and $ \theta \in \R $. If $ \max\menge{1,(1+\sigma)^{-1}} < p < \infty $, if $ l \geq 0 $, $ k \in \N_0 $ and $ \alpha \in \N_0^n $ satisfy the conditions $ l - k - |\alpha| = - 2 $ and $ 2 \, l - |\alpha| \leq -1 $, and if $ u $ is a $ \sigma $-solution to $ f $ on $ I \times \Hquer $ with zero initial value, then there exists a constant $ c = c(n,\sigma,l,k,\alpha,p,\theta) > 0 $ such that
	\begin{align*} 
		(1 + \sqrt{\yzeron})^{-2 l + |\alpha| - 1} \, |Q_1(y_0)|^{-\frac{1}{p}} \, \norm[\Lp(Q_1(y_0))]{\Pnl \, \ds^k \dy^\alpha u} \leq c \,  \Yon{f}.
	\end{align*}
\end{prop}
\begin{bew}
	The Calder\'{o}n-Zygmund estimate from Proposition \ref{prop_CZO} immediately yields \[ \norm[\Lp(Q_1(y_0))]{y_n^l \, \ds^k \dy^\alpha u} \lesssim \norm[\Lp(\widetilde{Q}_1(y_0) )]{\grady f} \]
	for the parameters stated. \\ 
	Consider now radii $ R_1 \ldef \frac{1}{\sqrt{2}}$, $ R_2 \ldef 1 $ and $ R_3 \ldef \frac{\sqrt{3}}{2}$. Independent of the location of $ y_0 $ in $ \Hquer $ we can find $ N_1 $ points $ \menge{z_{j_1}} $, $ N_2 $ points $ \menge{z_{j_2}} $ and $ N_3 $ points $ \menge{z_{j_3}} $, all of them contained in $ \widetilde{Q}_1(y_0) $ and with $ N_i $ only depending on $ n $ for $ i = 1,\, 2, \, 3 $, such that the collection $ \menge{Q_{R_i}(z_{j_i}) \mid j_i = 1, \ldots, N_i, \, i = 1, \, 2 \, ,3} $ covers $ \widetilde{Q}_1(y_0) $. It follows that \[ \norm[\Lp(\widetilde{Q}_1(x_0) )]{\grady f} \lesssim |Q_1(y_0)|^\frac{1}{p} \, \Yon{f}, \]
	where we used the doubling property to shrink balls into $ Q_1(y_0) $, and the constant depends on $ \theta $. This is the localisation of the Calder\'{o}n-Zygmund estimate. \\
	But for the set of exponents we consider, the only case for which $ - 2 \, l + |\alpha| - 1 \geq 0 $ is non-zero is the one with $ l = 0 $, $ k = 0 $ and $ |\alpha| = 2 $. Furthermore, if we fix $ y_0 $ with $ \sqrt{\yzeron} < 2 $, we have $ (1 + \sqrt{\yzeron})^{- 2 l + |\alpha| - 1} \eqsim 1 $ and the inequlity we showed is equivalent to the statement of this proposition. \\

	Let now $ l = 0 $, $ k = 0 $, $ |\alpha| = 2 $ and $ \sqrt{\yzeron} \geq 2 $. We study the kernels \[ K(s,y,\tau,z) = y_n^\frac{1}{2} \, \grady G_\sigma(s,y,\tau,z) \] for the unweighted $ \Ltwo $-operators that send $ f $ to $ \Pn^\frac{1}{2} \, \grady u $. For $ \sigma \geq 0 $, an application of Corollary \ref{kor_GreenLq} with $ q = 1 $ ensures that the $ \Lone $-norm of this kernels taken with respect to $ (s,y) $ as well as the ones taken with respect to $ (\tau,z) $ are bounded. Thus Schur's Lemma can be used to conclude that \[ \norm[\Lp(I \times H)]{\Pn^\frac{1}{2} \, \grady u} \lesssim \norm[\Lp(I \times H)]{f} \] for any $ 1 \leq p \leq \infty $. As before, we now use that $ \dyj u $ is a $ \sigma $-solution to $ \dyj f $ for $ j = 1, \ldots, n -1 $, while for the vertical direction $ \dyn u $ is a $ (1+\sigma) $-solution to $ \dyn f + (1+\sigma) \lapy' u $. Therefore, \[ \norm[\Lp(I \times H)]{\Pn^\frac{1}{2} \, D_y^2 u} \lesssim \norm[\Lp(I \times H)]{\grady f} + \norm[\Lp(I \times H)]{\lapy' u}, \] and the last summand is bounded by the second to last thanks to the unweighted Calder\'{o}n-Zygmund-estimate from Proposition \ref{prop_CZO}. 

	The resulting estimate can be localised as above, and we get 
	\begin{align*}
		\norm[\Lp(Q_1(y_0))]{\Pn^\frac{1}{2} \, D_y^2 u} &\lesssim |Q_1(y_0)|^\frac{1}{p} \, \Yon{f}.
	\end{align*}
	However, for any $ y \in B_1(y_0) $ we have that $ y_n \eqsim \yzeron $ if $ \sqrt{\yzeron} \geq 2 $ as in our present case. 
	Together with the first result for $ l = 0 $, $ k = 0 $ and $ |\alpha| = 2 $, this amounts to the improved inequality
	\begin{align*}
		(1 + \sqrt{\yzeron}) \, \norm[\Lp(Q_1(y_0))]{D^2_y u} \lesssim |Q_1(y_0)|^\frac{1}{p} \, \Yon{f} 
	\end{align*}
	in case of $ \sqrt{\yzeron} \geq 2 $, proving the statement for $ \sigma \geq 0 $. \\

	For the last remainig case, we now consider $ -1 < \sigma < 0 $ and $ \sqrt{\yzeron} \geq 2 $. The latter also means that we have $ \supp f \subset \menge{z \mid z_n \geq 1} $, and thus, also in the case where $ \sigma $ is not in the good range, Corollary \ref{kor_GreenLq} can be applied to the kernels $ K $ truncated to the support of $ f $ without the additional factor interfering. This makes it possible to invoke Schur's Lemma once more and get the same bound as before. The iteration in the tangential direction can be done as before. For a similar treatment of the vertical derivative, however, we lack the guarantee that the corresponding inhomogeneity $ \dyn f + (1 + \sigma) \, \lapy' u $ has support in the right region. But $ \dyn u $ is a $ (1 + \sigma) $-solution to this right hand side, and since $ (1 + \sigma) \geq 0 $ we can apply the result we obtained in case of positive $ \sigma $ onto this solution. Proceeding in the same fashion as before finishes the proof.
\end{bew}


We now turn to inhomogeneities that are supported away from the diagonal. 

\begin{prop} \label{prop_offDiagPointwise}
	Let $ I \ldef (0,1) $, $ \varepsilon_1 > 0 $ and $ \varepsilon_2 \geq 0 $. Fix $ y \in \Hquer $ as well as $ 0 < \delta \leq \frac{1}{2} $ and let $ \supp f \subset \left( \Iquer \times \Hquer \right) \ohne \left( (\delta,1) \times B_1(y) \right) $. If $ \max\menge{1,(1 + \sigma)^{-1}} < p < \infty $, if $ l \geq 0 $, $ k \in \N_0 $ and $ \alpha \in \N_0^n $ with $ 2 \, l - |\alpha| \leq 0 $, and if $ u $ is a $ \sigma $-solution to $ f $ on $ I \times \Hquer $ with zero initial value, then there exists a constant $ c = c(n,\sigma,l,k,\alpha,p,\delta,\varepsilon_2) > 0 $ such that 
	\begin{align*}
		(1 + \sqrt{y_n})^{-2l + |\alpha| - \varepsilon_2} \, y_n^l \, |\ds^k \dy^\alpha u(s,y)| \leq c \, \Yoff{f}
	\end{align*}
	for all $ s \in [2 \, \delta, 1] $.
\end{prop} 
\begin{bew} 
	We first consider the case $ \sigma \geq 0 $. Fix $ s \in [2 \, \delta,1] $ and denote \[ M_s \ldef \left( (0,s) \times \Hquer \right) \ohne \left( (\delta,s) \times B_1(y) \right). \] By Duhamel's principle \ref{bem_Duhamel} we get \[ y_n^l \, |\ds^k \dy^\alpha u(s,y)| \leq \inte{M_s}{y_n^l \, |\ds \dy^\alpha G_\sigma(s,y,\tau,z)| \, |f(\tau,z)|}{\Lebnp(\tau,z)}. \] This is only possible for $ \delta \leq s \leq 1 $, and the application of Corollary \ref{kor_GreenOffDiag} under the integral is justified for $ 2 \, \delta \leq s \leq 1 $ if $ 2 \, l - |\alpha| \leq 0 $. As a result we gain the upper bound \[ (1 + \sqrt{y_n})^{2 l - |\alpha| + \varepsilon_2}  \inte{(0,1) \times H}{(1 + \sqrt{z_n})^{- \varepsilon_2} \, |B_1(z)|^{-1} \, e^{- \frac{\dist(y,z)}{C}} \, |f(\tau,z)|}{\Lebnp(\tau,z)} \] for an arbitrary $ \varepsilon_2 \geq 0 $ and with a constant that depends on $ n, \, \sigma, \, l, \, k,\, \alpha $ and $\delta $. 
	For the computation of the integral we first cover $ \Hquer $ with a countable number of balls $ B_1(z_0) $. This leads to 
	\begin{align*}
		\inte{(0,1) \times H}{(1 + \sqrt{z_n})^{- \varepsilon_2} \, &|B_1(z)|^{-1} \, e^{- \frac{\dist(y,z)}{C}} \, |f(\tau,z)|}{\Lebnp(\tau,z)} \\
		&\leq \supre[z_0] \inte{(0,1) \times B_1(z_0)}{(1 + \sqrt{z_n})^{- \varepsilon_2} \, |B_1(z)|^{-1} \, |f(\tau,z)|}{\Lebnp(\tau,z)} \, \summe{z_0} e^{- \frac{\dist(y,z_0)}{C}}
	\end{align*}
	and the series in the back converges uniformely in $ y $. \\ 
	We now cover the time interval $ (0,1) $ by $ (\frac{1}{2} \, R_m^2,R_m^2) $, $ m \in N_0 $, where $ R_m \ldef 2^{-\frac{m}{2}} $. Furthermore, for any $ m \in \N_0 $ and any $ z_0 $ we cover $ B_1(z_0) $ with $ N(m) $ balls $ B_{R_m}(z_i) $, $ z_i \in B_1(y_0) $. This is possible by Vitali's covering lemma \citep{koch_quaderni}), which also ensures that \[ \summe[N(m)]{i=1} |B_{R_m}(z_i)|_\sigma \lesssim |B_1(z_0)|_\sigma \] independent of $ m $. Note that $ (2^{-m-1},2^{-m}) \times B_{R_m}(z_i) = Q_{R_m}(z_i) $ and we therefore now look at the expression \[ \supre[z_0] \summe{m \in \N_0} \summe[N(m)]{i = 1} \inte{Q_{R_m}(z_i)}{(1 + \sqrt{z_n})^{- \varepsilon_2} \, |B_1(z)|^{-1} \, |f(\tau,z)|}{\Lebnp(\tau,z)}. \]
	We consider the cases $ \sqrt{\zzeron} < 2 $ and $ \sqrt{\zzeron} \geq 2 $ separately, starting with the latter one away from the boundary. 
	Then $ |B_1(z)| \eqsim_n |B_1(z_0)| $ for any $ z \in B_1(z_0) $ 
	and thus also \[ (1 + \sqrt{z_n})^{- \varepsilon_2} \eqsim_n (R_m + \sqrt{z_{i,n}})^{- \varepsilon_2} \text{ for any } z \in B_{R_m}(z_i) \text{ with } z_i \in B_1(z_0). \] An application of Hölder's inequality shows that 
	\begin{align*}
		&\supre[\zzeron \geq 4] |B_1(z_0)|^{-1} \summe{m \in \N_0} \summe[N(m)]{i = 1} (R_m + \sqrt{z_{i,n}})^{- \varepsilon_2} \inte{Q_{R_m}(z_i)}{|f(\tau,z)|}{\Lebnp(\tau,z)} \\
		&\leq \supre[\zzeron \geq 4] |B_1(z_0)|^{-1} \summe{m \in \N_0} R_m^{\varepsilon_1} \summe[N(m)]{i = 1}  (R_m + \sqrt{z_{i,n}})^{- \varepsilon_2} \, |B_{R_m}(z_i)| \, R_m^{2 - \varepsilon_1} \, |Q_{R_m}(z_i)|^{-\frac{1}{p}} \, \norm[\Lp(Q_{R_m}(z_i))]{f}.
	\end{align*}
	Taking the supremum over $ 0 < R \leq 1 $ and $ z \in \Hquer $ then leads to the upper bound
	\begin{align*}
	\Yoff{f} \, \summe{m \in \N_0} R_m^{\varepsilon_1}
	\end{align*}
	by virtue of the Vitali-property discussed above. For any $ \varepsilon_1 > 0 $ also this series converges and can therefore be subsumed into the constant. \\
	In the case that $ z_0 $ is close to the boundary we vary the arguments slightly to obtain 
	\begin{align*}
		\supre[\zzeron < 4] \summe{m \in \N_0} \summe[N(m)]{i = 1} \inte{Q_{R_m}(z_i)}{|f(\tau,z)|}{\Lebnp(\tau,z)} \leq \supre[\zzeron < 4] |B_1(z_0)| \, \Yoff{f} \summe{m \in \N_0} R_m^{\varepsilon_1}.
	\end{align*}
	But for $ \sqrt{\zzeron} < 2 $ it is clear that $ |B_1(z_0)| \lesssim 1 $ and thus the same upper bound as above follows for any $ \varepsilon_1 > 0 $. \\
	It remains to consider the case $ - 1 < \sigma < 0 $, where the use of Corollary \ref{kor_GreenOffDiag} generates an additional factor $ z_n^\sigma $ under the integral if also $ z_n < 1 $. On $ B_1(z_0) $, this possibility only realises itself if $ \sqrt{\zzeron} < 2 $. The application of Hölder's inequality then results in \[ \supre[\zzeron < 4] \summe{m \in \N_0} R_m^{\varepsilon_1} \, \summe[N(m)]{i = 1} (R_m + \sqrt{z_{i,n}})^{- \varepsilon_2} \, R_m^{2(1 - \frac{1}{p}) - \varepsilon_1} \, |B_{R_m}(z_i)|_{\sigma\frac{p}{p-1}}^{1-\frac{1}{p}} \, \norm[\Lp(Q_{R_m}(z_i))]{f} \] if $ \sigma \, (1 - \frac{1}{p}) > - 1 $. The last condition is equivalent to $ p > (1 + \sigma)^{-1} $ for $ -1 < \sigma < 0 $, and does not pose a restriction on $ p $ if $ \sigma \geq 0 $. Now the formula for the measure of intrinsic balls asserts that \[ |B_{R_m}(z_i)|_{\sigma \frac{p}{p-1}}^{1-\frac{1}{p}} \eqsim |B_{R_m}(z_i)|^{1 - \frac{1}{p}} \, (R_m + \sqrt{z_i,n})^{2 \sigma} \eqsim |B_{R_m}(z_i)|_\sigma \, |B_{R_m}(z_i)|^{-\frac{1}{p}} \] and a reiteration of the steps above, this time using Vitali's lemma for the measure $ \musig $, implies the upper bound \[ \supre[\zzeron < 4] |B_1(z_0)|_\sigma \, \Yoff{f} \, \summe{m \in \N_0} R_m^{\varepsilon_1}. \] The fact that also $ |B_1(z_0)|_\sigma \lesssim 1 $ if $ \sqrt{\zzeron} < 2 $ concludes the proof.
\end{bew}


We have now finally reached the point where we can prove the second part of Theorem \ref{theo_linear}. To this end, we combine the off-diagonal pointwise estimate from Proposition \ref{prop_offDiagPointwise} with a pointwise estimate for an arbitrary inhomogeneity on $ (0,1) $, and integrate the off-diagonal pointwise estimate to get the off-diagonal complement of the on-diagonal $ \Lp $ estimate from Proposition \ref{prop_LpOnDiag}. Subsequently, both results are rescaled onto $ (0,S) $ by the invariant scaling $ A_\lambda $. In the following we still allow $ S = \infty $ and use the convention that $ r^2 \leq S $ means $ r^2 < \infty $ in this case. 

\begin{prop} \label{prop_XggY}
	If $ \max\menge{2 \, (n + 1),(1 + \sigma)^{-1}} <  p < \infty $, $ f \in Y(p) $, and $ u $ is a $ \sigma $-solution to $ f $ on $ (0,S) \times \Hquer $ with zero initial value, then there exists a constant $ c = c(n,\sigma,p) > 0 $ such that \[ \norm[X(p)]{u} \leq c \, \norm[Y(p)]{f}. \] 
	\vspace{-0.5cm}
\end{prop}
\begin{bew}
	Throughout this proof, the dependency of the constants on the parameters will be suppressed in the notation. \\ 
	Assume first $ S = 1 $ and let $ I \ldef (0,1) $. Fix an arbitrary $ y \in \Hquer $ and set \[ f_1 \ldef \chi_{Q_1(y)} \, f  \text{ with } \supp f_1 \subset \overline{Q_1(y)}, \] and \[ f_2 \ldef ( 1 - \chi_{Q_1(y)}) \, f \text{ with } \supp f_2 \subset \left( \Iquer \times \Hquer \right) \ohne Q_1(y). \] Denoting the $ \sigma $-solutions to $ f_1 $ and $ f_2 $ on $ I \times \Hquer $ with zero initial value by $ u_1 $ and $ u_2 $, respectively, we also get $ u = u_1 + u_2 $. \\ 

	For the on-diagonal part $ u_1 $ we obtain
	\begin{align*} 
		y_n^l \, |\dy^\alpha u_1(1,y)| &\leq \inte{I \times H}{|\dy^\alpha G_\sigma(1,y,\tau,z)| \, |f_1(\tau,z)|}{\Lebnp(\tau,z)} \\
		& \leq \norm[L^{\frac{p}{p-1}}(I \times H)]{\dy^\alpha G_\sigma(1,y, \, \cdot \, , \, \cdot \, )} \, \norm[\Lp(I \times H)]{f_1} 
	\end{align*} 
	with Hölder's inequality for any $ 1 \leq p \leq \infty $. But for $ 2 \, l - |\alpha| \leq 0 $ and $ l - |\alpha| > - 1 $, the first factor can be treated with Corollary \ref{kor_GreenLq} if both $ 1 \leq \frac{p}{p - 1} < \frac{n + 1}{n + |\alpha| - l} $ and $ \sigma \frac{p}{p - 1} > - 1 $. Iterating the notion of $ \sigma $-solution with respect to derivatives once more while using the Calder\'{o}n-Zygmund estimate from Proposition \ref{prop_CZO} to bound the norm of the Laplacian, we get \[ y_n^l \, |\dy^\alpha u_1(1,y)| \lesssim (1 + \sqrt{y_n})^{2 l - |\alpha| + 1} \, |B_1(y)|^{-\frac{1}{p}} \norm[\Lp(I \times H)]{\grady f} \] for $ 2 \, l - |\alpha| \leq - 1 $ and $ l - |\alpha| > - 2 $, and especially for $ l = 0 $, $ |\alpha| = 1 $ and $ l = \frac{1}{2} $, $ |\alpha| = 2 $. The conditions on $ p $ then become equivalent to $ \max\menge{2 \, (n + 1),(1 + \sigma)^{-1}} <  p < \infty $. All together we arrive at \[ |\grady u_1(1,y)| + y_n^\frac{1}{2} \, |D^2_y u_1(1,y)| \lesssim |B_1(y)|^{-\frac{1}{p}} \, \norm[\Lp(I \times H)]{\grady f_1}. \] Note that we did not use the on-diagonal support of $ f_1 $ in this calculation. \\
	In conjunction with the pointwise off-diagonal estimate from Proposition \ref{prop_offDiagPointwise}, the latter for $ \delta = \frac{1}{2} $, $ \varepsilon_2 = 1 $, and $ (l,k,|\alpha|) = (0,0,1) $ as well as $ (l,k,|\alpha|) = (\frac{1}{2},0,2) $, we then find 
	\begin{align*}
		|\grady u(1,y)| + y_n^\frac{1}{2} \, |D^2_y u(1,y)| 
		\lesssim |B_1(y)|^{-\frac{1}{p}} \, \norm[L^{p}(Q_1(y))]{\grady f_1} + \Yoffone{f_2} \leq \Yon{f} + \Yoffone{f} 
	\end{align*}
	for any $ \theta \in \R $, $ \varepsilon_1 > 0 $ and $ \max\menge{2 \, (n + 1),(1+\sigma)^{-1}} <  p < \infty $.

	For the $ \Lp $-estimate, we split $ f $ slightly differently. So fix an arbitrary $ y_0 \in \Hquer $ and consider \[ f = f \, \chi_{\widetilde{Q}_1(y_0)} + f \, (1 - \chi_{\widetilde{Q}_1(y_0)}) \rdef f_1 + f_2. \] Note that then \[ \supp f_1 \subset \overline{\widetilde{Q}_1(y_0)} \text{ and } \supp f_2 \subset \left( \Iquer \times \Hquer \right) \ohne \left( \left(\frac{1}{4},1 \right) \times B_2(y_0) \right). \] Again we study the $ \sigma $-solutions $ u_1 $ and $ u_2 $ to $ f_1 $ and $ f_2 $, respectively. \\
	For the on-diagonal part we use Proposition \ref{prop_LpOnDiag} and immediately get for almost any $ y_0 \in \Hquer $ that  
	\begin{align*}
		(1 + \sqrt{\yzeron})^{-2 l + |\alpha| - 1} \, |Q_1(y_0)|^{-\frac{1}{p}} \, \norm[L^{p}(Q_1(y_0))]{\Pnl \, \ds^k \dy^\alpha u_1} \lesssim \Yon{f_1}
	\end{align*}
	for any $ \theta \in \R $, if $ \max\menge{1,(1+\sigma)^{-1}} < p < \infty $ and both $ l - k - |\alpha| + 1 = - 1 $ and $ 2 \, l - |\alpha| + 1 \leq 0 $. \\
	Off-diagonal, on the other hand, we first note the triangle inequality implies $ B_1(y) \subset B_2(y_0) $ for any $ y \in \overline{B}_1(y_0) $. 
	For any such $ y $ and any $ s \in [\frac{1}{2},1] $, that is for $ (s,y) \in \overline{Q_1(y_0)} $ it is therefore clear that \[ \supp f_2 \subset \left( \Iquer \times \Hquer \right) \ohne \left( \left(\frac{1}{4},1 \right) \times B_1(y) \right) \] and thus the requirements needed for Proposition \ref{prop_offDiagPointwise} are met with $ \delta = \frac{1}{4} $. 
	An integration yields 
	\begin{align*}
		(1 + \sqrt{\yzeron})^{-2 l + |\alpha| - \varepsilon_2} \, |Q_1(y_0)|^{-\frac{1}{p}} \, \norm[L^{p}(Q_1(y_0))]{\Pnl \, \ds^k \dy^\alpha u_2} \lesssim \Yoff{f_2} 
	\end{align*}
	for any $ \varepsilon_2 \geq 0 $, $ \varepsilon_1 > 0 $ and any $ \max\menge{1,(1+\sigma)^{-1}} < p < \infty $, if $ 2 \, l - |\alpha| \leq 0 $. We can put this together with the on-diagonal estimates if we set $ \varepsilon_2 \ldef 1 $. \\

	Turning to the rescalation of these estimate, we let now $ I = (0,S) $ with $ S > 0 $ arbitrary and assume that $ u $ is a $ \sigma $-solution to $ f $ on $ I \times \Hquer $ with zero initial value. Then $ u $ is also a $ \sigma $-solution to $ f $ on $ (0,s) \times \Hquer $ for any $ s \in (0,S] $ with vanishing initial value. We denote the scaling function by $ A_\lambda: (\hat{s},\hat{y}) \mapsto (\lambda \, \hat{s}, \lambda \, \hat{y}) = (s,y) $ as above and start with the pointwise estimates. For $ (l,|\alpha|) = (0,1) $ as well as $ (l,|\alpha|) = (\frac{1}{2},2) $, and for any $ (s,y) \in I \times \Hquer $ we get 
	\begin{align*}
		y_n^l \, |\dy^\alpha u(s,y)| = \sqrt{\lambda}^{- |\alpha| - 1} \, \left(\frac{y_n}{\lambda} \right)^l \, \Bigl|\dyhut^\alpha (u \circ A_\lambda)(\lambda^{-1} \, s, \lambda^{-1} \, y)\Bigr|. 
	\end{align*}
	Now fix an $ s \in I $ and choose $ \lambda \ldef s $. Since $ u \circ A_\lambda $ is a $ \sigma $-solution to $ \lambda \, (f \circ A_\lambda) $ on $ (0,1) \times \Hquer $ with zero initial value, we can apply the above estimate to get 
	\begin{align*}
		y_n^l \, |\dy^\alpha u(s,y)| \lesssim \sqrt{s}^{1 - |\alpha|} \left(\sqrt{s}^{2 - \theta} \Yon{f} +  \sqrt{s}^{\varepsilon_1 - 1} \Yoffone{f} \right) 		
	\end{align*} 
	for any $ \max\menge{2 \, (n + 1),(1+\sigma)^{-1}} <  p < \infty $.
	We then set $ \varepsilon_1\ldef - \theta + 3  $, hereby restricting ourselves to $ \theta < 3 $, and get \[ |\grady u(s,y)| + \sqrt{s} \, \sqrt{y_n} |D^2_y u(s,y)| \lesssim \sqrt{s}^{2 - \theta} \left(\Yon{f} + \Yoffoneth{f} \right). \]
. \\

	Finally, we rescale the $ \Lp $-estimate by noting that $ u $ is also a $ \sigma $-solution to $ f $ on $ (0,r^2) \times \Hquer $ with zero initial value for any $ 0 < r^2 \leq S $. We then have \[ \norm[L^{p}(Q_r(y_0))]{\Pnl \, \ds^k \dy^\alpha u} = \lambda^{\frac{n+1}{p} + l - k - |\alpha|} \, \norm[L^{p}(A_\lambda^{-1}(Q_r(y_0)))]{\Pnlhut \, \dshut^k \dyhut^\alpha (u \circ A_\lambda)}. \] We can see that \[ A_\lambda^{-1}(Q_r(y_0)) \subset \vereinig[N]{i = 1} Q_{\frac{r}{\sqrt{\lambda}}}(\lambda^{-1} y_i) \text{ for } \lambda^{-1} y_i \in B_{4 \constd^2 \frac{r}{\sqrt{\lambda}}}(\lambda^{-1} y_0) \] and a number $ N $ only depending on the dimension $ n $. With $ \lambda = r^2 $ and the invariance of the scaling that makes $ u \circ A_{r^2} $ a $ \sigma $-solution to $ r^2 \, (f \circ A_{r^2}) $ on $ (0,1) \times \Hquer $, we can apply the statement proven above in every summand to get 
	\begin{align*}
		\norm[L^{p}(Q_r(y_0))]{\Pnl \, \ds^k \dy^\alpha u} \lesssim \, r^{\frac{2n+2}{p} + 2l - 2k - 2|\alpha| + 2} \, \summe[N]{i=1} &\left( 1 + \sqrt{\frac{y_{i,n}}{r^2}} \right)^{2l - |\alpha| + 1}  |Q_1(r^{-2} y_i)|^\frac{1}{p} \\
		&\bigl(\Yon{f \circ A_{r^2}} + \Yoffone{f \circ A_{r^2}} \bigr),
	\end{align*}
	where the last two expressions are understood with respect to the time interval $ (0,1) $ and the gradient in $ \Yon{\, \cdot \,} $ is taken with respect to $ \hat{y} $. Some calculations combined with the choice $ \varepsilon_1\ldef 3 - \theta $ lead to
	\begin{align*} 
		\supre[0 < r^2 \leq S \atop z \in H] r^{2k + |\alpha| - 3 + \theta} \, (r + \sqrt{z_n})^{-2l + |\alpha| - 1} \, |Q_r(z)|^{-\frac{1}{p}} \, &\norm[L^{p}(Q_r(z))]{\Pnl \, \ds^k \dy^\alpha u} \\
		&\lesssim  \Yon{f} + \Yoffoneth{f}  
	\end{align*}
	for any $ \theta < 3 $ and any $ \max\menge{1,(1+\sigma)^{-1}} < p < \infty $, if $ l - k - |\alpha| + 1 = - 1 $ and $ 2 \, l - |\alpha| + 1 \leq 0 $. These are exactly the exponents in the definition of $ X(p) $ if we set $ \theta \ldef 2 $, thus reaching the statement. 
\end{bew}
 
\section{Non-Linear Estimates} \label{section_nonlinear}
We come to the proof of Theorem \ref{theo_PE} that we state in a more detailed fashion here. An idea developed in \citet{angenent_semiflow} that was pushed further also by \citet{kochlamm} yields the analyticity in the temporal and tangential directions of the solution we construct. 

\begin{theo} \label{theo_ExInX} 
	Let $ \maxi\menge{2 \, (n + 1),(1+\sigma)^{-1}} <  p < \infty $. Then there exists an $ \varepsilon > 0 $, a $ \delta = \delta(n,\sigma) > 1 $ and a constant $ c_1 = c_1(n,\sigma,p) > 0 $ such that for any $ u_0 \in \homLip(\Hquer) $ with \[ \norm[\Linfty(H)]{\grady u_0} \leq \varepsilon \] we can find a $ \sigma $-solution $ u_* \in X(p) $ of the PE on $ (0,S) \times \Hquer $ with initial value $ u_0 $ satisfying \[ \norm[X(p)]{u_*} \leq c_1 \, \norm[\Linfty(H)]{\grady u_0} \] that is unique within \[ B_{\delta \varepsilon}^{X(p)} \ldef \menge{u \in X(p) \mid \norm[X(p)]{u} < \delta \, \varepsilon}. \] 
	Moreover, $ u_* $ depends analytically on the data $ u_0 $, we have $ u_* \in \glatt((0,S) \times \Hquer) $ and \[ \supre[(s,y) \in (0,S) \times H] s^{k+|\alpha|} \left|\ds^k \dy^\alpha \grady u_*(s,y) \right| \leq c_2 \, \norm[\Linfty(H)]{\grady u_0}  \] for any $ k \in \N_0 $ and $ \alpha \in \N_0^n $ with a constant $ c_2 = c_2(n,\sigma,k,\alpha) $. \\
	Furthermore, $ u_* $ is analytic in the temporal and tangential directions on $ (0,S) \times \Hquer $ with a $ \Lambda > 0 $ and a $ C = C(n) > 0 $ such that \[ \supre[(s,y) \in (0,S) \times H] s^{k + |\alpha'|} \, |\ds^k \dystrich^{\alpha'} \grady u_*(s,y)| \leq C \, \Lambda^{-k-|\alpha'|} \, k! \, \alpha'! \, \norm[\Linfty(H)]{\grady u_0} \] for any $ k \in \N_0 $ and $ \alpha' \in \N_0^{n-1} $ with $ k + |\alpha'| > 0 $.
\end{theo}
\begin{bew}
	Define the function \[ F: \homLip(\Hquer) \times X(p) \to X(p) \] by assigning $ (u_0,u) \in \homLip(\Hquer) \times X(p) $ to the $ \sigma $-solution to $ f[u] $ on $ (0,S) \times \Hquer $ with initial value $ u_0 $. By Theorem \ref{theo_linear} and Lemma \ref{lemma_YggX} we can find a $ R_1 > 0 $ and an $ \varepsilon_1 > 0 $ such that for any $ u_0 \in \overline{B}_{\varepsilon_1}^{\homLip(\Hquer)} $ the map $ F(u_0,\, \cdot\, ) $ is a contraction within $ \overline{B}_{R_1}^{X(p)} $. Hence the contraction mapping principle provides us with a unique fixed point $ u_* \in \overline{B}_{R_1}^{X(p)} $ that depends on $ u_0 $ in a Lipschitz continuous way. Thus $ u_* $ is the unique global $ \sigma $-solution for the non-linear equation with initial value $ u_0 $. The bound on $ \norm[X(p)]{u_*} $ especially implies that $ u_* $ is Lipschitz in time and space. By \citet{koch_habil} the smoothness of $ u_* $ on $ (0,S) \times \Hquer $ follows. \\
	But moreover, $ F $ and therefore also \[ G: \homLip(\Hquer) \times X(p) \to X(p), \,  G(u_0,u) \ldef u - F(u_0,u) \] are analytic on $ \homLip(\Hquer) \times \overline{B}_{R_1}^{X(p)} $. Since $ G(0,0) = 0 $ and $ D_u G(0,0) = \Id $, the analytic implicit function theorem on Banach spaces is applicable \citep{deimling}) and yields the existence of balls \[ \overline{B}_{\varepsilon_2}^{\homLip(\Hquer)} \subset \homLip(\Hquer) \text{ and } \overline{B}_{R_2}^{X(p)} \subset \overline{B}_{R_1}^{X(p)} \] alongside with an analytic function \[ A: B_{\varepsilon_2}^{\homLip(\Hquer)} \to B_{R_2}^{X(p)} \] that satisfies $ A(0) = 0 $ and $ F(u_0,u) = u $ for any $ u_0 \in B_{\varepsilon_2}^{\homLip(\Hquer)} $ and $ u \in B_{R_2}^{X(p)} $ if and only if $ u = A(u_0) $. By the uniqueness of the fixed point in $ \overline{B}_{R_1}^{X(p)} \supset B_{R_2}^{X(p)} $ above we conclude that $ A $ sends $ u_0 \in B_{\varepsilon_3}^{\homLip(\Hquer)} $ to $ u_* $ analytically, where $ \varepsilon_3 \ldef \min\menge{\varepsilon_1,\varepsilon_2} $. \\
	We now consider $ (\tau,\xi') \in B_{\kappa_1}^\R(1) \times B_{\delta_1}^{\Rnm} $ and define \[ f_{\tau,\xi'}[u] \ldef \tau \, f[u] - (1 - \tau) \, L_\sigma u - \xi' \cdot \grady' u. \] Note that then $ f_{1,0}[u] = f[u] $. Similar as above we define \[ \widetilde{F}: B_{\kappa_1}^\R(1) \times B_{\delta_1}^{\Rnm} \times \homLip(\Hquer) \times B_{R_2}^{X(p)} \to X(p) \] that maps $ (\tau,\xi',u_0,u) $ to the $ \sigma $-solution to $ f_{\tau,\xi'}[u] $ on $ (0,S) \times \Hquer $ with initial value $ u_0 $. 
	A straightforward calculation as before shows that we have \[ \norm[X(p)]{\widetilde{F}(\tau,\xi',u_0,u)} \lesssim \norm[\Linfty(H)]{\grady u_0} + \tau \, \norm[X(p)]{u}^2 + (|1 - \tau| + |\xi'|) \, \norm[X(p)]{u}. \] Hence we can conclude that $ D_u \widetilde{F}|_{(1,0,0,0)} $ vanishes. Applying the analytic implicit function theorem once more, this time on \[ \widetilde{G} \ldef \Id_{X(p)} - \widetilde{F} \] at the point $ (\tau,\xi',u_0,u) = (1,0,0,0) $ we then obtain the existence of $ \kappa_2 < \kappa_1 $, $ \delta_2 < \delta_1 $, $ R_3 < R_2 $ and $ \varepsilon_4 $ as well as a unique analytic function \[ \widetilde{A}: B_{\kappa_2}^\R(1) \times B_{\delta_2}^{\Rnm} \times B_{\varepsilon_4}^{\homLip(\Hquer)} \to B_{R_3}^{X(p)} \] such that \[ \widetilde{G}(\tau,\xi',u_0,\widetilde{A}(\tau,\xi',u_0)) = 0. \] 
	For $ (\tau,\xi') \in B_{\kappa_2}^\R(1) \times B_{\delta_2}^{\Rnm} $ and $ (s,y) \in (0,S) \times \Hquer $ let us now consider the transformation \[ U: (\tau,\xi',s,y) \mapsto (\tau s, y' - \xi' s, y_n). \] A simple calculation shows that $ u \circ U(\tau,\xi',\, \cdot\, ,\, \cdot\, ) $ is a $ \sigma $-solution to $ f_{\tau,\xi'}[u] $ on $ (0, \tau S) \times \Hquer $ with initial value $ u_0 $ if $ u $ is a $ \sigma $-solution to $ f[u] $ on $ (0,S) \times \Hquer $ with initial value $ u_0 $. Therefore we have that \[ \widetilde{G}(\tau,\xi', A(u_0) \circ U(\tau,\xi',\, \cdot\, ,\, \cdot\, )) = 0. \] It is also clear that \[ A(u_0) \circ U (\tau,\xi',0,\, \cdot\, ) = u_0 = \widetilde{A}(\tau,\xi',u_0)(0,\, \cdot\, ) \] and thus the above uniqueness results for $ \varepsilon \ldef \min\menge{\varepsilon_3,\varepsilon_4} $ and $ R_3 \rdef \delta \,m \varepsilon $ imply that $ A(\, \cdot\, ) \circ U = \widetilde{A} $. Especially, $ u_* \circ U(\, \cdot\, ,\, \cdot\, ,s,y) $ is analytic as a function of $ \tau $ and $ \xi' $ into $ X(p) $ near $ \tau = 1 $ and $ \xi' = 0 $ for any $ (s,y) \in (0,S) \times \Hquer $. 
	But we have
	\begin{align*}
		\bigl|\partial_\tau \grady (u_* \circ U(\tau,\xi',s,y))|_{(\tau,\xi') = (1,0)} \bigr| &= \bigl|s \, \ds \grady u_*(s,y)\bigr|, \\
		\bigl|\partial_{\xi_j} \grady (u_* \circ U(\tau,\xi',s,y))|_{(\tau,\xi') = (1,0)}\bigr| &= \bigl|- s  \, \dyj \grady u_*(s,y)\bigr| \text{ for } j = 1, \ldots, n-1
	\end{align*}
	with similar formulas for higher order and mixed derivatives. This implies the analyticity of $ u_* $ on $ I \times \Rn $ in $ s $ and $ y' $ with the formula \[ \supre[(s,y) \in (0,S) \times H] s^{k + |\alpha'|} \, |\ds^k \dystrich^{\alpha'} \grady u_*(s,y)| \lesssim \Lambda^{-k-|\alpha'|} \, k! \, \alpha'! \, \norm[\Linfty(H)]{\grady u_0} \] for any $ k \in \N_0 $ and $ \alpha' \in \N_0^{n-1} $, 
	where $ \Lambda > 0 $ and the constant depends only on $ n $. \\
		For the missing derivative in the tangential direction define the components of $ \vec{F}[u_*] $ by \[ F_j[u_*] \ldef \dyj u_* \text{ for } j = 1, \ldots, n - 1 \] and \[ F_n[u_*]\ldef - \frac{|\grady' u_*|^2 - \dyn u_*}{1 + \dyn u_*}, \] and rephrase the PE as \[ 0 = \ds u_* - L_\sigma u_* - f[u_*] = \ds u_* - \Pn^{-\sigma} \, \divy (\Pn^{1 + \sigma} \vec{F}[u_*]). \] We differentiate the equation with respect to the vertical direction and obtain \[ \dyn(\Pn^{2 + \sigma} \dyn F_n[u_*]) = \Pn^{1 + \sigma} \, (\ds \dyn u_* - \lapy' u_* - \Pn \lapy' \dyn u_*). \] An integration in $ y_n $ then shows that \[ |\dyn F_n[u_*](s,y)| \lesssim \norm[\Linfty(H)]{\ds \dyn u_*(s) - \lapy' u_*(s) - \Pn \, \lapy' \dyn u_*(s)} \] for $ (s,y) \in (0,S) \times H $ with a constant depending only on $ \sigma $. By the analyticity estimate, the first terms of the latter are bounded by $ s^{-1} \, \norm[\Linfty(H)]{\grady u_0} $ times an absolute constant. The third term can be handled by the pointwise term in the $ X^{(2)} $-norm. From there we get in a similar way as for the analyticity estimate that \[ |\ds^k \dystrich^{\alpha'} D^2_y u(s,y)| \lesssim s^{-k - |\alpha'| - \frac{1}{2}} \, y_n^{-\frac{1}{2}} \, \norm[\Linfty(H)]{\grady u_0} \] with a constant depending on $ n $, $ k $ and $ \alpha' $. For $ \sqrt{y_n} \leq \sqrt{s} $ this implies \[ y_n |\lapy' \dyn u_*(s)| \lesssim y_n^\frac{1}{2} \, s^{-\frac{3}{2}} \, \norm[\Linfty(H)]{\grady u_0} \leq s^{-1} \, \norm[\Linfty(H)]{\grady u_0}. \] But for $ \sqrt{y_n} > \sqrt{s} $ we are in the parabolic regime and get the same estimate directly.  However, we also have \[ |\dyn F_n[u_*]| \gtrsim |\dyn^2 u_*|. \] To finish the proof, this argument can be iterated, introducing a dependence on the order of derivatives into the constant.
\end{bew}

\addsubsection*{Acknowledgments}
This work is based on my PhD thesis. It would not have been possible without my advisor Herbert Koch. He introduced me to the different theories that fit together so beautifully in this field and showed an extraordinary degree of patience, kindness and intuition in guiding me. I also benefited from helpful conversations with my friend Dominik John. I would further like to express my gratitude towards the Bonn International Graduate School in Mathematics (BIGS) and Cusanuswerk for their support.




\ifthenelse{ 
		 \equal{\klassenart}{buch}
		}
		{ 
		 \backmatter 
		}
		{ 
		}



\ifthenelse{ 
		 \equal{\klassentyp}{standard}
		}
		{ 
		\ifthenelse{ 
				 \equal{\klassenart}{artikel}
				}
				{ 
				 \addcontentsline{toc}{section}{\refname} 
				}
				{ 
				 \addcontentsline{toc}{section}{\bibname} 
				}
		}
		{ 
		}

\bibliographystyle{\bibstil\zitierstil\sprache}

\bibliography{\journalnames,\authornames,\database}

\end{document}